\documentclass[pdflatex,sn-basic]{sn-jnl}
\usepackage{graphicx}%
\usepackage{multirow}%
\usepackage{amsmath,amssymb,amsfonts}%
\usepackage{amsthm}%
\usepackage[title]{appendix}%
\usepackage{xcolor}%
\usepackage{textcomp}%
\usepackage{booktabs}%
\usepackage{enumitem} 

\usepackage{fix-cm}
\usepackage[T1]{fontenc}
\usepackage[utf8]{inputenc}
\usepackage{bookmark}

\usepackage{subcaption}

\usepackage[normalem]{ulem}

\makeatletter
\renewcommand*{\top}{%
  {\mathpalette\@transpose{}}%
}
\newcommand*{\@transpose}[2]{%
  \raisebox{\depth}{$\m@th#1\scriptscriptstyle\mathsf{T}$}%
}
\makeatother

\newcommand{\zmespthing}{z_{\mbox{\normalfont\protect\tiny MESP}}}

\newcommand{\zspectralthing}{z_{\mbox{\normalfont\protect\tiny $\mathcal{S}$}}}

\newcommand{\zgfactplusthing}{z_{\mbox{\protect\tiny GFact$^{+}$}}}

\newcommand{\zdgfactplusthing}{z_{\mbox{\protect\tiny DGFact$^{+}$}}}

\newcommand{\zddgfactplusthing}{z_{\mbox{\protect\tiny DDGFact$^{+}$}}}

\newcommand{\ztracething}{z_{\mbox{\protect\tiny Trace$^{+}$}}}

\DeclareMathOperator{\rank}{rank}

\DeclareMathOperator{\diag}{\textstyle{diag}}
\DeclareMathOperator{\Diag}{\textstyle{Diag}}

\DeclareMathOperator{\tr}{\textstyle{tr}}

\DeclareMathOperator{\dom}{\textstyle{dom}}
\DeclareMathOperator{\ldet}{ldet}

\setcitestyle{authoryear,open={[},close={]},aysep={,}}

\theoremstyle{thmstyleone}%

\theoremstyle{thmstyletwo}%

\theoremstyle{thmstylethree}%
\newtheorem{theorem}{Theorem}%
\newtheorem{proposition}[theorem]{Proposition}%
\newtheorem{lemma}[theorem]{Lemma}
\newtheorem{corollary}[theorem]{Corollary}
\newtheorem{definition}[theorem]{Definition}%
\newtheorem{example}[theorem]{Example}%
\newtheorem{remark}[theorem]{Remark}%

\begin{document}

\title{Augmenting the generalized \break factorization bound for GMESP$^{\star}$}

\author[1]{\fnm{Ananias S.} \sur{Machado}}\email{ananiasmachado@poli.ufrj.br}

\author[1]{\fnm{Marcia} \sur{Fampa}}\email{fampa@cos.ufrj.br}

\author[2]{\fnm{Jon} \sur{Lee}}\email{jonxlee@umich.edu}

\author[2]{\fnm{Gabriel} \sur{Ponte}}\email{gabponte@umich.edu}

\affil[1]{\orgname{\small Universidade Federal do Rio de Janeiro},
\country{Brasil}}

\affil[2]{\orgname{\small University of Michigan},
\city{Ann Arbor}, 
\state{Michigan}, 
\country{USA}}

\abstract{Li proposed an augmented factorization upper bound for use in branch-and-bound (B\&B) algorithms
aimed at 
the maximum-entropy sampling problem (MESP). This bound improves on the factorization bound for MESP when the input covariance matrix is positive definite. GMESP, which selects $s$ of $n$ variables so as to maximize the product of the $t$ greatest eigenvalues of the corresponding principal submatrix of the covariance matrix, is a 
generalization of MESP (the case $t=s$), and the factorization bound for MESP generalizes to GMESP. We extend the augmentation technique of Li to improve on the generalized factorization bound for GMESP when the input covariance matrix is positive definite. We derive fundamental properties of the augmented bound,
and we develop primal--dual constructions 
for use within B\&B. We give checkable conditions under which the relaxation can be solved in closed form, which always hold when $t=1$, at the root and at every subproblem, and a condition under which the closed-form solution is optimal for GMESP. We compare the augmented bound with the spectral bound, establishing an exact comparison for $t=1$, an improved worst-case comparison for general $t$, and sufficient conditions for the dominance of each bound over the other. Computational experiments on benchmark instances show that augmentation can substantially reduce root gaps of the generalized factorization bound, while even reducing its computational time, and that it improves the performance of B\&B. Moreover, a hybrid 
strategy that uses the trace upper bound to prune nodes and the closed-form certificate to avoid solving the relaxation can significantly accelerate B\&B when $t$ is small.}


{\let\thefootnote\relax\footnotetext{$^{\star}$A very short preliminary version of this work is to appear in the issue of Lecture Notes in Operations Research associated to the conference CLAIO 2026.}}

\maketitle

\section{Introduction}\label{sec:intro}

Let $C$ be a symmetric positive-semidefinite matrix with rows and columns
indexed by $N:=\{1,2,\ldots,n\}$, with $n>1$.
We interpret $C$ as the covariance matrix of an  $n$-dimensional Gaussian random vector $\mathcal{Y}_N$\,.
Let $t$ and $s$ be integers satisfying
$0<t\le \rank(C)$ and $t\le s<n$. For a binary vector $x\in\{0,1\}^n$, let
$S(x)$ denote its support, let $C[S(x),S(x)]$ denote the corresponding principal
submatrix of $C$, let $\mathbf{e}$ denote the all-ones vector, let
$\lambda_\ell(X)$ denote the $\ell$-th greatest eigenvalue of a symmetric matrix
$X$, and let $\operatorname{ldet}(\cdot)$ denote the natural logarithm of the
determinant.

We consider the \emph{generalized maximum-entropy sampling problem}
\begin{equation}\tag{GMESP}\label{GMESP}
    \begin{array}{ll}
 z(C,s,t)\!:=&\!\max\!      \textstyle \left\{\sum_{\ell=1}^t \log( \lambda_\ell(C[S(x),S(x)])) \,:\, \mathbf{e}^\top x =s,~ x\in\{0,1\}^n\right\},
    \end{array}
\end{equation}
which seeks a subset of $s$ variables maximizing the geometric mean of the greatest
$t$ eigenvalues of the corresponding principal submatrix. \ref{GMESP} can be interpreted as a selection problem related to principal component analysis (PCA). The goal there is to choose the 
$s$-subvector $\mathcal{Y}_S$ (of $\mathcal{Y}_N$) having $|S|=s$,
to maximize the geometric mean of the variances of its $t$ leading principal components; see \citep*{GMESP_Alg}.

When $t=s$, \ref{GMESP} reduces to the classical and well-studied \emph{maximum-entropy sampling problem} 
\begin{equation*}\tag{MESP}\label{MESP}
    \hypertarget{zmesptarget}{\zmespthing}(C,s):= \max \left\{\ldet(C[S(x),S(x)]) ~:~\mathbf{e}^\top x =s,~ x\in\{0,1\}^n\right\},
\end{equation*}
whose objective is to maximize  Shannon's differential entropy of the selected  $s$-subvector $\mathcal{Y}_S$; 
see \citep*{SW} and \citep*{FLbook} (and the many references therein). 

Exact algorithms for \ref{GMESP} and its special cases are typically based on branch-and-bound (B\&B).
The efficiency of these algorithms depends critically on the availability of strong upper bounds for pruning nodes of the search tree.
Among the strongest known bounds are those obtained from convex relaxations of nonlinear formulations.
For \ref{MESP}, the factorization bound (DFact) has proved particularly effective, and its generalization to \ref{GMESP}, the generalized factorization bound (DGFact), has been shown to be effective within B\&B when $s-t$ is small
\citep*{GMESP_Alg}.
Besides providing high-quality upper bounds, these relaxations also yield valuable information for variable fixing and branching strategies.

Recently,
\citep*{li2025augmented}
introduced an augmented version of the DFact bound for \ref{MESP}, called DFact$^+${,}
which improves upon the DFact bound when $C$ is positive definite. We note that key benchmark data sets for \ref{MESP} have this property.
Li's analysis established several attractive theoretical properties of the augmented bound, including monotonicity with respect to the augmentation parameter and improved variable-fixing procedures.

The purpose of our present work is to extend Li's augmentation framework from \ref{MESP} to the more general \ref{GMESP}.
For this reason, we assume henceforth that $C$ is positive definite. 
We wish to emphasize that the extension of the augmentation technique from \ref{MESP} to \ref{GMESP} is not straightforward at all.
Several arguments in the analysis of DFact$^+$ rely fundamentally on the equality of $t$ and $s$, while \ref{GMESP} allows $t< s$.
Consequently, new theoretical tools are required to establish analogous properties for the generalized setting.

A very short preliminary version of this work, appears as \citep*{AugGFactShort}.

\medskip
\noindent {\bf Organization and contributions.} 

Our first contribution is the introduction of the \emph{augmented generalized factorization bound} (the DGFact$^+$ bound).
We derive its convex formulation and its Lagrangian dual, providing the basis for efficient computation, and we show that the bound is independent of the factorization of $C-\psi I$, so that it depends only on $C$ and on the augmentation parameter $\psi$. We show that the nonconvex relaxation is exact, in the sense that restricting its variables to be binary yields a correct formulation of GMESP, and we identify a spectral condition under which its convex bidual (the Lagrangian dual of its Lagrangian dual) is exact as well.
We then prove that the optimal value of the relaxation is monotonically nonincreasing with respect to the augmentation parameter, and we give a checkable sufficient condition for the decrease to be strict.
The proof for DFact$^+$ from \citep*{li2025augmented} relies on $t=s$ and does not extend directly; we obtain the generalization by combining a simple consequence of Li's majorization lemma with two properties of the function underlying the bound: it is componentwise nondecreasing{, which we derive from an explicit supergradient, and it is Schur-concave, which we derive from results of \citep*{Nikolov}}.

Our second contribution concerns the use of DGFact$^+$ within B\&B algorithms.
We extend the primal--dual construction of associated feasible solutions to the augmented setting with fixed variables, prove that the construction preserves optimality, and show that the vector $\diag(F\hat\Theta F^\top)$, computed from the constructed dual solution, is a supergradient of the objective of DDGFact$^+$ and gives its directional derivatives. 
The constructed dual solutions can also be used to fix variables, as we do in our B\&B experiments.
These results provide the theoretical foundation for efficiently incorporating DGFact$^+$ into B\&B algorithms.

Our third contribution establishes situations in which the relaxation admits closed-form solutions.
In particular, we show that the relaxation can be solved in closed form when $t=1$, both at the root node and at B\&B nodes with variables fixed to one.
We further derive sufficient conditions under which the same closed-form solution remains optimal for arbitrary values of $t$, and, under an additional multiplicity condition on the least eigenvalue of $C$, is optimal for GMESP itself.

Finally, we compare the DGFact$^+$ bound with the spectral bound.
We prove that the two bounds are incomparable in general, for every value of $t$, and we derive an exact comparison formula when $t=1$. For general $t$, we establish a  worst-case comparison that improves the corresponding result for the original DGFact bound, and identify sufficient conditions guaranteeing dominance of each bound over the other.
Our analysis demonstrates that augmentation improves not only the numerical strength of the factorization bound but also its theoretical properties. Although DGFact$^+$
 preserves the classical worst-case guarantee relative to the spectral bound, it admits a strictly better estimate whenever the augmentation parameter is positive and $s>t$, providing new insight into the effect of augmentation on determinant-maximization relaxations.

Our computational experiments on two benchmark covariance matrices show that augmentation substantially reduces the root gaps of the generalized factorization bound, especially for large $s$, while even reducing its computational time. Within B\&B, augmentation reduces both the root and the final gaps on all instances not solved within the time limit. Moreover, a hybrid strategy that uses the trace upper bound to prune nodes and the closed-form certificate to avoid solving the relaxation can significantly accelerate B\&B when $t$ is small.

The remainder of the paper is organized as follows.
Section~\ref{sec:gfactplus_convexity} introduces the augmented generalized factorization bound and derives its Lagrangian dual.
Section~\ref{sec:monotonicity_ddgfactplus} establishes its monotonicity with respect to the augmentation parameter.
Section~\ref{subsec:feasible_dual_sol_dgfactplus}
constructs primal--dual pairs of feasible solutions.
Section~\ref{sec:special_case_t=1}
derives closed-form solutions and optimality certificates.
Section~\ref{sec:comparison_spectral}
compares DGFact$^+$ with the spectral bound.
 Section~\ref{sec:computational_results}
reports computational experiments illustrating the practical performance of the proposed relaxation{, Section~\ref{sec:outlook} outlines future work, and Appendix~\ref{sec:A1} contains postponed proofs}.

\medskip
\noindent {\bf Literature.} 

\ref{GMESP} was introduced by \citep*{WilliamsPhD}; see
\citep*{LeeLind2020}. The latter developed a spectral bound for \ref{GMESP} and for its constrained version (which is not a straightforward extension, because the spectral bound is not based on a convex relaxation of \ref{GMESP}).
\citep*{SEA_proceedings} and 
\citep*{GMESP_Alg} developed a DGFact bound for \ref{GMESP}, extending the pre-existing DFact bound for \ref{MESP}; see \citep*{FLbook}. 
As we have mentioned, \citep*{li2025augmented}
improved the DFact bound for \ref{MESP}, when $C$ is positive definite. \citep*{ExtendedGMESP}  extended other bounds for \ref{MESP} to \ref{GMESP} (specifically the ``linx bound'', the ``NLP-Id bound'', and the ``complementary NLP-Id bound''), and compared such bounds and the spectral bound, both theoretically and computationally; the scaled glinx bound from that work serves as a benchmark in Section~\ref{sec:computational_results}.

Some further background and a wealth of references on
\ref{MESP} can be found in \citep*{FLbook} and
\citep*{FL_update}. 

A less obvious but very important special case of \ref{GMESP} is the 0/1 D-optimality problem, which arises by taking a full column rank ``design matrix'' $A\in\mathbb{R}^{n\times t}$, and letting $C:=AA^\top$; see \citep*{PonteFampaLeeMPB} and the references therein, for example. It is important to note that normally $n>t$, and then $C$ is not positive definite. So the results that we present do not apply to the 0/1 D-optimality problem. In any case, a study of mappings between \ref{MESP} instances and 0/1 D-optimality instances, and the implications for upper bounds on instances of these problems 
can be found in \citep*{MESP2DOPT}, which effectively presents a different bridge between these two problems than is afforded by \ref{GMESP}.

\medskip
\noindent {\bf Further notation.}
We let $\mathbb{S}^n$  (resp., $\mathbb{S}^n_+$, $\mathbb{S}^n_{++}$)
 denote the set of symmetric (resp., positive-semidefinite, positive-definite) matrices of order $n$. We let $\Diag(x)$ denote the $n\times n$ diagonal matrix with diagonal elements given by the components of $x\in \mathbb{R}^n$, and $\diag(X)$ denote the $n$-dimensional vector with elements given by the diagonal elements of $X\in\mathbb{R}^{n\times n}$.
We denote the identity matrix of order $n$ by $I_n$\,.  For a positive integer $q$, we denote by $\mathbb{I}_q$ a vector with the first $q$ elements equal to one and the others equal to zero (its dimension should be clear from the context).
For matrices $A$ and $B$ with compatible shapes,
$A\bullet B:=\tr(A^\top B)$ is the matrix dot-product, where $\tr$ denotes the trace. 
For matrix $X$,  $X[S,T]$\, is the submatrix  with row (column) indices $S$ ($T$);
 if $S$ ($T$) is the set of all row (column) indices of $X$, we write $X[\cdot, T]$ ($X[S,\cdot]$).
  We let  $\lambda(X):=(\lambda_{\max}(X)=\lambda_1(X),\lambda_2(X),\ldots,\allowbreak \lambda_n(X)=\lambda_{\min}(X))^\top$.  
 For $a\in\mathbb{R}$, $a_+:=\max\{0,a\}$. 

\section{Augmenting the generalized factorization bound}
\label{sec:gfactplus_convexity}

In this section we extend the augmented factorization framework of
\citep*{li2025augmented}, originally developed for \ref{MESP}, to \ref{GMESP}. This leads to a new upper bound for \ref{GMESP},  called the
augmented generalized factorization bound, which
strengthens the generalized factorization bound of
\citep*{GMESP_Alg} whenever the covariance matrix is
positive definite.

Our starting point is a nonconvex continuous relaxation of \ref{GMESP}. We then
derive its Lagrangian dual, denoted by DGFact$^+$, which forms the basis for the
structural properties established in the subsequent sections, including independence of the chosen factorization of  $C-\psi I$ and 
monotonicity with respect to the augmentation parameter.

Let
$0\le \psi\le\lambda_{\min}(C)$,
and consider an arbitrary  factorization
$C-\psi I = F(\psi) F(\psi)^\top$,
where
$F(\psi)\in\mathbb{R}^{n\times k}$,
with
{$\max\{t,\rank(C-\psi I)\}\le k\le n$. The requirement $k\ge t$ ensures that the functions defined below, which involve the $t$ greatest eigenvalues of $k\times k$ matrices, are well defined; it can always be met (e.g., by appending zero columns to $F(\psi)$), because $t\le s<n$.}
Whenever the dependence on $\psi$ is clear from the context, we simply write
$F:=F(\psi)$\,.
For $0<t\leq s$ and  $X \in \mathbb{S}^k_{+}$\,, let
\begin{equation*}
    f_t(X; \psi) := \textstyle \sum_{i=1}^{t}\log(\lambda_i(X) + \psi),
\end{equation*}
 and for $x \in \mathbb{R}^n_{+}$\,, let 
\begin{equation*}
    M_{\psi}(x) := F^\top \Diag(x) F.
\end{equation*}
Now, we define the nonconvex relaxation 
\begin{equation*}
\zgfactplusthing(C,s,t;\psi,F):=\max \left\{f_t(M_{\psi}(x);\psi) ~:~ \mathbf{e}^{\top}x=s,~x\in[0,1]^n \right\}.\tag{GFact$^{+}$}\label{GFactplus}
\end{equation*}

When $\psi=0$, formulation~\ref{GFactplus} reduces to the generalized
factorization relaxation GFact introduced in \citep*{GMESP_Alg}.

In what follows, a superscript ``$+$'' denotes the augmented version of a
formulation. The corresponding formulation without the superscript is recovered
by setting $\psi=0$ and coincides with the formulation introduced in
\citep*{GMESP_Alg}.

Throughout, we say that a continuous relaxation of \ref{GMESP} is \emph{exact}
when imposing $x\in\{0,1\}^n$ on it yields a correct formulation of
\ref{GMESP}, that is, when its objective function agrees with that of
\ref{GMESP} at every binary feasible solution. 

The following theorem extends \citep*[Thm.~1]{GMESP_Alg} by showing that
\ref{GFactplus} is an exact relaxation of \ref{GMESP}. 

\begin{theorem}\label{Gmesp2GFact} For every  feasible solution $\hat{x}$ of \ref{GMESP}, 
\[
\textstyle
\sum_{\ell=1}^t \log( \lambda_\ell(C[S(\hat x),S(\hat x)])) = f_t(M_{\psi}(\hat x);\psi).
\]
\end{theorem}

\begin{proof}
    Let $S := S(\hat x)$. By the definition of $M_\psi$\,, we obtain
    \begin{equation*}
    \textstyle
        M_{\psi}(\hat x) = \sum_{i \in S} F[i,\cdot]^\top F[i,\cdot] = F[S,\cdot]^\top F[S,\cdot].
    \end{equation*}
As $C[S,S] - \psi I = F[S,\cdot] F[S,\cdot]^\top$, and the positive-semidefinite matrices $F[S,\cdot]^\top F[S,\cdot]$ and $F[S,\cdot] F[S,\cdot]^\top$ have the same nonzero eigenvalues, hence the same $\ell$-th greatest eigenvalue for $\ell\le\min\{s,k\}$, and as $t\le s$ and $t\le k$, for every $\ell = 1,\dots,t$, we have
    \begin{equation*}
        \begin{aligned}
            \lambda_\ell(C[S,S]) & = \lambda_\ell(F[S,\cdot] F[S,\cdot]^\top) + \psi  = \lambda_\ell(M_{\psi}(\hat x)) + \psi.
        \end{aligned}
    \end{equation*}
    The result follows. 
\end{proof}

\subsection{Duality for  \texorpdfstring{\ref{GFactplus}}{GFactplus}} \label{subsec:dual_gfactplus}
Although \ref{GFactplus} is an exact relaxation of \ref{GMESP}, its nonconvexity prevents its direct use for computing upper bounds in a B\&B framework. Following the approach of \citep*{GMESP_Alg}, we derive the Lagrangian dual of \ref{GFactplus}, yielding a convex relaxation of \ref{GMESP} that can be solved efficiently. While our derivation follows the same general strategy as that of \citep*{GMESP_Alg}, the extension is not straightforward due to the additional augmentation term. The main challenge lies in obtaining a closed-form characterization of the dual function as a function only of the dual variables. To address this, we first establish the auxiliary results in Lemmas~\ref{lem:von_neumann} and \ref{lem:subproblem_W_multivariable}, which provide the key ingredients for the dual derivation. Building on these results, Theorem~\ref{thm:subproblem_W} yields a concise characterization of the dual subproblem, leading to a derivation that is shorter and more self-contained than the corresponding analysis in \citep*{GMESP_Alg}. 

Introducing an auxiliary matrix variable $W\in\mathbb{S}^k$, \ref{GFactplus} can be equivalently written~as
\begin{align*}\textstyle
    \max\left\{\sum_{\ell=1}^t \log (\lambda_\ell(W)+\psi) : W = M_{\psi}(x), \mathbf{e}^\top x = s, x \in [0,1]^n\right\}.
\end{align*}
The corresponding Lagrangian is  
\begin{align*}
    &\mathcal{L}(x, W, \Theta, \upsilon, \nu, \tau) := \\&\qquad\textstyle \sum_{i=1}^{t}\log(\lambda_i(W) + \psi)+ \Theta\bullet (M_{\psi}(x) - W) 
    + \upsilon^\top x + \nu^\top(\mathbf{e}-x) + \tau(s-\mathbf{e}^\top x), 
\end{align*}
with domain $\dom \mathcal{L} = \mathbb{R}^n \times \mathbb{S}^k_+ \times \mathbb{S}^k \times \mathbb{R}^n \times \mathbb{R}^n \times \mathbb{R}$. 
The associated dual function is
\begin{equation*}
    \mathcal{L}^*(\Theta, \upsilon, \nu, \tau) := \textstyle\sup_{x \in \mathbb{R}^n, W \in \mathbb{S}^k_+}\mathcal{L}(x, W, \Theta, \upsilon, \nu, \tau),
\end{equation*}
and the corresponding Lagrangian dual problem is
\begin{align*}
    \zdgfactplusthing
    (C,s,t;\psi,F) := \inf\{\mathcal{L}^*(\Theta, \upsilon, \nu, \tau) \;:\; \upsilon \geq 0,\; \nu \geq 0\}.
\end{align*}

We refer to $\zdgfactplusthing
    (C,s,t;\psi,F)$ as the \emph{augmented generalized factorization bound}.
As the Lagrangian is separable with respect to the primal variables $x$ and $W$, the dual function decomposes into
$\mathcal{L}^*(\Theta, \upsilon, \nu, \tau) 
= \textstyle\mathcal{L}_1^*(\Theta) +  \mathcal{L}_2^*(\Theta, \upsilon, \nu, \tau) $, with
\begin{align}
    &\mathcal{L}_1^*(\Theta) = \textstyle\sup_{W \in \mathbb{S}^k_+} \left\{\textstyle\sum_{i=1}^{t}\log(\lambda_i(W) + \psi) - \Theta \bullet W \right\},\label{prob:subprob_W}\\
    &\mathcal{L}_2^*(\Theta, \upsilon, \nu, \tau) = \textstyle\sup_{x \in \mathbb{R}^n}\left\{\Theta \bullet M_{\psi}(x) + \upsilon^\top x + \nu^\top(\mathbf{e}-x) + \tau(s-\mathbf{e}^\top x)\right\}.\nonumber
\end{align}

The following auxiliary results lead to a closed-form expression for $ \mathcal{L}_1^*(\Theta)$ established in Theorem~\ref{thm:subproblem_W}.

\begin{lemma}~\!\!\![Von Neumann's trace inequality; see {{\citep*{Theobald1975}}}] (proof in Appendix)
\label{lem:von_neumann}
    Let $A, B \in \mathbb{S}^n$. Then
    \begin{equation}\label{eq:von_neumann}
    \textstyle
    \sum_{i=1}^n \lambda_i(A) \lambda_{n-i+1}(B) \leq A\bullet B
    \le \sum_{i=1}^n \lambda_i(A) \lambda_i(B).
    \end{equation}
    Moreover, both bounds are attained: if $Q=[q_1,\ldots,q_n]$ is an orthogonal
    matrix with $A=Q\Diag(\lambda(A))Q^\top$, and $\beta\in\mathbb{R}^n$ is any
    vector whose components are a rearrangement of those of $\lambda(B)$, then
    $\bar B:=Q\Diag(\beta)Q^\top$ satisfies $\lambda(\bar B)=\lambda(B)$ and
    \begin{equation*}
        \textstyle A\bullet \bar B=\sum_{i=1}^{n}\lambda_i(A)\,\beta_i\,.
    \end{equation*}
    In particular, taking $\beta_i:=\lambda_{n-i+1}(B)$ attains the left-hand
    bound in \eqref{eq:von_neumann}, and taking $\beta_i:=\lambda_i(B)$ attains
    the right-hand bound.
\end{lemma}

\begin{lemma}\label{lem:subproblem_W_multivariable}~(proof in Appendix)
    Let $\psi\ge0$ and $w\in\mathbb{R}^t$ with $w_1\leq w_2\leq \cdots\leq w_t$\,.
    Consider the optimization problem
\begin{equation}\label{eq:solve_y_ordered_given_z_ordered}
        \textstyle\sup_{y_1 \geq y_2\geq \dots \geq y_t \geq 0 } \left\{ \textstyle \sum_{i=1}^{t}\log(y_i + \psi) - w^\top y \right\}.
    \end{equation}
If $w_1 \leq 0$,  \eqref{eq:solve_y_ordered_given_z_ordered}  is unbounded above.
    Otherwise, the unique maximizer is
    \begin{equation*}
        y_i^* = \max\left\{0, 1/w_i - \psi \right\}, 
    \end{equation*}
    for $i = 1, \dots, t$. Moreover, $y_i^*>0$ if and only if $0<w_i <1/\psi$, where by convention, $1/\psi=+\infty$ when $\psi=0$.
\end{lemma}

\begin{definition}\label{def:h_psi}
    For $\psi \geq 0$,  define the function $h_\psi\,:\,(0,\infty) \rightarrow \mathbb{R}$ by 
    \begin{equation*}
        h_\psi(w) := \begin{cases}
            -\log(w) + \psi w - 1, &  \mbox{if } 0 < w < 1/\psi; \\
            \log(\psi), & \mbox{if } w \geq 1/\psi,
        \end{cases}
    \end{equation*}
    where, by convention, $1/\psi = +\infty$ when $\psi = 0$, so that $h_0(w) = -\log(w) -1$.
\end{definition}

\begin{theorem}\label{thm:subproblem_W}
    For $\Theta\in\mathbb{S}^k$, we have
\begin{equation*}
   \mathcal{L}_1^*(\Theta) =
   \left\{\begin{array}{ll} 
        \textstyle \sum_{i=k-t+1}^{k}h_\psi(\lambda_{i}(\Theta)),
        & \mbox{if } \Theta \in \mathbb{S}^k_{++}\,;\\
        +\infty, &\mbox{otherwise}.
        \end{array}\right.
\end{equation*}
\end{theorem}

\begin{proof}
    Let $w_1 \leq w_2 \leq \cdots \leq w_k$ denote the eigenvalues of $\Theta$, and let $Q$ be an orthogonal matrix whose $i$-th column is a unit eigenvector of $\Theta$ associated with $w_i$\,. The objective of \eqref{prob:subprob_W} depends on $W$ only through $\lambda(W)$ and through $\Theta\bullet W$, and the left-hand bound of Lemma~\ref{lem:von_neumann} gives $\Theta\bullet W\ge\sum_{i=1}^{k}w_i\lambda_i(W)$ for every $W\in\mathbb{S}^k_+$\,. Hence the supremum in \eqref{prob:subprob_W} is not decreased by replacing any $W$ by $Q\Diag(\lambda(W))Q^\top$, for which the attainment clause of Lemma~\ref{lem:von_neumann} gives equality. Writing $y_i:=\lambda_i(W)$, so that $y_1 \geq y_2 \geq \cdots \geq y_k \geq 0$, problem \eqref{prob:subprob_W} therefore reduces to 
\begin{equation}\label{supdfact_reform}
       \textstyle\sup_{y_1 \geq y_2\geq \dots \geq y_k\geq 0}\left\{\textstyle\sum_{i=1}^{t}\left(\log(y_i + \psi) - w_i y_i\right) - \sum_{i=t+1}^k w_iy_i\right\}.
    \end{equation}
    
    If $w_1 \leq 0$ (equivalently, $\Theta$ is not positive definite), then \eqref{supdfact_reform} is unbounded above. Indeed, fix $y_2=\dots=y_t=1$ and $y_{t+1}=\dots=y_k=0$, and let $y_1\ge1$; these values satisfy the ordering constraints, and the objective of \eqref{supdfact_reform} becomes $\log(y_1+\psi)-w_1y_1+\sum_{i=2}^{t}\left(\log(1+\psi)-w_i\right)$, where the sum is a finite constant. As $w_1\le0$, the objective tends to $+\infty$ as $y_1\to+\infty$. Hence, $\mathcal{L}_1^*(\Theta) =+\infty$.
    
    Now suppose that $w_1 > 0$ (equivalently, $\Theta \in \mathbb{S}^k_{++}$).  As the objective function of \eqref{supdfact_reform} is separable, the terms corresponding to $i = t+1,\dots,k$ are maximized at $y_i^* = 0$. For $i =1,\dots,t$,  Lemma~\ref{lem:subproblem_W_multivariable} yields the unique maximizer $y_i^* = \max\left\{0, 1/w_i - \psi \right\}$.
    Substituting these optimal values into the objective function of \eqref{supdfact_reform} gives $\mathcal{L}_1^*(\Theta) =\textstyle \sum_{i=1}^{t}h_\psi(w_i)$. Recalling that $w_i=\lambda_{k-i+1}(\Theta)$, the result follows.
\end{proof}

The following theorem provides a closed-form expression for $\mathcal{L}_2^*(\Theta, \upsilon, \nu, \tau)$.

\begin{theorem}\label{dfactlin}
    For $(\Theta,\upsilon,\nu,\tau)\in \mathbb{S}^k\times \mathbb{R}^{n}\times \mathbb{R}^{n}\times \mathbb{R}$, we have
    \begin{align*}
        &\textstyle\sup_{x}\left\{\Theta\bullet M_{\psi}(x) + \upsilon^\top x - \nu^\top x   - \tau \mathbf{e}^\top x +  \nu^\top \mathbf{e}  +\tau s \right\}\\
        &\quad = \left\{\begin{array}{ll}
        \nu^\top \mathbf{e} +\tau s,& \quad \mbox{if} ~ 
        \diag(F \Theta F^\top) + \upsilon - \nu   - \tau\mathbf{e}=0;\\
        +\infty, &\quad \mbox{otherwise.}
        \end{array}\right.
    \end{align*}
\end{theorem}
\begin{proof}
    The result  follows from the fact that an affine function is bounded above only when its linear part vanishes.
\end{proof}

Combining Theorems \ref{thm:subproblem_W} and \ref{dfactlin} yields the following explicit formulation of the  Lagrangian dual  of \ref{GFactplus}.
\begin{align}
    \begin{array}{lll}
    \zdgfactplusthing
    (C,s,t;\psi,F)= &\min&\sum_{\ell=k-t+1}^{k}h_\psi(\lambda_{\ell}(\Theta))
    + \nu^\top \mathbf{e} +\tau s\\[2pt]
    &\mbox{s.t.}&
         \diag(F \Theta F^\top) + \upsilon - \nu  - \tau\mathbf{e}=0,\\[1pt]
    &&\Theta\in \mathbb{S}^k_{++}\,, ~\upsilon\in\mathbb{R}^n_+\,, ~\nu\in\mathbb{R}^n_+\,,~\tau\in\mathbb{R}.
    \end{array}
    \tag{DGFact$^{+}$}\label{DGFactplus}
\end{align}

We next show that the optimal value of \ref{DGFactplus} is invariant under the choice of factorization
$C-\psi I = FF^\top$. 
The proof follows the approach of \citep*[Thm.~2.2]{FactPaper}, which proved the corresponding invariance property for the DFact bound for \ref{MESP}. As there, we first show, using eigenvalue interlacing, that appending zero columns to a factorization does not change the optimal value, and then relate any two factorizations with $n$ columns by an orthogonal change of variables. The only difference is that the objective of \ref{DGFactplus} applies $h_\psi$ to the $t$ least eigenvalues of $\Theta$, whereas that of DFact applies $-\log$ to its $s$ least eigenvalues; as $h_\psi$ is also nonincreasing, the interlacing argument carries over.

\begin{theorem}~\!\!\![Factorization invariance of DGFact$^+$]
\label{thm:dgfactplus_factorization_invariance}
    Let $C \in \mathbb{S}_ {++}^n$\,,  $0 \leq \psi \leq \lambda_{\min}(C)$, and let
$
        C-\psi I=F_1F_1^\top=F_2F_2^\top
$    be two factorizations, where $F_j\in\mathbb R^{n\times k_j}$, with $\max\{t,\rank(C-\psi I)\}\le k_j\le n$, for $j=1,2$. Then,
    \begin{equation*}
        \zdgfactplusthing(C,s,t;\psi,F_1)
        =
        \zdgfactplusthing(C,s,t;\psi,F_2).
    \end{equation*}
\end{theorem}

\begin{proof}
First, we show that appending zero columns does not change the optimal value of \ref{DGFactplus}.
Let $\bar F:=[F,0]$, obtained by appending $\bar{k}-k$ zero columns to $F\in\mathbb{R}^{n\times k}$, for any $\bar{k}$ with $k<\bar{k}\le n$. Let $(\Theta,\upsilon,\nu,\tau)$ be feasible for \ref{DGFactplus} associated with $F$. Choose any
$
L\ge\lambda_1(\Theta)
$
and define
\[
\bar\Theta:=
\begin{pmatrix}
\Theta&0\\
0&LI
\end{pmatrix} \in \mathbb{S}^{\bar{k}}.
\]
Then $\bar\Theta\in \mathbb{S}^{\bar{k}}_{++}$ and $\bar F\bar\Theta\bar F^\top
=
F\Theta F^\top$.
As all newly introduced eigenvalues are equal to $L$, they are not less than any eigenvalue of $\Theta$. Hence the $t$ least eigenvalues of $\bar\Theta$ coincide with those of $\Theta$, implying
\[
\textstyle\sum_{\ell=\bar k-t+1}^{\bar k}
h_\psi(\lambda_\ell(\bar\Theta))
=
\sum_{\ell=k-t+1}^{k}
h_\psi(\lambda_\ell(\Theta)).
\]
Therefore, $(\bar\Theta,\upsilon,\nu,\tau)$ is feasible for the formulation associated with $\bar F$ and has the same objective value.

Conversely, let $(\bar\Theta,\upsilon,\nu,\tau)$ be feasible for the formulation associated with $\bar F$, and partition
\[
\bar\Theta=
\begin{pmatrix}
\Theta_{11}&\Theta_{12}\\
\Theta_{12}^\top&\Theta_{22}
\end{pmatrix},
\]
where $\Theta_{11}\in\mathbb{S}^k$. As $\bar\Theta\in \mathbb{S}^{\bar{k}}_{++}$\,, we also have $\Theta_{11}\in \mathbb{S}^k_{++}$\,, and $
\bar F\bar\Theta\bar F^\top
=
F\Theta_{11}F^\top$.
Therefore, $(\Theta_{11},\upsilon,\nu,\tau)$ is feasible for the formulation associated with $F$.

As only the $t$ least eigenvalues of $\Theta_{11}$ are relevant for the objective value, it is necessary to compare the spectra of $\Theta_{11}$ and $\bar\Theta$. By the Cauchy interlacing theorem,
\[
\lambda_{k-t+i}(\Theta_{11})
\ge
\lambda_{\bar k-t+i}(\bar\Theta),
\qquad
i=1,\ldots,t.
\]
As $h_\psi$ is nonincreasing on $(0,\infty)$,
\[
\textstyle \sum_{\ell=k-t+1}^{k}
h_\psi(\lambda_\ell(\Theta_{11}))
\le
\sum_{\ell=\bar k-t+1}^{\bar k}
h_\psi(\lambda_\ell(\bar\Theta)).
\]
Thus the feasible solution obtained by deleting the zero columns has objective value no greater than the original one.

The two preceding constructions show that appending zero columns leaves the optimal value of \ref{DGFactplus} unchanged. Therefore, without loss of generality, we may assume that both factorizations $F_1$ and $F_2$ have $n$ columns.

As $F_1F_1^\top
=
F_2F_2^\top$,
we may choose singular-value decompositions
$
F_j=U\Sigma V_j^\top$, for 
$j=1,2$,
with the same matrices $U$ and $\Sigma$. Hence $
F_1=F_2Q$, with $
Q:=V_2V_1^\top$,
where $Q$ is orthogonal.

Now let $(\Theta_1,\upsilon,\nu,\tau)$ be feasible for \ref{DGFactplus} associated with $F_1$, and define
$
\Theta_2:=Q\Theta_1Q^\top$.
As orthogonal similarity preserves positive definiteness and eigenvalues,
$
\Theta_2\in \mathbb{S}^n_{++}$\,, and $
\lambda_\ell(\Theta_2)=\lambda_\ell(\Theta_1)$, for $
\ell=1,\ldots,n.
$
Moreover,
$
F_2\Theta_2F_2^\top
=
F_2Q\Theta_1Q^\top F_2^\top
=
F_1\Theta_1F_1^\top,
$
so that
$
\diag(F_2\Theta_2F_2^\top)
=
\diag(F_1\Theta_1F_1^\top).
$
Therefore $(\Theta_2,\upsilon,\nu,\tau)$ is feasible for the formulation associated with $F_2$, and has exactly the same objective value as $(\Theta_1,\upsilon,\nu,\tau)$.
Thus every feasible solution for the formulation associated with $F_1$ corresponds to a feasible solution for the formulation associated with $F_2$ with identical objective value. Reversing the roles of $F_1$ and $F_2$ yields the converse correspondence, completing the proof.
\end{proof}

In view of the invariance of the bound with respect to the factorization established in the previous theorem, we assume henceforth that $k=n$ and set $F:=(C-\psi I_n)^{\scriptscriptstyle 1/2}$, unless specifically stated otherwise. Accordingly, as the bound depends on the augmentation parameter but not on the particular factorization, from now on we write $\zdgfactplusthing(C,s,t;\psi)$ in place of $\zdgfactplusthing(C,s,t;\psi,F)$, and likewise for the other bounds introduced below.

\subsection{Duality for \texorpdfstring{\ref{DGFactplus}}{DGFactplus}}\label{subsec:dual_dgfactplus}
Although \ref{DGFactplus} can be solved directly to compute the augmented generalized factorization bound, it is computationally more efficient to work with its Lagrangian dual, which is derived in this subsection. The Lagrangian dual of DGFact 
was derived by \citep*{GMESP_Alg}; however, the presence of the augmentation term requires a substantially different analysis. Rather than adapting their derivation, we develop a new  and more direct approach. As in the derivation of \ref{DGFactplus}, von Neumann's trace inequality (Lemma~\ref{lem:von_neumann}) plays a central role. Combined with our analysis, it leads to Theorem~\ref{thm:subproblem_Theta}, which provides a concise characterization of the dual subproblem. Similarly to the derivation of  \ref{DGFactplus},  the  derivation presented here is shorter and more self-contained than that of \citep*{GMESP_Alg}.

After eliminating $\upsilon$, the Lagrangian associated with \ref{DGFactplus} is
\begin{align*}
    &\mathcal{L}(\Theta,\nu,\tau,x,y):=
    \textstyle 
    \sum_{\ell=k-t+1}^k h_\psi(\lambda_\ell(\Theta)) + \nu^\top \mathbf{e} +\tau s +x^\top\left(\diag(F \Theta F^\top)  \!-\! \nu   \!-\! \tau\mathbf{e}\right) -y^\top\nu ,
\end{align*}
with $\mbox{ dom }\mathcal{L} =  \mathbb{S}_{++}^k\times  \mathbb{R}^{n}\times  \mathbb{R} \times \mathbb{R}^{n}\times \mathbb{R}^{n}$.
The corresponding dual function is
\begin{equation*} 
    \mathcal{L}^*(x,y):=\textstyle\inf_{\Theta\in~\mathbb{S}^{k}_{++}\,,\nu,\tau}\mathcal{L}(\Theta,\nu,\tau,x,y),
\end{equation*}
and the  Lagrangian dual problem is
\begin{equation*}
    \zddgfactplusthing(C,s,t;\psi):=\max\{\mathcal{L}^*(x,y)~:~x\geq 0,~y\geq 0\}.
\end{equation*}
We note that \ref{DGFactplus} has a strictly-feasible solution (e.g., given by $(\hat \Theta :=I,\,\hat \upsilon:=\mathbf{e},\,\hat \nu:=\diag(C -\psi I)+\mathbf{e},\,\hat \tau:=0)$). Then, Slater's condition holds for \ref{DGFactplus} and   we are justified to use maximum in the formulation of the Lagrangian dual problem, rather than supremum, as the optimal value of the Lagrangian dual problem is attained.
The dual function separates as
\begin{align*}
    \textstyle\inf_{\Theta\in\mathbb{S}^{k}_{++}\,,\nu,\tau}\mathcal{L}(\Theta,\nu,\tau,x,y) &=  
\textstyle\inf_{\Theta\in\mathbb{S}^{k}_{++}} \left\{ 
    \textstyle 
    \sum_{\ell=k-t+1}^k h_\psi(\lambda_{\ell}(\Theta))
     +x^\top \diag(F \Theta F^\top)\right\}
     \\
    &\qquad\quad + 
     \textstyle\inf_{\nu,\tau}
     \left\{\nu^\top (\mathbf{e} - x - y)   +\tau (s - \mathbf{e}^\top x)\right\}.
\end{align*}
Furthermore, we have
\begin{align}
    \textstyle\inf_{\Theta\in\mathbb{S}^{k}_{++}} &\left\{ 
    \textstyle 
    \sum_{\ell=k-t+1}^k h_\psi(\lambda_{\ell}(\Theta))
     +x^\top \diag(F \Theta F^\top)\right\}\nonumber\\
     &= \textstyle\inf_{\Theta\in\mathbb{S}^{k}_{++}} \left\{ 
    \textstyle 
    \sum_{\ell=k-t+1}^k h_\psi(\lambda_{\ell}(\Theta))
     + (F^\top \Diag(x) F) \bullet  \Theta\right\}\nonumber\\
      &= \textstyle\inf_{0 < w_1 \leq \dots \leq w_k} \left\{ 
    \textstyle 
    \sum_{\ell=1}^t h_\psi(w_\ell)
     + \sum_{\ell=1}^k \lambda_{\ell}(F^\top \Diag(x) F) w_\ell\right\},
     \label{eq:inf_Theta_to_w}
\end{align}
where the last equality follows from Lemma~\ref{lem:von_neumann}, with $A:=F^\top \Diag(x) F$. Indeed, for every $\Theta\in\mathbb{S}^k_{++}$\,, writing its eigenvalues in nondecreasing order as $w_1\le\cdots\le w_k$\,, the left-hand bound gives $A\bullet\Theta\ge\sum_{\ell=1}^k\lambda_\ell(A)w_\ell$\,, so the objective at $\Theta$ is at least the objective at $w$; hence the left-hand infimum is at least the right-hand one. Conversely, for every $w$ with $0<w_1\le\cdots\le w_k$\,, let $Q$ be an orthogonal matrix with $A=Q\Diag(\lambda(A))Q^\top$, and let $\Theta:=Q\Diag(w)Q^\top\in\mathbb{S}^k_{++}$\,; the attainment clause gives $A\bullet\Theta=\sum_{\ell=1}^k\lambda_\ell(A)w_\ell$\,, so the objective at $\Theta$ equals the objective at $w$; hence the left-hand infimum is at most the right-hand one.

We next analyze the above optimization problem by presenting {a} candidate primal solution and Lagrange multipliers and showing that they satisfy the KKT conditions. We begin by presenting a key result and definition that will be used throughout the remainder of the paper.

\begin{lemma}~\label{lem:index_iota}(see \citep*[Lem.~13]{Nikolov})
     Let $y\in\mathbb{R}_+^k$ with $y_1\geq y_2\geq \cdots\geq y_k$\,, and let $0<t\leq k$. There exists a unique integer $\hat\iota$, with $0\leq \hat\iota< t$, such that
     \begin{equation*}
     \textstyle
     y_{\hat\iota}>\frac{1}{t-\hat\iota}
     \textstyle 
     \sum_{j=\hat\iota+1}^k y_{j}\geq y_{\hat\iota +1}\,,
     \end{equation*}
     with the convention $y_0=+\infty$.
\end{lemma}

\begin{definition}~\label{def:phi_t}(see \citep*[Def.~2]{Nikolov})
     Let $y\in\mathbb{R}_+^k$ with $y_1\geq y_2\geq \cdots\geq y_k$ and let $0<t\leq k$. Let 
   $\hat\iota$ be the unique integer  defined by Lemma~\ref{lem:index_iota}. We define
    \begin{equation}
    \textstyle
        \phi_t(y) := \sum_{j=1}^{\hat\iota}\log(y_j) + (t-\hat\iota)\log\!\left(\textstyle\frac{1}{t-\hat\iota}
     \sum_{j=\hat\iota+1}^k y_{j}\right).
    \end{equation}
    Furthermore, for $\psi\ge0$ and $X\in\mathbb{S}_{+}^k$, extending the definition of \citep*{li2025augmented} from $t=s$ to $t\le s$, we define
$\Gamma_t^+(X; \psi) := \phi_t(\lambda(X) + \psi\mathbb{I}_t)$. 
\end{definition}

\begin{theorem}\label{thm:subproblem_Theta}
    Let $\lambda\in\mathbb{R}^k_+$ with $
        \lambda_1 \geq \cdots \geq \lambda_k$\,. Let $0<t\leq k$ and  either $\psi>0$, or $\psi=0$ and $\lambda_t>0$. 
Let
    \begin{align*}
&        \hat y_\ell :=
        \begin{cases}
            \lambda_\ell+\psi, & \mbox{for } \ell = 1,\dots,t;\\
            \lambda_\ell\,, & \mbox{for } \ell = t+1,\dots,k,
        \end{cases}
&        \hat w_\ell :=
        \begin{cases}
           \textstyle 1/\hat{y}_\ell\,, & \mbox{for } \ell = 1,\dots,\hat\iota;\\
           \textstyle 1/\bar{y}, & \mbox{for } \ell = \hat\iota+1,\dots,k,
        \end{cases}
    \end{align*}
where $\hat\iota$ is the unique index corresponding to $y:=\hat y$ that satisfies Lemma~\ref{lem:index_iota} and $\bar y:=
\textstyle\frac{1}{t-\hat\iota}
     \sum_{j=\hat\iota+1}^k \hat{y}_{j}$\,.  Then, $\hat w$ is an optimal solution to 
    \begin{equation}\label{probinw}
       \textstyle\inf_{0<w_1\le\cdots\le w_k}
        \left\{
        \sum_{\ell=1}^t h_\psi(w_\ell)
        +
        \sum_{\ell=1}^k \lambda_\ell w_\ell
        \right\}{,}
    \end{equation}
with objective value equal to $\phi_t(\hat y)$.
\end{theorem}

\begin{proof}

We show that $\hat w$ is optimal by exhibiting multipliers $\hat\eta$ for the ordering constraints such that $(\hat w,\hat\eta)$ satisfies KKT conditions that are sufficient for optimality in the convex problem \eqref{probinw} (steps~1--3), and then we evaluate the objective at $\hat w$ (step~4).

{
{1. Convexity and differentiability.}
The function $h_\psi$ is continuously differentiable and convex on $(0,\infty)$. Indeed, for $\psi>0$, we have $h_\psi'(w)=\psi-1/w$ for $0<w<1/\psi$ and $h_\psi'(w)=0$ for $w>1/\psi$, and both one-sided derivatives at $w=1/\psi$ are equal to zero. Hence
$
h_\psi'(w)=\min\{\psi-1/w,\,0\}
$
for all $w>0$ (also when $\psi=0$, by the convention $1/\psi=+\infty$), which is continuous and nondecreasing in $w$. Therefore, the objective function $f$ of \eqref{probinw} is convex and continuously differentiable on the convex feasible set
$
\mathcal{W}:=\{w\in\mathbb{R}^k:\,0<w_1\le\cdots\le w_k\}
$,
with
$\partial f/\partial w_\ell=h_\psi'(w_\ell)+\lambda_\ell$\,, for $\ell=1,\ldots,t$, and $\partial f/\partial w_\ell=\lambda_\ell$\,, for $\ell=t+1,\ldots,k$.

{2. Sufficient optimality conditions.}
Let $\hat w\in\mathcal{W}$, and suppose that there exist $\hat\eta_1,\ldots,\hat\eta_{k-1}\ge0$ such that, with $\hat\eta_0=\hat\eta_k:=0$,
\begin{align*}
&\textstyle\frac{\partial f}{\partial w_\ell}(\hat w)+\hat\eta_\ell-\hat\eta_{\ell-1}=0, \quad \ell=1,\ldots,k, &&\mbox{(stationarity)}\\
&\hat\eta_\ell(\hat w_\ell-\hat w_{\ell+1})=0, \quad \ell=1,\ldots,k-1. &&\mbox{(complementary slackness)}
\end{align*}
Then $\hat w$ is an optimal solution of \eqref{probinw}. Indeed, for every $w\in\mathcal{W}$, the convexity of $f$, stationarity and complementary slackness give
\begin{align*}
f(w)&\textstyle\ge f(\hat w)+\sum_{\ell=1}^{k}\frac{\partial f}{\partial w_\ell}(\hat w)(w_\ell-\hat w_\ell)
=f(\hat w)-\sum_{\ell=1}^{k-1}\hat\eta_\ell\left((w_\ell-w_{\ell+1})-(\hat w_\ell-\hat w_{\ell+1})\right)\\
&\textstyle=f(\hat w)-\sum_{\ell=1}^{k-1}\hat\eta_\ell(w_\ell-w_{\ell+1})\ge f(\hat w),
\end{align*}
because $\hat\eta_\ell\ge0$ and $w_\ell\le w_{\ell+1}$\,. These are the KKT conditions for \eqref{probinw}, with multipliers $\eta_\ell$ associated with the ordering constraints $w_\ell-w_{\ell+1}\le0$; as only their sufficiency is used, no constraint qualification is needed. Moreover, if $\hat w_\ell\le1/\psi$ for $\ell=1,\ldots,t$, then stationarity reads
$
-\frac{1}{\hat w_\ell}+\psi+\lambda_\ell+\hat\eta_\ell-\hat\eta_{\ell-1}=0$, for
$\ell=1,\ldots,t$, and $
\lambda_\ell+\hat\eta_\ell-\hat\eta_{\ell-1}=0$, for
$\ell=t+1,\ldots,k$.
}

{3. Verification that the candidate solution satisfies the KKT conditions.}
Define the Lagrange multipliers
\[
\hat\eta_\ell:=
\begin{cases}
0,&\ell=1,\dots,\hat\iota;\\
\textstyle
\sum_{j=\hat\iota+1}^{\ell}
\left(\bar{y}-\hat y_j\right),
&\ell=\hat\iota+1,\dots,t;\\[1mm]
\textstyle
\sum_{j=\hat\iota+1}^{t}
\left(\bar{y}-\hat y_j\right)
-
\textstyle\sum_{j=t+1}^{\ell}\hat y_j\,,
&\ell=t+1,\dots,k.
\end{cases}
\]
In step~4 below, we show that $\hat w_\ell\le1/\psi$, for $\ell=1,\ldots,t$. The verification that $(\hat w,\hat\eta)$ satisfies the conditions of step~2 is routine and is given in Lemma~\ref{lem:kkt_verification} in the Appendix.

    {4. Optimal value.}
 It remains to evaluate the objective at $\hat w$.
For $\ell\le\hat\iota$, we have
$\hat w_\ell=1/\hat y_\ell$\,, where
$\hat y_\ell=\lambda_\ell+\psi$. We first observe that $\bar y\ge\psi$  and $\bar y>0$. Indeed, $\hat y_j\ge\psi$ for $\hat\iota<j\le t$ and $\hat y_j\ge0$ for $j>t$, so that 
$\sum_{j=\hat\iota+1}^{k}\hat y_j\ge(t-\hat\iota)\psi$; if $\psi>0$, this gives $\bar y>0$, and if $\psi=0$, then $\hat y_j=\lambda_j\ge\lambda_t>0$ for $\hat\iota<j\le t$, so that $\bar y\ge\lambda_t>0$. Lemma~\ref{lem:index_iota} then gives
$\lambda_\ell+\psi=\hat y_\ell\ge \hat y_{\hat\iota}>\bar y\ge\psi$, for $\ell\le\hat\iota$, so that $\lambda_\ell>0$ and therefore $\hat w_\ell=1/\hat y_\ell< 1/\psi$
and 
\[
\begin{aligned}
&h_\psi(\hat w_\ell)+\lambda_\ell\hat w_\ell =
-\log\!\left(\textstyle\frac1{\hat y_\ell}\right)
+\textstyle\frac{\psi+\lambda_\ell}{\hat y_\ell}
-1 = 
\log\hat y_\ell\,.
\end{aligned}
\]
Therefore,
\begin{equation}\label{part1}
\textstyle \sum_{\ell=1}^{\hat\iota}
\bigl(h_\psi(\hat w_\ell)+\lambda_\ell\hat w_\ell\bigr)
=
\sum_{\ell=1}^{\hat\iota}\log\hat y_\ell\,.
\end{equation}
For $\hat\iota<\ell\le k$, we have
$\hat w_\ell=1/\bar y$, and $\bar y\ge\psi$ gives $\hat w_\ell\le1/\psi$. As $h_\psi(w)=-\log(w)+\psi w-1$ holds at $w=1/\psi$ as well, by continuity, we obtain 
\begin{align}
&\textstyle\sum_{\ell=\hat\iota+1}^{t}h_\psi(\hat w_\ell)
+
\textstyle\sum_{\ell=\hat\iota+1}^{k}\lambda_\ell\hat w_\ell
=
(t-\hat\iota)(-\log\!\left(\textstyle\frac1{\bar y}\right)+\frac\psi{\bar y}-1)
+\frac1{\bar y}\textstyle\sum_{\ell=\hat\iota+1}^{k}\lambda_\ell\nonumber\\
&\qquad=
(t-\hat\iota)\textstyle\left(-\log\left(\frac1{\bar y}\right)\right)
+\textstyle\frac1{\bar y}\sum_{\ell=\hat\iota+1}^{k}\hat y_\ell
-(t-\hat\iota)
=(t-\hat\iota)
        \log \bar y,\label{part2}
\end{align}
where the last equality follows from the identity
$
\frac1{\bar y}
\sum_{\ell=\hat\iota+1}^{k}\hat y_\ell
=
t-\hat\iota$.

Combining \eqref{part1} and \eqref{part2} completes the proof.
\end{proof}

We next derive a closed-form expression for the linear subproblem in the Lagrangian dual problem.

\begin{theorem}\label{thm:subproblem_xy}
    For $(x,y)\in  \mathbb{R}^{n}\times \mathbb{R}^{n}$, we have
    \begin{align*}
        &\inf_{\nu,\tau}
         \nu^\top (\mathbf{e} - x - y)   +\tau (s - \mathbf{e}^\top x)= \left\{\begin{array}{cl}
        0,& \quad \mbox{if  } ~ \mathbf{e} - x - y=0,~s - \mathbf{e}^\top x=0;\\
        -\infty, &\quad \mbox{otherwise.}
        \end{array}\right.
    \end{align*}
\end{theorem}
\begin{proof}
    The result  follows from the fact that an affine function is bounded below only when its linear part vanishes.
\end{proof}

Considering  Theorems \ref{thm:subproblem_Theta} and \ref{thm:subproblem_xy}, the Lagrangian dual of \ref{DGFactplus} is given  by 

\begin{align}
    \begin{array}{lll}
    \zddgfactplusthing
    (C,s,t;\psi)= &\max&   \Gamma_t^+(M_\psi(x); \psi)\\[1pt]
    &\mbox{s.t.}&\!\!
         \mathbf{e}^\top x = s, 0 \leq x \leq \mathbf{e}.
    \end{array}
    \tag{DDGFact$^{+}$}\label{DDGFactplus}
\end{align}
{When $\psi=0$, the hypothesis $\lambda_t>0$ of Theorem~\ref{thm:subproblem_Theta} holds for $\lambda:=\lambda(M_0(x))$ at every $x$ satisfying $\mathbf{e}^\top x=s$ and $0\le x\le\mathbf{e}$ (for any other $x\ge0$, the dual function equals $-\infty$ by Theorem~\ref{thm:subproblem_xy}). Indeed, as $C\in\mathbb{S}^n_{++}$\,, $F$ is nonsingular, and such an $x$ has at least $s\ge t$ positive components, so that $\rank(M_0(x))=\rank(\Diag(x))\ge t$. The same observation applies to every feasible solution of the fixed-variable formulations considered in Section~\ref{subsec:feasible_dual_sol_dgfactplus}.}

From Lagrangian duality, we have that \ref{DDGFactplus} is a convex optimization problem. Moreover, because \ref{DGFactplus} has a strictly feasible solution, and \ref{DDGFactplus} has a feasible solution with finite objective value (e.g., $\hat x:=\frac{s}{n}\mathbf{e}$, for which $\Gamma_t^+(M_\psi(\hat x);\psi)$ 
is finite because $\lambda_t(M_\psi(\hat x))+\psi>0$: if $\psi=0$, then $M_0(\hat x)=\frac sn C$ is positive definite, and if $\psi>0$, the inequality is immediate), strong
duality holds between \ref{DDGFactplus} and \ref{DGFactplus}{, that is, $\zdgfactplusthing(C,s,t;\psi)=\zddgfactplusthing(C,s,t;\psi)$}.

We note that, unlike the nonconvex relaxation \ref{GFactplus},  \ref{DDGFactplus} is not an exact relaxation in general. However, in the following particular case,  we show that it is.

\begin{theorem}~\!\!\![Sufficient condition under which  \ref{DDGFactplus} is an exact relaxation]
\label{thm:gamma_t_collapse_multiplicity}
    Let $C\in\mathbb S_ {++}^n$\,, let $1\le t\le s< n$, and let $\psi:=\lambda_{\min}(C)$\,. Assume that the multiplicity $\mu_{\min}$ of the least eigenvalue of $C$ satisfies
$
\mu_{\min}\ge n-t.
$ 
 Then, for every feasible solution $\hat x$ of \ref{GFactplus}, $\Gamma_t^+(M_\psi(\hat x);\psi)
        =
        f_t(M_\psi(\hat x);\psi)$.
\end{theorem}

\begin{proof}
    Let
$
C-\psi I=FF^\top
$
be any factorization with 
$F\in\mathbb{R}^{n \times k}$,
for $\max\{t,\rank(C-\psi I)\}\le k\le n$ (here we do not use the normalization $k=n$). The multiplicity assumption implies $\rank(F)=\rank(C-\psi I)\leq t$.     Then,
    \begin{equation*}
        \rank(M_\psi(\hat x))
        =
        \rank(F^\top\Diag(\hat x)F)
        \leq
        \rank(F)
        \leq
        t.
    \end{equation*}
    Let $y:=\lambda(M_\psi(\hat x))+\psi\mathbb{I}_t$, and let $\hat\iota$ be the
    index associated with $y$ in Lemma~\ref{lem:index_iota}, so that
    $0\le\hat\iota<t$ and, in particular, $t-\hat\iota\ge1$.
    As $\rank(M_\psi(\hat x))\le t$, we have
    $\lambda_{t+1}(M_\psi(\hat x))=\cdots=\lambda_{k}(M_\psi(\hat x))=0$, and
    therefore $y_{t+1}=\cdots=y_k=0$. Consequently,
    \begin{equation*}
        \textstyle\frac{1}{t-\hat\iota}\sum_{j=\hat\iota+1}^{k} y_j
        =
        \textstyle\frac{1}{t-\hat\iota}\sum_{j=\hat\iota+1}^{t} y_j
    \end{equation*}
    is the arithmetic mean of the $t-\hat\iota$ numbers
    $y_{\hat\iota+1},\ldots,y_t$\,, each of which is at most $y_{\hat\iota+1}$
    because the components of $y$ are sorted in nonincreasing order. Hence
    \begin{equation*}
        y_{\hat\iota +1}
        \le
        \textstyle\frac{1}{t-\hat\iota}\sum_{j=\hat\iota+1}^{k} y_j
        =
        \textstyle\frac{1}{t-\hat\iota}\sum_{j=\hat\iota+1}^{t} y_j
        \le
        y_{\hat\iota +1}\,,
    \end{equation*}
    where the first inequality follows from Lemma~\ref{lem:index_iota}.
    Therefore all the inequalities above hold with equality; as the mean of
    $y_{\hat\iota+1},\ldots,y_t$ equals their greatest element, we conclude that
    \begin{equation*}
        y_{\hat\iota+1}
        =
        \cdots
        =
        y_t
        =
        \textstyle\frac{1}{t-\hat\iota}\sum_{j=\hat\iota+1}^{k} y_j.
    \end{equation*}
    Substituting this equality into the definition of $\phi_t$ yields
    \begin{align*}
        \Gamma_t^+(M_\psi(\hat x);\psi)
        &=
        \textstyle\sum_{j=1}^{\hat\iota}\log(y_j)
        +
        (t-\hat\iota)
        \textstyle\log\left(
            \textstyle\frac{\textstyle\sum_{j=\hat\iota+1}^k y_j}{t-\hat\iota} 
        \right)=
        \textstyle\sum_{j=1}^t\log(y_j)=
        f_t(M_\psi(\hat x);\psi),
    \end{align*}
    which completes the proof.
\end{proof}

\section{Monotonicity in the augmentation parameter}
\label{sec:monotonicity_ddgfactplus}

In this section we prove that the optimal value of \ref{DDGFactplus} is monotonically nonincreasing with respect to the augmentation parameter $\psi$. This implies that the best choice is $\psi=\lambda_{\min}(C)$. 

The proof combines three ingredients. First, we show that the function $\phi_t$ is componentwise nondecreasing. Next, we derive from \citep*[Lem.~1]{li2025augmented} a majorization result adapted to the setting $t\le s$. Finally, we show that $\phi_t$ is Schur-concave (Lemma~\ref{lem:phi_t_majorization}). Combining these three ingredients, we establish the desired monotonicity result in Theorem~\ref{thm:ddgfact_plus_monotonicity}. Here, our proof differs from that of \citep*{li2025augmented}. Their argument relies on the identity $t=s$. Consequently, the original proof does not extend directly. We overcome this by applying Li's majorization lemma in the modified form of Lemma~\ref{lem:lambda_modified_shift}, in which the shift on the components $t+1,\ldots,s$ is removed from both vectors, and by compensating for the removed shift through the componentwise monotonicity of $\phi_t$ (Corollary~\ref{cor:phi_t_monotone}).

\begin{lemma}~\!\!\![Supergradient of $\phi_t$]\label{lem:phi_t_subgrad}~(proof in Appendix{; the result can also be derived from \citep*[Lem.~11, Lem.~16 and the proof of Thm.~14]{Nikolov}; cf.} \citep*[Remark~2]{li2025augmented})
Let $y\in\mathbb{R}_+^k$ with $y_1\geq y_2\geq \cdots\geq y_k$\,, let $0<t\leq k$, and assume that $y_t>0$. Let
$\hat\iota$ be the unique integer defined by Lemma~\ref{lem:index_iota} and $\bar y:=\frac{1}{t-\hat\iota}\sum_{j=\hat\iota+1}^k y_{j}$\,, and let
\begin{align*}
g_i := \begin{cases}
1/y_i\,, &\mbox{for }i=1,\dots,\hat\iota;\\
1/\bar{y}, &\mbox{for }i=\hat\iota+1,\dots,k.
\end{cases}
\end{align*}
Then $g$ is well defined, $0<g_1\le\cdots\le g_k$\,, and
\[
\phi_t(z)\le\phi_t(y)+g^\top(z-y),
\]
for every $z\in\mathbb{R}^k_+$ with $z_1\geq\cdots\geq z_k$\,, where, with the convention $\log0=-\infty$, $\phi_t(z)=-\infty$ whenever $z_t=0$.
\end{lemma}

\begin{corollary}~\!\!\![Componentwise monotonicity of $\phi_t$]
\label{cor:phi_t_monotone}
Let $y,z \in \mathbb{R}_{+}^k$ with $y_1\geq y_2\geq \cdots\geq y_k$ and $z_1\geq z_2\geq \cdots\geq z_k$\,. Assume that 
$z_i\le y_i$\,, for $i=1,\ldots,k$. Let $0<t\leq k$. 
Then
$\phi_t(z)\le\phi_t(y)$.
Hence, $\phi_t$ is componentwise nondecreasing.
\end{corollary}
\begin{proof}
If $y_t=0$, then $z_t=0$, so that $\phi_t(z)=-\infty\le\phi_t(y)$ by Lemma~\ref{lem:phi_t_subgrad}. Otherwise, Lemma~\ref{lem:phi_t_subgrad} gives a vector $g\ge0$ such that
$\phi_t(z)
\le
\phi_t(y)+g^\top(z-y)$. Because $z\le y$ componentwise and $g\ge 0$, we have
$g^\top(z-y)\le0$. Therefore, 
$\phi_t(z)\le\phi_t(y)$.
\end{proof}

Besides Corollary~\ref{cor:phi_t_monotone}, the proof of Theorem~\ref{thm:ddgfact_plus_monotonicity} uses majorization and Schur-concavity, which we now recall.

\begin{definition}~\!\!\![Majorization; see \citep*[Ch. 1, Def. A.1]{marshall2011inequalities}]For $v\in\mathbb{R}_+^m$, we denote by $v^{\downarrow}$ the vector with the components of $v$ sorted in nonincreasing order. For $x,y\in\mathbb{R}_+^m$, we say that $y$ majorizes $x$, and write
$y \succ x$,
if
\[ 
\begin{array}{l}
\textstyle \sum_{i=1}^{\ell}x^{\downarrow}_i
\le
\sum_{i=1}^{\ell}y^{\downarrow}_i\,,
\qquad
\ell=1,\ldots,m-1,\\[2pt]
\textstyle\sum_{i=1}^{m}x_i
=
\sum_{i=1}^{m}y_i\,.
\end{array}
\]
In particular, if $x$ and $y$ are sorted in nonincreasing order, then $x^{\downarrow}=x$ and $y^{\downarrow}=y$.
\end{definition}

A function $f:\mathbb{R}^k_+\to[-\infty,\infty)$ is Schur-concave if $f(y)\le f(x)$ whenever $y\succ x$; see \citep*[Ch.~3, Def.~A.1]{marshall2011inequalities}.

\begin{lemma}~\!\!\![Schur-concavity of $\phi_t$]\label{lem:phi_t_majorization}
Let $0<t\le k$, and define $\tilde\phi_t:\mathbb{R}^k_+\to[-\infty,\infty)$ by $\tilde\phi_t(y):=\phi_t(y^{\downarrow})$, with the convention $\log0=-\infty$. Then $\tilde\phi_t$ is symmetric, concave and Schur-concave. In particular, if $a,b\in\mathbb{R}^k_+$ satisfy $a_1\ge\cdots\ge a_k$\,, $b_1\ge\cdots\ge b_k$ and $a\succ b$, then $\phi_t(a)\le\phi_t(b)$.
\end{lemma}
\begin{proof}
For $X\in\mathbb{S}^k_+$, let $\Gamma_t(X):=\phi_t(\lambda(X))$ $(=\Gamma_t^+(X;0))$, which is the function $\Gamma_j$ of \citep*[Def.~2]{Nikolov} with $j:=t$ and $d:=k$. By \citep*[Lem.~11, Lem.~16 and the proof of Thm.~14]{Nikolov}, for every $X\in\mathbb{S}^k_+$,
\[
\Gamma_t(X)=\inf_{W\in\mathbb{S}^k_{++}}\big\{\Delta_t(W)+X\bullet W\big\}-t,
\]
where $\Delta_t(W)$ is minus the sum of the logarithms of the $t$ least eigenvalues of $W$ \citep*[Def.~1]{Nikolov}, and where both sides equal $-\infty$ when $\rank(X)<t$. Hence $\Gamma_t$, as an infimum of affine functions of $X$, is concave on $\mathbb{S}^k_+$. As $\lambda(\Diag(y))=y^{\downarrow}$, we have $\tilde\phi_t(y)=\Gamma_t(\Diag(y))$, so that $\tilde\phi_t$ is concave on $\mathbb{R}^k_+$, and it is clearly symmetric. Now let $x,y\in\mathbb{R}^k_+$ with $y\succ x$. By \citep*[Ch.~2, Thm.~B.2 and Thm.~A.2]{marshall2011inequalities}, $x=\sum_{P}\alpha_P\,Py$, where the sum is over the permutation matrices $P$ of order $k$, $\alpha_P\ge0$ and $\sum_P\alpha_P=1$. By the concavity and symmetry of $\tilde\phi_t$,
\[
\textstyle\tilde\phi_t(x)\ge\sum_P\alpha_P\,\tilde\phi_t(Py)=\tilde\phi_t(y).
\]
This is the argument of \citep*[Ch.~3, Prop.~C.2]{marshall2011inequalities}, which we repeat because $\tilde\phi_t$ may take the value $-\infty$. The last statement follows because $\tilde\phi_t(a)=\phi_t(a)$ and $\tilde\phi_t(b)=\phi_t(b)$ for sorted $a$ and $b$.
\end{proof}

\begin{lemma}~\!\!\![Shift]
\label{lem:lambda-shift-s}
(see \citep*[Lem.~1]{li2025augmented})
Let $C\in\mathbb{S}^n_{++}$. Let $0\!\le\!\psi_1\!\le\!\psi_2\!\le\!\lambda_{\min}(C)$ and $x\in[0,1]^n$ satisfy $\mathbf e^\top x=s$. Let $\pi$ be a permutation of $N$ such that
$
x_{\pi(1)}
\ge
x_{\pi(2)}
\ge
\cdots
\ge
x_{\pi(n)}\,.
$
Then,
\begin{enumerate}
\item[($i$)]
$
\lambda(M_{\psi_2}(x))
+\psi_2\mathbb I_s
\succ
\lambda(M_{\psi_1}(x))
+\psi_1\mathbb I_s$;

\item[($ii$)]
$
\textstyle\sum_{i=1}^{\ell}\lambda_i(M_{\psi_2}(x))
+
(\psi_2-\psi_1)
\sum_{i=1}^{\ell}x_{\pi(i)}
\ge
\textstyle\sum_{i=1}^{\ell}\lambda_i(M_{\psi_1}(x)), \qquad \mbox{for}~ \ell=1,\ldots,s$;
\item[($iii$)]
$
\lambda_i(M_{\psi_2}(x))
+\psi_2-\psi_1
\ge
\lambda_i(M_{\psi_1}(x)), \qquad \mbox{for}~ i=1,\ldots,s$.
\end{enumerate}
\end{lemma}

The following lemma is a direct consequence of part~($i$) of Lemma~\ref{lem:lambda-shift-s}; together with Corollary~\ref{cor:phi_t_monotone}, it is what allows us to extend the monotonicity argument of \citep*{li2025augmented} from \ref{MESP} to \ref{GMESP}.

\begin{lemma}~\!\!\![Modified shift]
\label{lem:lambda_modified_shift}
Let $C\in\mathbb{S}^n_{++}$. Let $0\!\le\!\psi_1\!\le\!\psi_2\!\le\!\lambda_{\min}(C)$ and $x\in[0,1]^n$ satisfy $\mathbf e^\top x=s$.
Then,
\[
\lambda(M_{\psi_2}(x))
+\psi_2\mathbb I_t
+(\psi_2-\psi_1)(\mathbb I_s-\mathbb I_t)
\succ
\lambda(M_{\psi_1}(x))
+\psi_1\mathbb I_t\,.
\]
\end{lemma}

\begin{proof}
By part~($i$) of Lemma~\ref{lem:lambda-shift-s}, we have $\lambda(M_{\psi_2}(x))+\psi_2\mathbb I_s\succ\lambda(M_{\psi_1}(x))+\psi_1\mathbb I_s$. Subtracting $\psi_1\mathbb I_s$ from both sides preserves majorization, yielding $\lambda(M_{\psi_2}(x))+(\psi_2-\psi_1)\mathbb I_s\succ\lambda(M_{\psi_1}(x))$, as both resulting vectors are nonnegative and sorted in nonincreasing order (the first because $\psi_2\ge\psi_1$), so that, for each $\ell=1,\ldots,n$, the sum of the first $\ell$ components of each vector decreases by the same amount $\psi_1\min\{\ell,s\}$. Using $\mathbb I_s=\mathbb I_t+(\mathbb I_s-\mathbb I_t)$ and adding $\psi_1\mathbb I_t$ to both sides also preserves majorization, as both vectors are sorted in nonincreasing order and $\psi_1\ge0$, so their ordering and nonnegativity are  maintained and, for each $\ell=1,\ldots,n$, the sum of the first $\ell$ components of each vector increases by the same amount $\psi_1\min\{\ell,t\}$. Combining the terms gives $\lambda(M_{\psi_2}(x))+\psi_2\mathbb I_t+(\psi_2-\psi_1)(\mathbb I_s-\mathbb I_t)\succ\lambda(M_{\psi_1}(x))+\psi_1\mathbb I_t$, as claimed.
\end{proof}

We are now in a position to establish the monotonicity of the \ref{DDGFactplus} bound. 

\begin{theorem}~\!\!\![Monotonicity of the \ref{DDGFactplus} bound] \label{thm:ddgfact_plus_monotonicity} Let $C\in\mathbb S_ {++}^n$\,, let $1\le t\le s< n$, and let $0\le\psi_1\le\psi_2\le\lambda_{\min}(C)$. Then, \[  \zddgfactplusthing
    (C,s,t;\psi_2) \le  \zddgfactplusthing
    (C,s,t;\psi_1). \] \end{theorem} 

\begin{proof} Fix an arbitrary $ x\in[0,1]^n$ such that $\mathbf e^\top x=s$. By Lemma~\ref{lem:lambda_modified_shift}, 
\[ 
\lambda(M_{\psi_2}(x)) +\psi_2\mathbb I_t + (\psi_2-\psi_1)(\mathbb I_s-\mathbb I_t) \succ \lambda(M_{\psi_1}(x)) +\psi_1\mathbb I_t\,. 
\] 
As both vectors are sorted in nonincreasing order and nonnegative (see the proof of Lemma~\ref{lem:lambda_modified_shift}), Lemma~\ref{lem:phi_t_majorization} gives 
\begin{equation} \label{eq:hatphi_monotone1} \phi_t\!\left( \lambda(M_{\psi_2}(x)) +\psi_2\mathbb I_t + (\psi_2-\psi_1)(\mathbb I_s-\mathbb I_t) \right) \le \phi_t\!\left( \lambda(M_{\psi_1}(x)) +\psi_1\mathbb I_t \right). \end{equation} 
Moreover, 
$ 
(\psi_2-\psi_1)(\mathbb I_s-\mathbb I_t)\ge0 
$  
componentwise. As $\phi_t$ is componentwise nondecreasing by Corollary~\ref{cor:phi_t_monotone}, \begin{equation} \label{eq:hatphi_monotone2} \phi_t\!\left( \lambda(M_{\psi_2}(x)) +\psi_2\mathbb I_t \right) \le \phi_t\!\left( \lambda(M_{\psi_2}(x)) +\psi_2\mathbb I_t + (\psi_2-\psi_1)(\mathbb I_s-\mathbb I_t) \right). \end{equation} 
Combining \eqref{eq:hatphi_monotone1} and \eqref{eq:hatphi_monotone2} yields 
\[ 
\phi_t\!\left( \lambda(M_{\psi_2}(x)) +\psi_2\mathbb I_t \right) \le \phi_t\!\left( \lambda(M_{\psi_1}(x)) +\psi_1\mathbb I_t \right). 
\] 
Equivalently,
\[ 
\Gamma_t^+(M_{\psi_2}(x);\psi_2) \le \Gamma_t^+(M_{\psi_1}(x);\psi_1). 
\] 
As the above inequality holds for every feasible vector $x$, we  have
\[ 
\begin{aligned}  \zddgfactplusthing
    (C,s,t;\psi_2) &= \max \left\{ \Gamma_t^+(M_{\psi_2}(x);\psi_2): x\in[0,1]^n,\, \mathbf e^\top x=s \right\} \\ &\le \max \left\{ \Gamma_t^+(M_{\psi_1}(x);\psi_1): x\in[0,1]^n,\, \mathbf e^\top x=s \right\} 
\\ &
=  \zddgfactplusthing
    (C,s,t;\psi_1), \end{aligned} 
\] 
which completes the proof. 
\end{proof}

\begin{theorem}~\!\!\![Strict monotonicity of the \ref{DDGFactplus} bound under {{a}} sufficient condition]\label{thm:ddgfact_plus_strict_monotonicity}
    Let $C\in\mathbb S_ {++}^n$\,, let $1\le t\le s< n$, and let $0\le\psi_1<\psi_2\le\lambda_{\min}(C)$. Let $x^*$ be an optimal solution of  \ref{DDGFactplus}  with $\psi:=\psi_2$. Define $\beta^*
        :=
        \lambda(M_{\psi_1}(x^*))
        +
        \psi_1\mathbb I_t $\,.
    Let $\hat\iota$ be the index given by Lemma \ref{lem:index_iota} with respect to $\beta^*$. Let $\pi$ be a permutation of $N$ such that
$
x^*_{\pi(1)}
\ge
x^*_{\pi(2)}
\ge
\cdots
\ge
x^*_{\pi(n)}\,$. If $\hat\iota \geq 1$ and $x^*_{\pi(\hat\iota)} < 1$,
    then $ \zddgfactplusthing
    (C,s,t;\psi_2) <  \zddgfactplusthing
    (C,s,t;\psi_1)$.
\end{theorem}

We note that the hypothesis of Theorem~\ref{thm:ddgfact_plus_strict_monotonicity} is a posteriori: it is expressed in terms of an optimal solution $x^*$ of \ref{DDGFactplus} with $\psi:=\psi_2$. It is, however, immediately checkable once that relaxation has been solved, which is precisely the situation at a node of a B\&B algorithm, where one may then certify that the bound at $\psi_2$ is strictly better than the bound that would have been obtained at a given $\psi_1<\psi_2$ (e.g., $\psi_1=0$, which corresponds to the DDGFact bound).

\begin{proof}
Let $\alpha^*:=\lambda(M_{\psi_2}(x^*))+\psi_2\mathbb I_t$, so that $\zddgfactplusthing(C,s,t;\psi_2)=\phi_t(\alpha^*)$ and $\phi_t(\beta^*)\le\zddgfactplusthing(C,s,t;\psi_1)$, as $x^*$ is feasible. It suffices to show that $\phi_t(\alpha^*)<\phi_t(\beta^*)$. Let $\widetilde\alpha^*:=\alpha^*+(\psi_2-\psi_1)(\mathbb I_s-\mathbb I_t)$ and
\[
\theta^*:=\widetilde\alpha^*-\beta^*=\lambda(M_{\psi_2}(x^*))-\lambda(M_{\psi_1}(x^*))+(\psi_2-\psi_1)\mathbb I_s\,.
\]
By Lemma~\ref{lem:lambda-shift-s}, part~($i$) gives $\sum_{i=1}^n\theta_i^*=0$, part~($iii$) gives $\theta_i^*\ge0$ for $i\le s$, and part~($ii$), with $\ell:=\hat\iota\le t$, gives
\[
\textstyle\sum_{i=1}^{\hat\iota}\theta_i^*\ge(\psi_2-\psi_1)\big(\hat\iota-\sum_{i=1}^{\hat\iota}x^*_{\pi(i)}\big)>0,
\]
because $\hat\iota\ge1$ and $x^*_{\pi(\hat\iota)}<1$.

Let $g$ be the vector of Lemma~\ref{lem:phi_t_subgrad} at $y:=\beta^*$; the lemma applies because $\beta^*_t>0$ (this is clear if $\psi_1>0$, and follows from the discussion after \eqref{DDGFactplus} if $\psi_1=0$). Then $g_1\le\cdots\le g_{\hat\iota}<g_{\hat\iota+1}=\cdots=g_n$\,, where the strict inequality is $1/\beta^*_{\hat\iota}<(t-\hat\iota)/\sum_{j>\hat\iota}\beta^*_j$, by Lemma~\ref{lem:index_iota}. As $\widetilde\alpha^*$ is nonnegative and sorted in nonincreasing order (see the proof of Lemma~\ref{lem:lambda_modified_shift}), Lemma~\ref{lem:phi_t_subgrad} and the three facts above give
\[
\textstyle\phi_t(\widetilde\alpha^*)-\phi_t(\beta^*)\le g^\top\theta^*
=\sum_{i=1}^{\hat\iota}(g_i-g_{\hat\iota+1})\theta_i^*
\le(g_{\hat\iota}-g_{\hat\iota+1})\sum_{i=1}^{\hat\iota}\theta_i^*<0.
\]
Finally, as $\alpha^*\le\widetilde\alpha^*$, Corollary~\ref{cor:phi_t_monotone} gives $\phi_t(\alpha^*)\le\phi_t(\widetilde\alpha^*)<\phi_t(\beta^*)$.
\end{proof}

\section{Constructing  dual solutions}
\label{subsec:feasible_dual_sol_dgfactplus}

Lagrangian duality provides not only computable upper bounds but also a
systematic way of constructing dual-feasible solutions from
primal-feasible ones. In a B\&B framework, such dual-feasible
solutions serve two important purposes. First, they provide rigorous upper
bounds even when the relaxation is not solved to optimality. Second,  they enable variable
fixing (see Section~\ref{sec:computational_results} and Theorem~\ref{thm:variable_fixing_dgfact_plus} in the Appendix). For the classical factorization bound for
\ref{MESP}, \citep*[Sec.~3.4.4]{FLbook} introduced such a primal--dual
construction, and \citep*[Thm.~21]{GMESP_Alg} extended it to the
generalized factorization bound for \ref{GMESP}, showing that an optimal
primal solution yields an optimal dual solution.

In this section, we extend these results to the augmented generalized
factorization bound in the presence of fixed variables, a setting that
naturally arises in B\&B algorithms. In a B\&B algorithm for \ref{GMESP}, each node corresponds to a subproblem in which some variables are fixed at zero or one. We represent such a subproblem by vectors $l,c\in\{0,1\}^n$ with $l\le c$: $x_i$ is fixed at one if $l_i=1$, and at zero if $c_i=0$. The root node corresponds to the original problem, with no variables fixed, that is, $l=\mathbf{0}$ and $c=\mathbf{e}$; its relaxation is called the root relaxation. We first introduce the
corresponding fixed-variable primal and dual formulations. Following the approach of  \citep*{GMESP_Alg},
we  show
how to construct a feasible solution of the dual problem from any feasible
solution of the primal problem, we show that the vector $\diag(F\hat\Theta F^\top)$, computed from the constructed dual solution, is a supergradient of the objective of DDGFact$^+$ and gives its directional derivatives, which enables us to prove that this
construction preserves optimality.

\subsection{Fixed-variable formulations}

\begin{definition}
For $l,c\in\{0,1\}^n$ satisfying $l\le c$, let
DDGFact$^+(l,c)$
denote formulation \ref{DDGFactplus} with the additional constraints
$
l\le x\le c$. Specifically, 
    \begin{align}
    \textstyle\max\left\{   \Gamma_t^+(M_\psi(x); \psi)~:~
         \mathbf{e}^\top x = s,~ l \leq x \leq c\right\}.
    \tag{DDGFact$^{+}(l,c)$}\label{DDGFactpluslc}
\end{align}
\end{definition}
Following the development of Subsections \ref{subsec:dual_gfactplus} and \ref{subsec:dual_dgfactplus}, we can verify  that \ref{DDGFactpluslc} is the Lagrangian dual problem of 
\begin{equation}
    \begin{array}{lll}
    &\min&  ~~\sum_{\ell=k-t+1}^{k}h_\psi(\lambda_{\ell}(\Theta))
    -\upsilon^\top l+ \nu^\top c +\tau s\\[3pt]
    &\mbox{s.t.}&
         ~~\diag(F \Theta F^\top) + \upsilon - \nu  - \tau\mathbf{e}=0,\\[1pt]
    && ~~\Theta\in \mathbb{S}^k_{++}\,, ~\upsilon\in\mathbb{R}^n_+\,, ~\nu\in\mathbb{R}^n_+\,,~\tau\in\mathbb{R}.
    \end{array}
    \tag{DGFact$^{+}(l,c)$}\label{DGFactpluslc}
\end{equation}

\subsection{Constructing feasible dual solutions}

We begin by showing that every feasible solution of the fixed-variable
primal problem naturally induces a feasible solution of the corresponding
dual problem.

\begin{proposition}~\!\!\![Constructed solution for \ref{DGFactpluslc}]
\label{prop:dgfactplus_associated_solution_fixed}
Let $l,c\in\{0,1\}^n$ satisfy $l\le c$, let $\psi\ge0$, and let
$\hat x$ be a feasible solution of \ref{DDGFactpluslc}.
Let $F^\top\Diag(\hat x)F=Q\Diag(\lambda)Q^\top$ be a spectral decomposition
with $\lambda_1\ge\cdots\ge\lambda_k\ge0$, and let \mbox{$\hat w\in\mathbb{R}^k$} be
the vector associated with $\lambda$ in Theorem~\ref{thm:subproblem_Theta}.
Let
$
\hat\Theta:=Q\Diag(\hat w)Q^\top$ and $\hat d:=\diag(F\hat\Theta F^\top)$.
Let $\pi$ be a permutation of $N$ such that
$
        \hat d_{\pi(1)}
        \geq
        \hat d_{\pi(2)}
        \geq
        \cdots
        \geq
        \hat d_{\pi(n)}\,.
$ If
$
c_{\pi(1)}
+
\textstyle \sum_{i=2}^{n}l_{\pi(i)}
>
s$,
let $q:=0$.
Otherwise, let
$
q
:=
\textstyle \max
\left\{
j\in N:
\sum_{i=1}^{j}c_{\pi(i)}
+
\sum_{i=j+1}^{n}l_{\pi(i)}
\le s
\right\}$,
and
\[
\hat\tau:=
\begin{cases}
\hat d_{\pi(q+1)}\,,
&
q<n;
\\
0,
&
q=n.
\end{cases}
\]
Let $
P
:=
\{\pi(1),\ldots,\pi(q)\}$ and $
R
:=
\{\pi(q+2),\ldots,\pi(n)\}$. Finally, define
\[
\hat\nu_i:=
\begin{cases}
\hat d_i-\hat\tau,
&i\in P;
\\
0,
&i\notin P,
\end{cases}
\qquad
\hat\upsilon_i:=
\begin{cases}
\hat\tau-\hat d_i\,{,}
&i\in R;
\\
0,
&i\notin R.
\end{cases}
\]
Then
$(\hat\Theta,\hat\upsilon,\hat\nu,\hat\tau)$
is feasible for 
\ref{DGFactpluslc}.

\end{proposition}

\begin{proof}
Note that $\hat\nu=(\hat d - \hat\tau\mathbf{e})_+$ and $\hat\upsilon=(\hat\tau\mathbf{e}-\hat d)_+$\,. It follows that  $\hat d +\hat\upsilon-\hat\nu-\hat\tau\mathbf{e}=0$ and, therefore, $(\hat\Theta,\hat\upsilon,\hat\nu,\hat\tau)$  is feasible for \ref{DGFactpluslc}.
\end{proof}

The constructed dual solution also yields a supergradient of the objective function of \ref{DDGFactpluslc}.

\begin{lemma}~\!\!\![Supergradient of the \ref{DDGFactplus} objective]
\label{lem:gamma_t_supergradient}
Let $l$, $c$, $\psi$, $\hat x$, $\hat\Theta$ and $\hat d$ be as in
Proposition~\ref{prop:dgfactplus_associated_solution_fixed}. Then $\hat d$ is a
supergradient of the concave function $x\mapsto\Gamma_t^+(M_\psi(x);\psi)$ at
$\hat x$; that is,
\[
\Gamma_t^+(M_\psi(x);\psi)\le\Gamma_t^+(M_\psi(\hat x);\psi)+\hat d^\top(x-\hat x),
\qquad\text{for all } x\in[0,1]^n.
\]
\end{lemma}

\begin{proof}
By \eqref{eq:inf_Theta_to_w} and Theorem~\ref{thm:subproblem_Theta}, for every
$x\in[0,1]^n$,
\[
\Gamma_t^+(M_\psi(x);\psi)
=
\inf_{\Theta\in\mathbb{S}^k_{++}}
\Big\{\textstyle\sum_{\ell=k-t+1}^{k}h_\psi(\lambda_\ell(\Theta))
+\big(F^\top\Diag(x)F\big)\bullet\Theta\Big\},
\]
(If $\psi=0$ and $\lambda_t(M_0(x))=0$, Theorem~\ref{thm:subproblem_Theta} does not apply, but this equality still holds, as both sides equal $-\infty$. Indeed, the left-hand side equals $-\infty$ by Lemma~\ref{lem:phi_t_subgrad}. For the right-hand side, let $r:=\rank(M_0(x))<t$, let $Q$ be an orthogonal matrix with $M_0(x)=Q\Diag(\lambda(M_0(x)))Q^\top$, and, for $\mu\ge1$, let $\Theta_\mu:=Q\Diag(\mathbb{I}_r+\mu(\mathbf{e}-\mathbb{I}_r))Q^\top$. Then $M_0(x)\bullet\Theta_\mu=\sum_{\ell=1}^r\lambda_\ell(M_0(x))$ does not depend on $\mu$, and the $t$ least eigenvalues of $\Theta_\mu$ consist of $r$ eigenvalues equal to $1$ and $t-r\ge1$ eigenvalues equal to $\mu$, so that the objective inside the braces equals $-r+(t-r)(-\log\mu-1)+\sum_{\ell=1}^r\lambda_\ell(M_0(x))$, which tends to $-\infty$ as $\mu\to\infty$.)
and the infimum is attained at $\hat\Theta$ when $x=\hat x$. 
Indeed, as $\hat w$ is nondecreasing while $\lambda$ is nonincreasing, $\hat\Theta$ and $F^\top\Diag(\hat x)F$ share the eigenbasis $Q$, with the least eigenvalue of $\hat\Theta$ paired with the greatest eigenvalue of $F^\top\Diag(\hat x)F$, which is the reverse pairing required by Lemma~\ref{lem:von_neumann}. Hence, the attainment clause of Lemma~\ref{lem:von_neumann} and Theorem~\ref{thm:subproblem_Theta} show that the objective inside the braces at $\Theta=\hat\Theta$ equals $\Gamma_t^+(M_\psi(\hat x);\psi)$.
For each fixed
$\Theta$ the function inside the braces is affine in $x$, with gradient
$\diag(F\Theta F^\top)$, because
$(F^\top\Diag(x)F)\bullet\Theta=\sum_{i\in N}x_i\,F[i,\cdot]\Theta F[i,\cdot]^\top$.
Hence, for every $x\in[0,1]^n$,
\[
\Gamma_t^+(M_\psi(x);\psi)
\le
\textstyle\sum_{\ell=k-t+1}^{k}h_\psi(\lambda_\ell(\hat\Theta))
+\big(F^\top\Diag(x)F\big)\bullet\hat\Theta
=
\Gamma_t^+(M_\psi(\hat x);\psi)+\hat d^\top(x-\hat x),
\]
which is the supergradient inequality. Note, moreover, that $\hat w$ is exactly
the supergradient $g$ of $\phi_t$ at $\lambda+\psi\mathbb{I}_t$ given by
Lemma~\ref{lem:phi_t_subgrad}, so that
$\hat d_i=F[i,\cdot]\,Q\Diag(g)\,Q^\top F[i,\cdot]^\top$, for $i\in N$; see the form of the subdifferentials of convex spectral functions in \citep*[Lem.~10]{Nikolov}.
\end{proof}

{
\begin{lemma}~\!\!\![Directional derivatives of the \ref{DDGFactplus} objective]
\label{lem:gamma_t_directional}
(proof in Appendix)
Let $l$, $c$, $\psi$, $\hat x$, $\hat\Theta$ and $\hat d$ be as in
Proposition~\ref{prop:dgfactplus_associated_solution_fixed}, and let
$f(x):=\Gamma_t^+(M_\psi(x);\psi)$, for $x\in\mathbb{R}^n_+$. Then
\[
\lim_{x\to\hat x,\;x\in\mathbb{R}^n_+\setminus\{\hat x\}}
\frac{f(x)-f(\hat x)-\hat d^\top(x-\hat x)}{\|x-\hat x\|}=0,
\]
where $\|\cdot\|$ denotes the Euclidean norm. In particular, $f'(\hat x;v)=\hat d^\top v$ for every $v\in\mathbb{R}^n$ such that $\hat x+\delta v\in\mathbb{R}^n_+$ for all sufficiently small $\delta>0$.
\end{lemma}
For $\psi=0$, Lemma~\ref{lem:gamma_t_directional} also supplies a self-contained justification of the directional-derivative property used in the proof of \citep*[Lem.~20]{GMESP_Alg}, where it is taken from \citep*[Prop.~11]{gscale}.
}

\subsection{Optimal associated dual solutions} 

To prove that the constructed  dual solution presented in the previous subsection is optimal whenever the primal
solution is optimal, we adapt the proof of
\citep*[Thm.~21]{GMESP_Alg} to the augmented generalized
factorization bound. The key ingredient is Lemma~\ref{lem:gamma_t_directional}, which shows that $\hat d$ gives the directional derivatives of the primal objective at $\hat x$.
The proof follows the same
principles as that of \citep*[Thm.~21]{GMESP_Alg} and is included in the Appendix for
completeness.

\begin{theorem}~\!\!\![Optimal associated DGFact$^+$ solution]
\label{thm:dgfactplus_associated_solution_optimal_fixed}
(proof in Appendix) 
Let $l,c\in\{0,1\}^n$ satisfy $l\le c$, let $\psi\ge0$, and let $\hat x$ be an optimal solution of \ref{DDGFactpluslc}.
Construct
$(\hat\Theta,\hat\upsilon,\hat\nu,\hat\tau)$
from $\hat x$ as in
Proposition~\ref{prop:dgfactplus_associated_solution_fixed}.
Then
$(\hat\Theta,\hat\upsilon,\hat\nu,\hat\tau)$
is an optimal solution of 
\ref{DGFactpluslc}.

\end{theorem}

\section{Closed-form solutions for special cases}
\label{sec:special_case_t=1}

In this section, we derive a closed-form solution for  \ref{DDGFactplus} under easily verifiable sufficient conditions. The main ingredient is a trace upper bound for the function $\Gamma_t^+$. Combined with the observation that maximizing the trace amounts to selecting the $s$ greatest diagonal entries of $C$, this bound yields a sufficient condition under which \ref{DDGFactplus} admits a closed-form optimal solution. We then extend this result to B\&B subproblems. 
When $t=1$, the sufficient condition is automatically satisfied, providing a closed-form solution for the relaxation at the root  and every B\&B subproblem.
Finally, we show that, under an additional eigenvalue multiplicity condition, the closed-form solution of \ref{DDGFactplus} is also optimal for  \ref{GMESP}. 

\begin{lemma}~\!\!\![Trace upper bound for $\Gamma_t^+$]
\label{lem:gamma_trace_upper_bound}
    Let $X\in\mathbb{S}_+^k$,  $\psi\ge 0$, and  $0< t\leq k$. Then
    \begin{equation*}
        \textstyle\Gamma_t^+(X;\psi)
        \leq
        t\log\left(
            \frac{\tr(X)}{t}+\psi
        \right).
    \end{equation*}
    Moreover, equality holds whenever  $\tr(X)\geq t\lambda_{\max}(X)$.
\end{lemma}

\begin{proof}
    Let $y \in \mathbb{R}^k_+$ be
    \begin{equation*}
       \textstyle y_i :=
        \begin{cases}
            \lambda_i(X)+\psi, & i=1,\dots,t;\\
            \lambda_i(X), & i=t+1,\dots,k.
        \end{cases}
    \end{equation*}
    Then $y_1\geq \cdots\geq y_k\geq 0$ and
\begin{equation}\label{trace}      \textstyle\textstyle\frac{\sum_{i=1}^k y_i}{t}
        =
        \frac{\tr(X)}{t}+\psi.
\end{equation}    
Let $\hat\iota$ be the unique index corresponding to $y$ that satisfies Lemma~\ref{lem:index_iota} and let $\bar y:=
\textstyle\frac{1}{t-\hat\iota}
\sum_{j=\hat\iota+1}^k {y}_{j}$\,.
    Then
    \begin{equation*}
        \Gamma_t^+(X;\psi)
        =
        \log\left(
            \left(\textstyle\prod_{i=1}^{\hat\iota} y_i\right)
            \bar{y}^{t-\hat\iota}
        \right).
    \end{equation*}   
    Applying the arithmetic--geometric mean inequality to the $t$ numbers
$y_1,\ldots,y_{\hat\iota},\bar y,\ldots,\bar y$ gives
    \begin{equation*}
        \textstyle\left(\prod_{i=1}^{\hat\iota} y_i\right)
        \bar{y}^{t-\hat\iota}
        \leq
        \left(
            \frac{
                \sum_{i=1}^{\hat\iota}y_i
                +
                (t-\hat\iota)\bar{y}
            }{t}
        \right)^t.
    \end{equation*}
   As $(t-\hat\iota)\bar y = \sum_{i=\hat\iota+1}^k y_i$\,, it follows that 
    \begin{align*}
        \Gamma_t^+(X;\psi)
        \leq      
        t\log\left(
            \textstyle\frac{
                \sum_{i=1}^{k}y_i}{t}\right) 
=
        t\log\left(
            \textstyle\frac{\tr(X)}{t}+\psi
        \right).
    \end{align*}
    Moreover, by  \eqref{trace} and Lemma~\ref{lem:index_iota}, we have 
    \begin{equation*}
        \textstyle\frac{\tr(X)}{t}
        \geq
        \textstyle\lambda_{\max}(X) \quad\Longleftrightarrow\quad\textstyle\frac{1}{t}\sum_{i=1}^k y_i \geq y_1 \quad\Longleftrightarrow\quad \hat\iota = 0.
    \end{equation*}
Finally,  when $\hat\iota = 0$, we have 
    \begin{equation*}
        \textstyle \Gamma_t^+(X;\psi)
        =
        t\log\left(
            \frac{\sum_{i=1}^k y_i}{t}
        \right)
        =
        t\log\left(
            \frac{\tr(X)}{t}+\psi
        \right).
    \end{equation*}
    This completes the proof.
\end{proof}

The trace upper bound of
Lemma~\ref{lem:gamma_trace_upper_bound}
is a key ingredient to the derivation of the following result.

\begin{theorem}~\!\!\![Closed-form solution  for \ref{DDGFactplus} under sufficient conditions]
\label{thm:ddgfact_plus_closed_form_iota_zero_certificate}
Let $C\in\mathbb S_{++}^n$\,, let $1\le t\le s< n$, and let
$0\le\psi\le\lambda_{\min}(C)$.
Let $\bar{S}\subset N$ be any subset of cardinality $s$
corresponding to the $s$ greatest diagonal entries of $C$, and let
$\bar{x}$ be its characteristic vector.
Suppose that
\[
\textstyle\tr((C-\psi I)[\bar{S},\bar{S}])
\ge
t\,\lambda_{\max}\!\left((C-\psi I)[\bar{S},\bar{S}]\right).
\]
Then $\bar{x}$ is an optimal solution of \ref{DDGFactplus}, with objective value
\[
\textstyle t\log\!\left(
\frac1t\tr((C-\psi I)[\bar{S},\bar{S}])
+\psi
\right).
\]
\end{theorem}

\begin{proof}
Let
$
C-\psi I=FF^\top
$
be any factorization. As
\[
\textstyle \tr(M_\psi(x))
=
\sum_{i=1}^n x_i(C_{ii}-\psi),
\]
and every feasible solution of \ref{DDGFactplus} satisfies $\mathbf e^\top x=s$, the vector $\bar{x}$
maximizes $\tr(M_\psi(x))$ over the feasible region. Hence,
\[
\tr(M_\psi(x))
\le
\tr(M_\psi(\bar{x}))
\]
for every feasible $x$.
Then, applying Lemma~\ref{lem:gamma_trace_upper_bound} and the monotonicity of the log function gives
\[
\textstyle \Gamma_t^+(M_\psi(x);\psi)
\le
t\log\!\left(
\frac{\tr(M_\psi(x))}{t}+\psi
\right)
\le
t\log\!\left(
\frac{\tr(M_\psi(\bar{x}))}{t}+\psi
\right).
\]
Moreover,
\[
M_\psi(\bar{x})
=
F[\bar{S},\cdot]^{\top}F[\bar{S},\cdot],
\qquad
(C-\psi I)[\bar{S},\bar{S}]
=
F[\bar{S},\cdot]F[\bar{S},\cdot]^{\top},
\]
and these matrices have the same nonzero eigenvalues. Therefore, the assumed
condition is equivalent to
\[
\tr(M_\psi(\bar{x}))
\ge
t\,\lambda_{\max}(M_\psi(\bar{x})).
\]
Consequently, Lemma~\ref{lem:gamma_trace_upper_bound}  yields
\[
\textstyle \Gamma_t^+(M_\psi(\bar{x});\psi)
=
t\log\!\left(
\frac{\tr(M_\psi(\bar{x}))}{t}+\psi
\right),
\]
which proves that
\[
\Gamma_t^+(M_\psi(x);\psi)
\le
\Gamma_t^+(M_\psi(\bar{x});\psi)
\]
for every feasible $x$. Thus, $\bar{x}$ is optimal.
Finally,
\[
\textstyle \tr(M_\psi(\bar{x}))
=
\tr((C-\psi I )[\bar{S},\bar{S}]),
\]
and therefore
\[
\Gamma_t^+(M_\psi(\bar{x});\psi)
=
\textstyle t\log\!\left(
\frac1t\tr((C-\psi I )[\bar{S},\bar{S}])
+\psi
\right),
\]
which is the claimed objective value.
\end{proof}
The previous theorem applies to the root relaxation.
The same argument immediately extends to B\&B subproblems. In the following, without loss of generality, we only consider variables fixed at one in the subproblems (variables fixed at zero are easily handled by eliminating rows/columns of $C$).

\begin{corollary}~\!\!\![Closed-form solution  for \ref{DDGFactpluslc}  with variables fixed at one under sufficient conditions]
\label{cor:ddgfact_plus_closed_form_iota_zero_node}
Let $C\in\mathbb S_ {++}^n$\,,  $1\le t\le s< n$, and  $0\le\psi\le\lambda_{\min}(C)$. Let
$J\subset N$ satisfy $|J|\le s$, and let $\bar{S}\subset N\setminus J$
be any subset of cardinality
$s-|J|$
corresponding to the
$s-|J|$
greatest diagonal elements of $C$
over
$N\setminus J$. Let $l$ be the characteristic vector of $J$ and $c:=\mathbf{e}$.
Let $\bar{x}$ be the characteristic vector of $J\cup \bar{S}$. 
    Suppose that
    \begin{equation*}
       \textstyle \tr((C-\psi I)[ J\cup \bar{S}, J\cup \bar{S}])
        \geq
        t\,\lambda_{\max}\left(
             (C-\psi I)[J\cup \bar{S},J\cup \bar{S}]
        \right).
    \end{equation*}
    Then $\bar{x}$ is an optimal solution of \ref{DDGFactpluslc} with objective value 
    \begin{equation*}
        t\log\left(
            \textstyle\frac{1}{t}\tr((C-\psi I)[ J\cup \bar{S}, J\cup \bar{S}]) +\psi
        \right).
    \end{equation*}
\end{corollary}

\begin{proof}
    The proof follows exactly as in Theorem~\ref{thm:ddgfact_plus_closed_form_iota_zero_certificate}, replacing the trace maximizer over $N$ by the trace maximizer over $N\setminus J$. As the contribution of the fixed variables is constant, maximizing the trace reduces to selecting the $s-|J|$ greatest diagonal elements of $C$  among the remaining indices. The assumed inequality guarantees that equality holds in the trace upper bound, yielding optimality.
\end{proof}

\begin{remark}
    For the case where $t=1$, the assumptions in Theorem \ref{thm:ddgfact_plus_closed_form_iota_zero_certificate} and Corollary \ref{cor:ddgfact_plus_closed_form_iota_zero_node} always hold, and therefore, we obtain a closed-form solution for \ref{DDGFactplus} at the root node and for every subproblem as well.
\end{remark}

Theorem
\ref{thm:ddgfact_plus_closed_form_iota_zero_certificate} 
shows that under a particular assumption \ref{DDGFactplus}  has a closed-form solution.
A natural question is whether the resulting binary solution is also optimal for the
original problem \ref{GMESP}. This matters because, when the relaxation is exact at a binary optimal solution, the upper bound coincides with the objective value of a feasible solution of \ref{GMESP}, which is therefore certified as optimal. In B\&B, the corresponding node can be pruned after updating the incumbent, without branching. The following corollary provides additional conditions under which this optimality is certified.

\begin{corollary}~\!\!\![Closed-form solution of \ref{GMESP} under sufficient conditions]
\label{cor:optimal_gmesp_closed_form}
Assume that the hypothesis of Theorem~\ref{thm:ddgfact_plus_closed_form_iota_zero_certificate} holds with $\psi:=\lambda_{\min}(C)$, and let $x^*$ be the closed-form optimal solution of \ref{DDGFactplus} given by that theorem. If, in addition, the multiplicity $\mu_{\min}$ of the least eigenvalue of $C$ satisfies $\mu_{\min}\ge n-t$, then $x^*$ is also an optimal solution of \ref{GMESP}.
\end{corollary}

\begin{proof}
    By Theorem~\ref{thm:ddgfact_plus_closed_form_iota_zero_certificate}, $x^*$ is
    an optimal solution of \ref{DDGFactplus}. As $\mu_{\min}\ge n-t$ and
    $\psi=\lambda_{\min}(C)$, Theorem~\ref{thm:gamma_t_collapse_multiplicity}
    gives
    $\Gamma_t^+(M_\psi(x^*);\psi)=f_t(M_\psi(x^*);\psi)$, and, as $x^*$ is
    binary, Theorem~\ref{Gmesp2GFact} shows that the latter is the objective
    value of $x^*$ for \ref{GMESP}. Finally, \ref{DDGFactplus} is a relaxation of
    \ref{GMESP}, so its optimal value is an upper bound on $z(C,s,t)$; this upper
    bound is attained by the feasible solution $x^*$, which is therefore optimal
    for \ref{GMESP}.
\end{proof}

\begin{remark}\label{rmk:closed_form_condition}
    The two hypotheses of
    Corollary~\ref{cor:optimal_gmesp_closed_form} pull in opposite directions in
    $\psi$ when $t\ge2${ and $s>t$}. Indeed, as
    $\tr((C-\psi I)[\bar S,\bar S])=\tr(C[\bar S,\bar S])-s\psi$ and
    $\lambda_{\max}((C-\psi I)[\bar S,\bar S])=\lambda_{\max}(C[\bar S,\bar S])-\psi$,
    the hypothesis of
    Theorem~\ref{thm:ddgfact_plus_closed_form_iota_zero_certificate} is
    equivalent to
    \begin{equation}\label{eq:trace_condition_lmin}
        \tr\!\left(C[\bar S,\bar S]\right)
        -
        t\,\lambda_{\max}\!\left(C[\bar S,\bar S]\right)
        \ \ge\
        (s-t)\,\psi,
    \end{equation}
    whose left-hand side does not depend on $\psi$, while its right-hand side is
    increasing in $\psi$; and the corollary requires
    $\psi=\lambda_{\min}(C)$, the greatest admissible value. Moreover, under the multiplicity hypothesis, $\rank((C-\psi I)[\bar S,\bar S])\le\rank(C-\psi I)\le t$, so that $\tr((C-\psi I)[\bar S,\bar S])$ is the sum of the $t$ greatest eigenvalues of $(C-\psi I)[\bar S,\bar S]$; hence, the hypothesis of Theorem~\ref{thm:ddgfact_plus_closed_form_iota_zero_certificate} holds if and only if these $t$ eigenvalues are all equal, which is restrictive when $t\ge2$. For $t=1$, however,
    \eqref{eq:trace_condition_lmin} holds automatically, as
    $\lambda_\ell(C[\bar S,\bar S])\ge\lambda_{\min}(C)$, for $\ell=1,\ldots,s$,
    by eigenvalue interlacing.
\end{remark}

\section{Comparison with the spectral bound}\label{sec:comparison_spectral}

In this section, we investigate the relationship between the augmented generalized factorization bound and the spectral bound. We first characterize this relationship exactly when $t=1$, obtaining a necessary and sufficient condition for one bound to dominate the other. We  show that, in general, the two bounds are incomparable, although the augmentation parameter can significantly change their ordering. Finally, for general $t$, we derive worst-case comparison guarantees and sufficient conditions under which either bound provably dominates the other.

Throughout, for two upper bounds on $z(C,s,t)$, we say that one dominates the other if it is never greater, so that it is at least as tight. We begin by recalling the definition of the spectral bound.

\begin{definition}\label{def:spectral_bound}
The \emph{spectral bound} for \ref{GMESP} (see
\citep*{LeeLind2020}; see also \citep*{KLQ}) is
$
    \zspectralthing(C,t)
    :=
    \textstyle\sum_{\ell=1}^{t}\log\bigl(\lambda_\ell(C)\bigr)$.
\end{definition}

\subsection{Exact comparison for \texorpdfstring{$t=1$}{t=1}}\label{subsec:exact_comparison_t1}

In this subsection, we completely characterize the relationship between the \ref{DDGFactplus} bound and the spectral bound  when $t=1$. In particular, we  see that whether
the \ref{DDGFactplus} bound dominates the spectral bound depends only on the sum of the
$s$ greatest diagonal entries of $C$, the greatest eigenvalue of $C$, and the
augmentation parameter $\psi$. 

\begin{proposition}~\!\!\![Exact comparison for $t=1$]
\label{prop:ddgfactplus_t1_exact_spectral_comparison}
Let $C\in\mathbb S_ {++}^n$\,,  $1\le s< n$, $t=1$, and  $0\le\psi\le\lambda_{\min}(C)$. 
Let $\bar{S}\subset N$  be any subset of cardinality $s$
of the indices of the $s$ greatest diagonal elements of $C$. Then,
\[
\zddgfactplusthing(C,s,1;\psi)
-
\zspectralthing(C,1)
=
\log \left(
\textstyle\frac{
\tr((C-\psi I)[\bar S,\bar S]) +\psi
}
{\lambda_{\max}(C)}
\right).
\]
\end{proposition}

\begin{proof}
We have that $\zspectralthing(C,1)=\log(\lambda_{\max}(C))$. Then the result directly follows from  Theorem~\ref{thm:ddgfact_plus_closed_form_iota_zero_certificate}, which gives 
$
\textstyle \zddgfactplusthing(C,s,1;\psi)
=
\textstyle \log\!\left(
\tr((C-\psi I)[\bar{S},\bar{S}])
+\psi
\right).
$
\end{proof}

From Proposition \ref{prop:ddgfactplus_t1_exact_spectral_comparison}, we see that the quantity obtained by subtracting the spectral bound from the \ref{DDGFactplus} bound is  nondecreasing in $s$ (a set-wise larger subset $\bar{S}$ yields a larger trace of the corresponding submatrix of $C-\psi I$). Exploiting this monotonicity, the following example shows that, in general, neither the  \ref{DDGFactplus} bound  nor the
spectral bound dominates the other. 

\begin{example}
\label{ex:ddgfactplus_t1_nondominance}
Let $t=1$ and $C :=
\left(\begin{smallmatrix}
4 & 2 & 2 & 1\\
2 & 4 & 1 & 2\\
2 & 1 & 4 & 2\\
1 & 2 & 2 & 4
\end{smallmatrix}\right)\in \mathbb{S}^4_{++}$\,.
We can verify that
$\lambda_{\min}(C) = 1$ and $\lambda_{\max}(C) = 9$.
Then, from Proposition \ref{prop:ddgfactplus_t1_exact_spectral_comparison}, we have
$
\zddgfactplusthing(C,s,1;\psi) - \zspectralthing(C,1) =
\log\!\left(\tfrac{4s-(s-1)\psi}{9}\right)$, for $1\le s< 4$ and  $0\le\psi\le\lambda_{\min}(C)$.
Thus, for every
$\psi \in [0,\lambda_{\min}(C)]$, the  \ref{DDGFactplus} bound dominates the spectral bound when $s=1,2$, whereas the reverse holds when $s=3$.\hfill $\clubsuit$
\end{example}

Finally, the next example demonstrates that the introduction of the augmentation 
parameter can change the ordering between the factorization and the spectral 
bounds. 

\begin{example}
Let
\[
C:=(1-\rho)I_n+\rho\,\mathbf e\mathbf e^\top,
\qquad
0<\rho<\textstyle\frac{s-1}{n-1},
\]
with $1< s<n$.
Then
$
\lambda_{\max}(C)=1+(n-1)\rho$ and 
$\lambda_2(C)=\cdots=\lambda_n(C) (=\lambda_{\min}(C)) =1-\rho$.
Choosing $\psi=\lambda_{\min}(C)$,
Proposition~\ref{prop:ddgfactplus_t1_exact_spectral_comparison}
gives
$
\zddgfactplusthing(C,s,1;\lambda_{\min}(C))
=
\log\!\left(1+(s-1)\rho\right),
$
whereas
$
\zspectralthing(C,1)
=
\log\!\left(1+(n-1)\rho\right).
$
As $s<n$, we have
$
\zddgfactplusthing(C,s,1;\lambda_{\min}(C))
<
\zspectralthing(C,1).
$

In contrast, when $\psi=0$ (the original DGFact bound),
we have $
\zddgfactplusthing(C,s,1;0)=\log(s),
$
and the assumption
$
\rho<\frac{s-1}{n-1}
$
implies that 
$
1+(n-1)\rho<s,
$
so that
$
\zspectralthing(C,1)
<
\zddgfactplusthing(C,s,1;0).
$

Therefore,
$
\zddgfactplusthing(C,s,1;\lambda_{\min}(C))
<
\zspectralthing(C,1)
<
\zddgfactplusthing(C,s,1;0).
$
\hfill $\clubsuit$
\end{example}

\subsection{Worst-case comparison}

In this subsection, we investigate how far the
augmented bound can deviate from the spectral bound in the worst case.
For the original DGFact bound,
\citep*[Thm.~7]{GMESP_Alg} established that
$
\zdgfactplusthing(C,s,t;0)
-
\zspectralthing(C,t)
\le
t\log(s/t),
$
which provides a worst-case guarantee that is independent of $C$.

Our first result derives the analogous estimate for the augmented
formulation. A direct adaptation of the proof of
\citep*[Thm.~7]{GMESP_Alg}
yields an upper bound containing an additional correction term depending on
the augmentation parameter $\psi$. 
The main contribution of this subsection is to show that this additional
term can be eliminated by constructing a different dual feasible
solution. This leads to a better upper bound, which not only
recovers the classical worst-case guarantee $t\log(s/t)$ for
\ref{DGFactplus}, but also strengthens it by exhibiting a nonnegative
correction term that strictly improves the estimate whenever $\psi>0$  and $s>t$.

The proof of the preliminary comparison theorem is presented in
the Appendix, as it closely follows the arguments used for the
original DGFact bound.

\begin{theorem}~\!\!\![A preliminary comparison] 
\label{thm:dgfactplus_spectral_comparison}
(proof in Appendix) 
Let $C\in\mathbb S_ {++}^n$\,, 
$0<\psi\le\lambda_{\min}(C)$, and 
$1\le t\le s<n$.
Let $C-\psi I=FF^\top$.
Then
\[
\zdgfactplusthing(C,s,t;\psi)
-
\zspectralthing(C,t)
\le
t\log(s/t)
+
\textstyle\frac{t\psi}{s}
\sum_{\ell=1}^{t}\frac1{\lambda_\ell(C)}.
\]
\end{theorem}

In Theorem~\ref{thm:dgfactplus_improved_spectral_comparison}, we strengthen the above result by recovering the same worst-case comparison
with the spectral bound as the one established for the original DGFact bound. The proof relies on a consequence of the Schur--Horn theorem.

\begin{theorem}~\!\!\![An improved comparison]
\label{thm:dgfactplus_improved_spectral_comparison}
Let $C\in\mathbb S_{++}^{n}$\,, 
$0<\psi\le\lambda_{\min}(C)$,
and 
$1\le t\le s<n$. Let $C-\psi I=FF^\top$. 
Then
\[
\zdgfactplusthing(C,s,t;\psi)
-
\zspectralthing(C,t)
\le
\textstyle t\log(s/t)
-
\frac{t\psi}{s}
\sum_{\ell=t+1}^{s}
\textstyle\frac1{\lambda_\ell(C)}.
\]
\end{theorem}

\begin{proof}
By Theorem~\ref{thm:dgfactplus_factorization_invariance}, we may assume
$
F=(C-\psi I)^{1/2}$.
Let
$
C=\textstyle\sum_{\ell=1}^{n}\lambda_\ell(C)u_\ell u_\ell^\top$ with 
$\lambda_1(C)\ge\cdots\ge\lambda_n(C)>0$,
and define
$
\textstyle\alpha:=\frac{t}{s}$ and $\hat\Theta:=\alpha C^{-1}$.

Unlike the proof of
Theorem~\ref{thm:dgfactplus_spectral_comparison},
we construct the remaining dual variables so as to exploit the
Schur--Horn majorization theorem.

Let
$
\hat d:=\diag(F\hat\Theta F^\top)$. Let $\pi$ be a permutation of $N$ such that
$
        \hat d_{\pi(1)}
        \geq
        \hat d_{\pi(2)}
        \geq
        \cdots
        \geq
        \hat d_{\pi(n)}\,.
$
Choose
$\hat\tau
\in
[\hat d_{\pi(s+1)},\hat d_{\pi(s)}]$ 
 and define
$
\hat\nu_i
=
(\hat d_i-\hat\tau)_+$\,, and $
\hat\upsilon_i
=
(\hat\tau-\hat d_i)_+$\,, for $
i\in N$.
Then
$
\hat d+\hat\upsilon-\hat\nu-\hat\tau\mathbf e=0$,
with
$\hat\upsilon,\hat\nu\ge0$,
so that
$(\hat\Theta,\hat\upsilon,\hat\nu,\hat\tau)$
is feasible for \ref{DGFactplus}.

We first bound the linear part of the dual objective.
By construction,
\[
\hat\tau s+\hat\nu^\top\mathbf e
=
\textstyle\sum_{\ell=1}^{s}\hat d_{\pi(\ell)}\,.
\]
Moreover,
\[
F\hat\Theta F^\top
=
\alpha(C-\psi I)^{1/2}C^{-1}(C-\psi I)^{1/2}
=
\alpha(I-\psi C^{-1}),
\]
whose eigenvalues are
\[
\textstyle\alpha\left(1-\frac{\psi}{\lambda_\ell(C)}\right),
\qquad
\ell=1,\ldots,n.
\]
Applying to $F\hat\Theta F^\top$ the fact that the sum of the $s$ greatest diagonal entries of a symmetric matrix is at most the sum of its $s$ greatest eigenvalues (a consequence of the Schur--Horn theorem; see \citep*[Thm.~4.3.26]{HJBook}), we have
\begin{equation}\label{est1}
\hat\tau s+\hat\nu^\top\mathbf e
\le
\alpha
\textstyle\sum_{\ell=1}^{s}
\left(
1-\frac{\psi}{\lambda_\ell(C)}
\right)
=
t
-
\textstyle\frac{t\psi}{s}
\sum_{\ell=1}^{s}
\frac1{\lambda_\ell(C)}.
\end{equation}
Next, we evaluate the spectral part.
The $t$ least eigenvalues of
$\hat\Theta=\alpha C^{-1}$ are
$
\textstyle\frac{\alpha}{\lambda_1(C)},
\ldots,
\frac{\alpha}{\lambda_t(C)}.
$
As
$\alpha\le1$
and
$\psi\le\lambda_\ell(C)$, we have 
$
\textstyle\frac{\alpha}{\lambda_\ell(C)}
\le
\frac1\psi$, for $\ell=1,\ldots,t,
$
and therefore
\[
\textstyle h_\psi\!\left(
\frac{\alpha}{\lambda_\ell(C)}
\right)
=
-
\log\!\left(
\frac{\alpha}{\lambda_\ell(C)}
\right)
+
\frac{\alpha\psi}{\lambda_\ell(C)}
-
1.
\]
Hence
\begin{equation}\label{est2}
\textstyle\sum_{\ell=k-t+1}^{k}
h_\psi(\lambda_\ell(\hat\Theta))
=
\textstyle\sum_{\ell=1}^{t}
\log\lambda_\ell(C)
-
t\log\alpha
+
\textstyle\frac{t\psi}{s}
\sum_{\ell=1}^{t}
\frac1{\lambda_\ell(C)}
-
t.
\end{equation}
Combining \eqref{est1} and \eqref{est2},
the objective value of the constructed feasible solution is at most
\[
\begin{aligned}
&
\textstyle\sum_{\ell=1}^{t}\log\lambda_\ell(C)
-
t\log\alpha
-
\frac{t\psi}{s}
\sum_{\ell=t+1}^{s}
\frac1{\lambda_\ell(C)}.
\end{aligned}
\]
Substituting
$\alpha=t/s$
yields
\[
\zspectralthing(C,t)
+
t\log(s/t)
-\textstyle
\frac{t\psi}{s}
\sum_{\ell=t+1}^{s}
\textstyle\frac1{\lambda_\ell(C)}.
\]
As \ref{DGFactplus} is a minimization problem, its optimal value
cannot exceed the objective value of any feasible solution, proving the
desired inequality.
\end{proof}

\begin{remark}
We note that
$
\textstyle \frac{t\psi}{s}\sum_{\ell=t+1}^{s}\frac{1}{\lambda_\ell(C)}
\ge 0.
$
Therefore, Theorem~\ref{thm:dgfactplus_improved_spectral_comparison} immediately
implies
$
\zdgfactplusthing(C,s,t;\psi)
-
\zspectralthing(C,t)
\le
t\log(s/t).
$
Hence, \ref{DGFactplus} satisfies the same worst-case guarantee relative to the
spectral bound as the original DGFact bound established in
\citep*[Thm.~7]{GMESP_Alg}. Moreover, for $\psi>0$ and $s>t$, the correction term is strictly positive, which confirms that   Theorem~\ref{thm:dgfactplus_improved_spectral_comparison} yields a strictly better worst-case comparison for the \ref{DGFactplus} bound than the corresponding comparison known for the DGFact bound.
\end{remark}

\subsection{Sufficient conditions for dominance, and incomparability for general \texorpdfstring{$t$}{t}}

The results of the previous subsections show that the spectral bound and
the \ref{DGFactplus} bound are, in general, incomparable. Nevertheless, it is possible to identify
 conditions based on the multiplicity of the greatest and least eigenvalues of $C$, under which each bound  provably dominates the other. 
 
 In the next theorem we give a sufficient condition under which  the spectral bound dominates the  \ref{DDGFactplus} bound. The theorem is an adaptation of
\citep*[Thm.~33]{ExtendedGMESP} to the augmented formulation. 

\begin{theorem}~\!\!\![Sufficient conditions for dominance of the spectral bound]
\label{thm:ddgfact_plus_spectral_dominance}
Let $C \in \mathbb S_{++}^{n}$\,, $0 < t \leq s < n$, and $0 \leq \psi \leq \lambda_{\min}(C)$. If the multiplicity $\mu_{\max}$ of the greatest eigenvalue of $C$ satisfies
$
\mu_{\max}\ge\frac{nt}{s} 
$, then
$
\zspectralthing(C,t)
\le
\zdgfactplusthing(C,s,t;\psi).
$

\end{theorem}
\begin{proof}
Let
$
B_\psi :=C-\psi I=FF^\top.
$
As
$
\lambda_{\max}(B_\psi)
=
\lambda_{\max}(C)-\psi,
$
the multiplicity of the greatest eigenvalue of  $B_\psi$ is
$\mu_{\max}$\,. Then,
\begin{align}
&\textstyle \tr(B_\psi)\ge \mu_{\max}\lambda_{\max}(B_\psi)
\ge
\frac{nt}{s}
\bigl(\lambda_{\max}(C)-\psi\bigr)\Longrightarrow s\,\tr(B_\psi)+nt\psi
\ge
nt\lambda_{\max}(C)\Longrightarrow\nonumber\\
&\textstyle \frac{
nt\,\lambda_\ell(C)}
{s\,\tr(B_\psi)+nt\psi}
\le1,
\quad
\forall \ell=1,\ldots,t \Longrightarrow \textstyle \sum_{\ell=1}^{t}
\log\!\left(
\frac{
nt\,\lambda_\ell(C)}
{s\,\tr(B_\psi)+nt\psi}
\right)
\le0.\label{eq:ineq}
\end{align}
Moreover, as $\mu_{\max}\ge nt/s\ge t$ and $\lambda_{\max}(B_\psi)\ge0$, we have 
$
\tr(B_\psi)
\ge\mu_{\max}\lambda_{\max}(B_\psi)\ge
t\lambda_{\max}(B_\psi),
$
 which implies 
$
\textstyle\tr\!\left(\frac{s}{n}B_\psi\right)
\ge
t\lambda_{\max}\!\left(\frac{s}{n}B_\psi\right).
$  Then, Lemma~\ref{lem:gamma_trace_upper_bound}  gives
\[
\textstyle\Gamma_t^+\!\left(
\frac{s}{n}B_\psi;\psi
\right)=
t\log\!\left(
\frac{\tr(\frac sn\,B_\psi)}{t}+\psi
\right)
=
t\log\!\left(
\frac{s\,\tr(B_\psi)+nt\psi}{nt}
\right).
\]
Therefore, we have 
\begin{equation}\label{eq:dif}
\textstyle\sum_{\ell=1}^{t}\log\lambda_\ell(C)
-
\Gamma^+_t\!\left(
\frac{s}{n}B_\psi;\psi
\right) = \textstyle \sum_{\ell=1}^{t}
\log\!\left(
\frac{
nt\,\lambda_\ell(C)}
{s\,\tr(B_\psi)+nt\psi}
\right).
\end{equation}
Let
$
\hat x:=\frac{s}{n}\mathbf e.
$
As $\hat x$ is feasible for \ref{DDGFactplus}, we have 
$
\zdgfactplusthing(C,s,t;\psi)
\ge
\Gamma^+_t(M_\psi(\hat x);\psi).
$
Moreover,
$
\textstyle M_\psi(\hat x)
=
F^\top\Diag(\hat x)F
=
\frac{s}{n}F^\top F.
$
As $F^\top F$ and $FF^\top=B_\psi$ have the same nonzero eigenvalues,
$
\textstyle \Gamma^+_t(M_\psi(\hat x);\psi)
=
\Gamma^+_t\!\left(
\frac{s}{n}B_\psi;\psi
\right).
$
Hence
$
\textstyle \zdgfactplusthing(C,s,t;\psi)
\ge
\Gamma^+_t\!\left(
\frac{s}{n}B_\psi;\psi
\right),
$
and 
\[
\textstyle \zspectralthing(C,t)
-
\zdgfactplusthing(C,s,t;\psi)
\le
\sum_{\ell=1}^{t}\log\lambda_\ell(C)
-
\Gamma^+_t\!\left(
\frac{s}{n}B_\psi;\psi
\right)\le 0,
\]
where the last inequality follows from \eqref{eq:ineq} and \eqref{eq:dif}.
\end{proof}

Finally, we give a sufficient condition under which the  \ref{DDGFactplus} bound dominates the spectral bound when $\psi:=\lambda_{\min}(C)$.

\begin{proposition}~\!\!\![Dominance of \ref{DDGFactplus} over the spectral bound under {{a}} sufficient condition]
\label{prop:ddgfactplus_dominates_spectral}
Let $C\in\mathbb{S}_{++}^n$, $1\le t\le s<n$, and set
$
\psi:=\lambda_{\min}(C).
$
If the multiplicity $\mu_{\min}$ of the least eigenvalue of $C$ satisfies
$
\mu_{\min}\ge n-t,
$
then
$
\zddgfactplusthing(C,s,t;\lambda_{\min}(C))
\le
\zspectralthing(C,t).
$
\end{proposition}

\begin{proof}
Let $\hat x$ be any feasible solution of \ref{DDGFactplus}. As
$
\mu_{\min}\ge n-t
$, by
Theorem~\ref{thm:gamma_t_collapse_multiplicity}, we have
$
\Gamma_t^+(M_\psi(\hat x);\psi)
=
\textstyle \sum_{\ell=1}^t
\log\!\left(\lambda_\ell(F^\top\Diag(\hat x)F)+\psi\right).
$
As $0\preceq\Diag(\hat x)\preceq I$,
$
F^\top\Diag(\hat x)F\preceq F^\top F,
$
and hence, 
$
\lambda_\ell(F^\top\Diag(\hat x)F)
\le
\lambda_\ell(F^\top F)$, for $
\ell=1,\ldots,t.
$
Moreover, $F^\top F$ and $FF^\top=C-\psi I$ have the same positive eigenvalues. Because
$
\mu_{\min}\ge n-t,
$
we have
$
\rank(C-\psi I)\le t.
$
Hence,
$
\lambda_\ell(F^\top F)+\psi
=
\lambda_\ell(C)$, for $
\ell=1,\ldots,t$. 
Therefore,
$
\Gamma_t^+(M_\psi(\hat x);\psi)
\le
\textstyle\sum_{\ell=1}^t
\log\!\bigl(\lambda_\ell(C)\bigr)
=
\zspectralthing(C,t).
$
As the inequality holds for every feasible $\hat x$, taking the maximum over all feasible solutions of \ref{DDGFactplus} gives the inequality claimed.
\end{proof}

The examples of Subsection~\ref{subsec:exact_comparison_t1} concern the case
$t=1$. Theorem~\ref{thm:ddgfact_plus_spectral_dominance} and
Proposition~\ref{prop:ddgfactplus_dominates_spectral} allow us to exhibit
instances on which each bound is strictly stronger than the other. We 
give two families of such instances, for every $t\ge1$, in Example~\ref{ex:incomparable_general_t}, showing that the two bounds are incomparable for every value of $t$.

{
\begin{example}
\label{ex:incomparable_general_t}
Let $t\ge1$, and take $\psi:=\lambda_{\min}(C)$ in each of the following two families of instances. 

($a$) Let $n:=2t$, let $s$ satisfy $t\le s<n$, and let $C:=I_n+\sum_{j=1}^t u_ju_j^\top$, where $u_j\in\mathbb{R}^n$ is the characteristic vector of $\{2j-1,2j\}$, for $j=1,\ldots,t$. The eigenvalues of $C$ are $3$ and $1$, each with multiplicity $t$, so that $\psi=1$ and $\mu_{\min}=t=n-t$, and Proposition~\ref{prop:ddgfactplus_dominates_spectral} applies. We may take $F:=[u_1,\ldots,u_t]\in\mathbb{R}^{n\times t}$, as $k=t=\max\{t,\rank(C-\psi I)\}$, and then $M_\psi(x)=\Diag(p)$, where $p_j:=x_{2j-1}+x_{2j}$\,, for $j=1,\ldots,t$. By Theorem~\ref{thm:gamma_t_collapse_multiplicity}, $\Gamma_t^+(M_\psi(x);\psi)=f_t(M_\psi(x);\psi)=\sum_{j=1}^t\log(p_j+1)$ for every feasible $x$. The feasible region of \ref{DDGFactplus} maps onto $\{p\in[0,2]^t:\sum_{j=1}^tp_j=s\}$, and, by the concavity of the logarithm, $\sum_{j=1}^t\log(p_j+1)\le t\log\big(\tfrac st+1\big)$ on this set, with equality at $p_j=\tfrac st\in[0,2]$, for $j=1,\ldots,t$. Hence $\zddgfactplusthing(C,s,t;\psi)=t\log\big(\tfrac st+1\big)<t\log3=\zspectralthing(C,t)$, as $s<2t$.

($b$) Let $n:=t+2$, $s:=t+1$, and $C:=\alpha I_n-\mathbf{e}\mathbf{e}^\top$, with $\alpha>n$. The eigenvalues of $C$ are $\alpha$, with multiplicity $n-1$, and $\alpha-n>0$, so that $\psi=\alpha-n$ and $C-\psi I=nI_n-\mathbf{e}\mathbf{e}^\top$. The point $\hat x:=\frac sn\mathbf{e}$ is feasible for \ref{DDGFactplus}, and, with $F:=(C-\psi I)^{1/2}$, we have $M_\psi(\hat x)=\frac sn(nI_n-\mathbf{e}\mathbf{e}^\top)$, whose eigenvalues are $s$, with multiplicity $n-1$, and $0$. Thus $\tr(M_\psi(\hat x))=s(n-1)\ge ts=t\,\lambda_{\max}(M_\psi(\hat x))$, as $n-1=t+1$, and Lemma~\ref{lem:gamma_trace_upper_bound} gives $\Gamma_t^+(M_\psi(\hat x);\psi)=t\log\big(\tfrac{s(n-1)}{t}+\alpha-n\big)$. As $s(n-1)=(t+1)^2>t(t+2)=tn$, we obtain
\[
\zddgfactplusthing(C,s,t;\psi)\ge t\log\Big(\frac{s(n-1)}{t}+\alpha-n\Big)>t\log\alpha=\zspectralthing(C,t).
\]
\hfill $\clubsuit$
\end{example}
}

\section{Computational results} \label{sec:computational_results}

In this section, we report on two experiments, with the main goal of evaluating the relative performance of \ref{DDGFactplus} and DDGFact. For our first experiment, at the root node, we compare the \ref{DDGFactplus} bound with the generalized factorization bound DDGFact of \citep*{SEA_proceedings,GMESP_Alg}.
To put this experiment into a broader context, we compare also with the scaled glinx bound of \citep*{ExtendedGMESP}, and the spectral bound of \citep*{LeeLind2020}, in terms of root gaps and runtimes. 
In our second experiment, to see further how \ref{DDGFactplus} can improve on DDGFact, we embed both DDGFact and \ref{DDGFactplus} in a B\&B algorithm for \ref{GMESP} and compare their performance. For each of them, we also evaluate two enhancements: the trace upper bound of Lemma~\ref{lem:gamma_trace_upper_bound}, which is used to prune nodes, and the closed-form certificate of Theorem~\ref{thm:ddgfact_plus_closed_form_iota_zero_certificate}, which avoids solving the relaxation at the nodes where it applies. Throughout this section, we refer to the augmented generalized factorization bound as the \ref{DDGFactplus} bound, because \ref{DDGFactplus} is the problem that we solve with KNITRO; by strong duality, it gives the same bound as \ref{DGFactplus}.

In all B\&B experiments with DDGFact and \ref{DDGFactplus}, we fix variables, a well-established technique for accelerating B\&B algorithms that has proved particularly effective for \ref{MESP} and \ref{GMESP}; see, e.g., \citep*{AFLW_IPCO,AFLW_Using,GMESP_Alg}. At each node, from the solution $\hat x$ of DDGFact$^+(l,c)$ returned by the solver, we construct the feasible solution $(\hat\Theta,\hat\upsilon,\hat\nu,\hat\tau)$ of DGFact$^+(l,c)$ given by Proposition~\ref{prop:dgfactplus_associated_solution_fixed}, with objective value $\hat\zeta$, and we fix $x_i$ at $0$ (resp., at $1$) whenever $\hat\zeta-{\rm LB}<\hat\upsilon_i$ (resp., $\hat\zeta-{\rm LB}<\hat\nu_i$), where ${\rm LB}$ is the value of the best feasible solution of \ref{GMESP} known. The validity of this rule follows from Theorem~\ref{thm:variable_fixing_dgfact_plus} in the Appendix, whose proof applies, with the obvious changes, to the subproblem at each node. For DDGFact, we apply the same rule with $\psi:=0$.

\subsection{Programming environment}\label{subsec:prog_env}

We used the general-purpose nonlinear-programming solver KNITRO 14.0.0 (see \citep*{KNITRO}), via the Julia wrapper KNITRO.jl v0.14.8, to solve \ref{DDGFactplus}, and the conic solver Hypatia 0.9.2 (see \citep*{coey2022solving}) to solve the semidefinite programs arising in the benchmark bounds. For KNITRO, we used the settings FEASTOL $=10^{-6}$ (feasibility tolerance) and OPTTOL $=10^{-6}$ (absolute optimality tolerance). For Hypatia we used the default parameters. To reduce the runtime of our experiments, we used 16 threads using the linear algebra backend MKL v0.9.1.

We used the programming language Julia v1.12.6 (with VSCodium v1.126.04524), and we ran the experiments on ``zebratoo'', a 32-core machine (running Windows Server 2022 Standard) with two Intel Xeon Gold 6444Y processors running at 3.60GHz, with 16 cores each, and 128 GB of memory.

\subsection{\texorpdfstring{{Root-gap comparison on benchmark instances}}{Root-gap comparison on benchmark instances}}\label{subsec:root-gap_comparison}

For our instances, we used two benchmark positive-definite covariance matrices $C$, of orders $n=63$ and $n=90$, originally obtained from J.\ Zidek (University of British Columbia). The first arises from the redesign of an environmental monitoring network, and the second is based on temperature data from monitoring stations in the Pacific Northwest of the United States; see \citep*{Guttorp-Le-Sampson-Zidek1993, HLW}. These matrices have been used extensively to test and develop algorithms for \ref{MESP}; see \citep*{KLQ, LeeConstrained, AFLW_Using, LeeWilliamsILP, AnstreicherLee_Masked, BurerLee, Anstreicher_BQP_entropy, Kurt_linx, Mixing, FactPaper}.

For each matrix, we considered $\kappa\in\{2,4,6,8\}$, $s\in\{10,20,\ldots,60\}$ for $n=63$ and $s\in\{10,20,\ldots,80\}$ for $n=90$, and $t:=s-\kappa$; each instance thus corresponds to a tuple $(n,s,t)$. In applications to principal component analysis (PCA), large values of $\kappa$ are the most relevant; however, for instances with large $\kappa$, the scaled glinx bound of \citep*{ExtendedGMESP} gives small gaps, but it requires solving a semidefinite program, which limits its use within B\&B for larger instances; the \ref{DDGFactplus} bound is a contribution towards cheaper bounds in this regime. As observed in \citep*{SEA_proceedings,GMESP_Alg}, larger values of $\kappa$ correspond to harder instances, so we expect the gaps of all methods to grow with $\kappa$. For each instance and each method, we computed the root gap, that is, the difference between the upper bound produced by the method at the root node and a lower bound given by the objective value of a feasible solution of \ref{GMESP}. The lower bounds were obtained with the local-search procedures of \citep*{GMESP_Alg}, to which we added the vector $\bar x$ of Theorem~\ref{thm:ddgfact_plus_closed_form_iota_zero_certificate} as an additional initial solution.

Our benchmark methods are the generalized factorization bound of \citep*{SEA_proceedings,GMESP_Alg}, the scaled glinx bound of \citep*{ExtendedGMESP}, which we refer to simply as the glinx bound in the remainder of this section, and the spectral bound of \citep*{LeeLind2020}. The glinx bound is computed by solving a semidefinite program, whereas both factorization bounds are convex problems in a vector variable, so we expect a trade-off between the quality of the bound and the time needed to compute it. Regarding the generalized factorization bound, Theorem~\ref{thm:ddgfact_plus_monotonicity} already guarantees that the augmented bound is at least as strong for every $\psi>0$; the purpose of this comparison is to measure how large the gap reduction is in practice. As the bound is nonincreasing in $\psi\in[0,\lambda_{\min}(C)]$, we use $\psi:=\lambda_{\min}(C)$ in all experiments with \ref{DDGFactplus}. 
 The size of the resulting improvement cannot be inferred from the value of $\lambda_{\min}(C)$ alone. Indeed, replacing $C$ and $\psi$ by $\alpha C$ and $\alpha\psi$, for $\alpha>0$, adds $t\log\alpha$ to both factorization bounds and to the optimal value of \ref{GMESP}, so the gaps do not change. Moreover, the proof of Theorem~\ref{thm:ddgfact_plus_strict_monotonicity}, with $\psi_1:=0$ and $\psi_2:=\lambda_{\min}(C)$, shows that, when $\hat\iota\ge1$, the decrease of the bound is at least $\lambda_{\min}(C)\,(g_{\hat\iota+1}-g_{\hat\iota})\big(\hat\iota-\sum_{i=1}^{\hat\iota}x^*_{\pi(i)}\big)$, where the components of $g$ are reciprocals of components of $\beta^*=\lambda(M_0(x^*))$. This lower bound depends on $\lambda_{\min}(C)$ relative to the eigenvalues of $M_0(x^*)$, and also on the optimal solution $x^*$ of \ref{DDGFactplus}, which differs from instance to instance.
For our instances with $n=63$ and $n=90$, we have $\lambda_{\min}(C)\approx0.0317$ and $\lambda_{\min}(C)\approx19.0$, respectively.

\begin{figure}[htbp]
    \centering
    \begin{subfigure}[t]{0.49\textwidth}
        \centering
        \includegraphics[width=\textwidth]{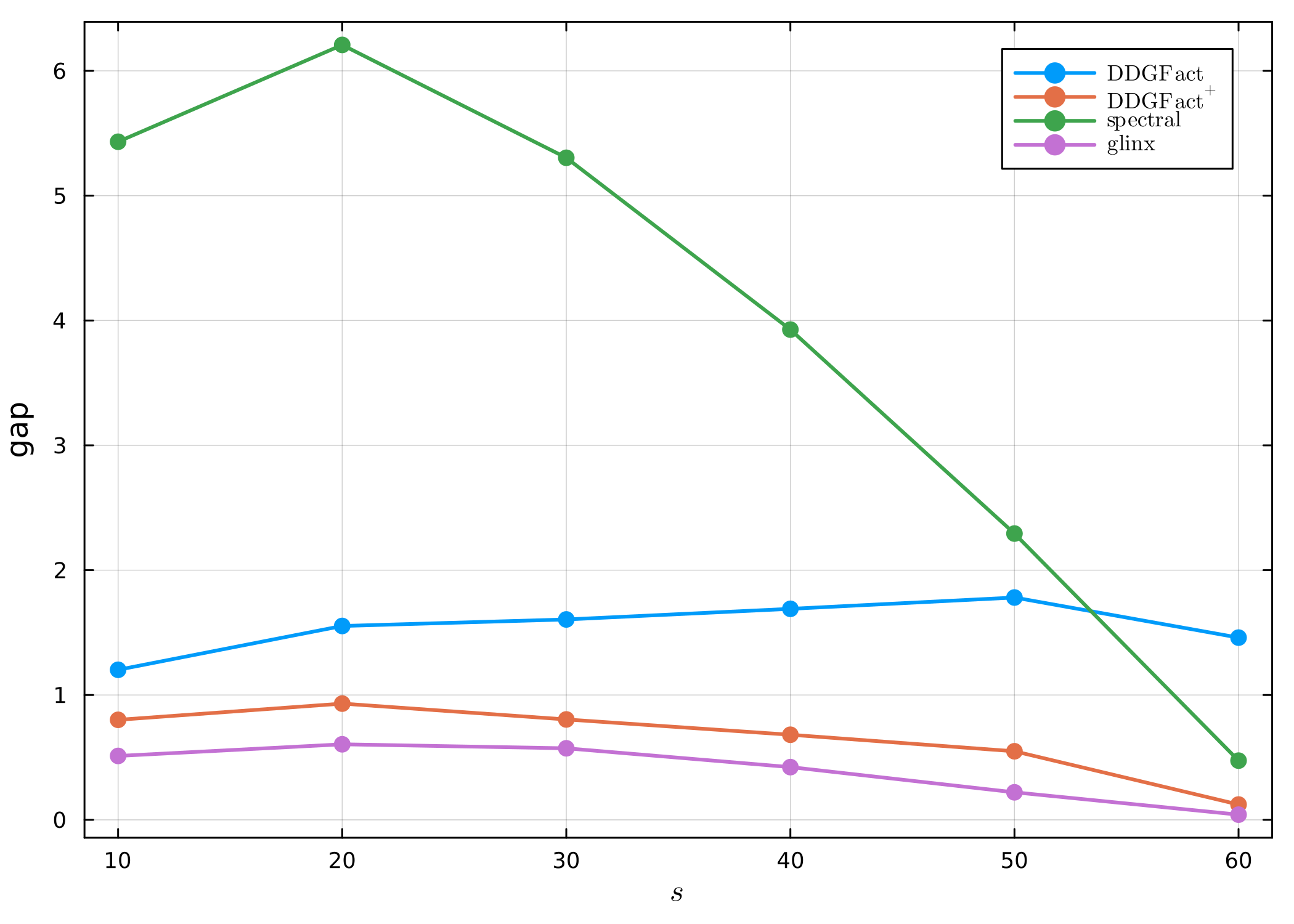}
        \caption{$\kappa = 2$.}
        \label{fig:plot_n63_kappa2}
    \end{subfigure}
    \hfill
    \begin{subfigure}[t]{0.49\textwidth}
        \centering
        \includegraphics[width=\textwidth]{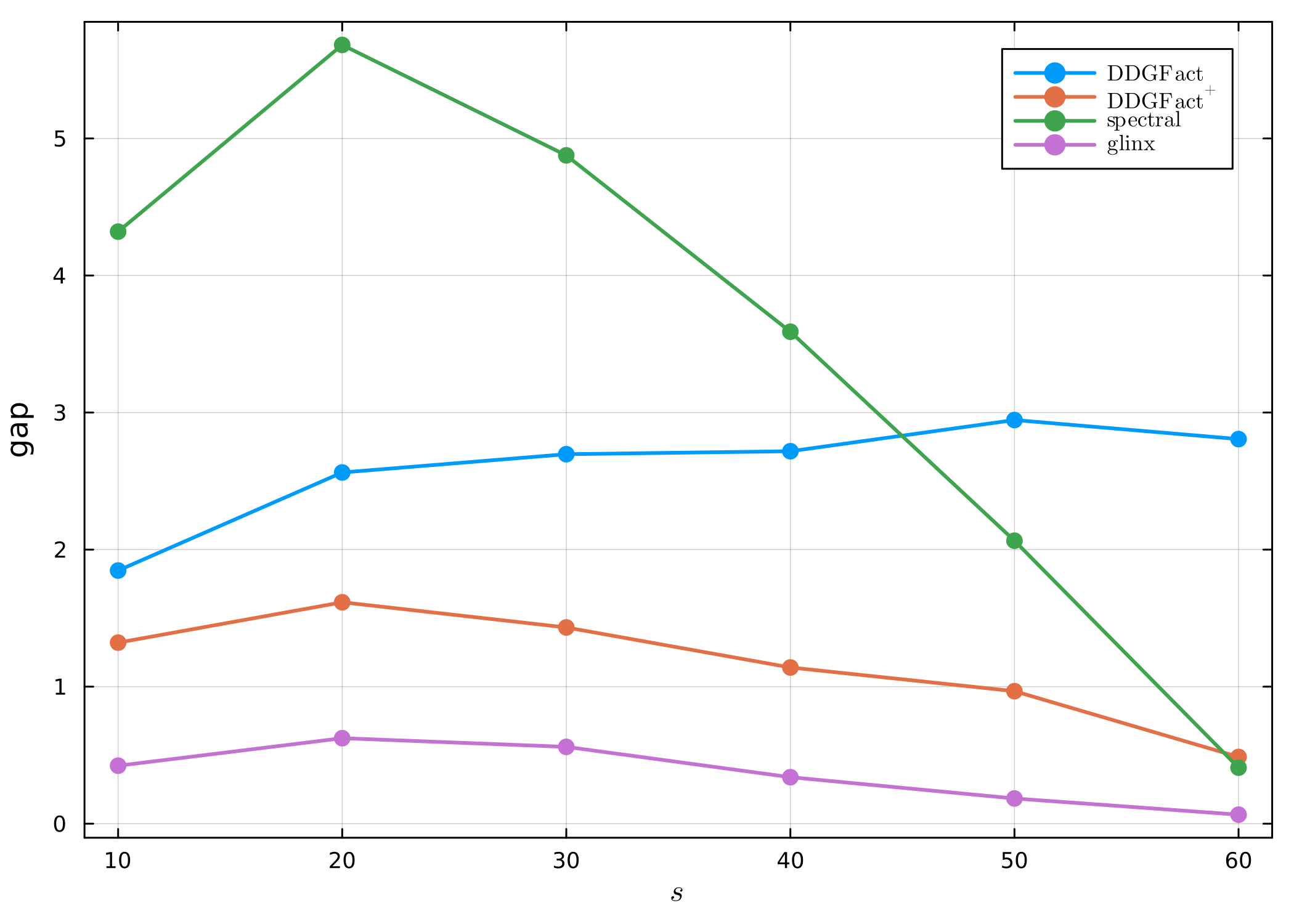}
        \caption{$\kappa = 4$.}
        \label{fig:plot_n63_kappa4}
    \end{subfigure}

    \medskip

    \begin{subfigure}[t]{0.49\textwidth}
        \centering
        \includegraphics[width=\textwidth]{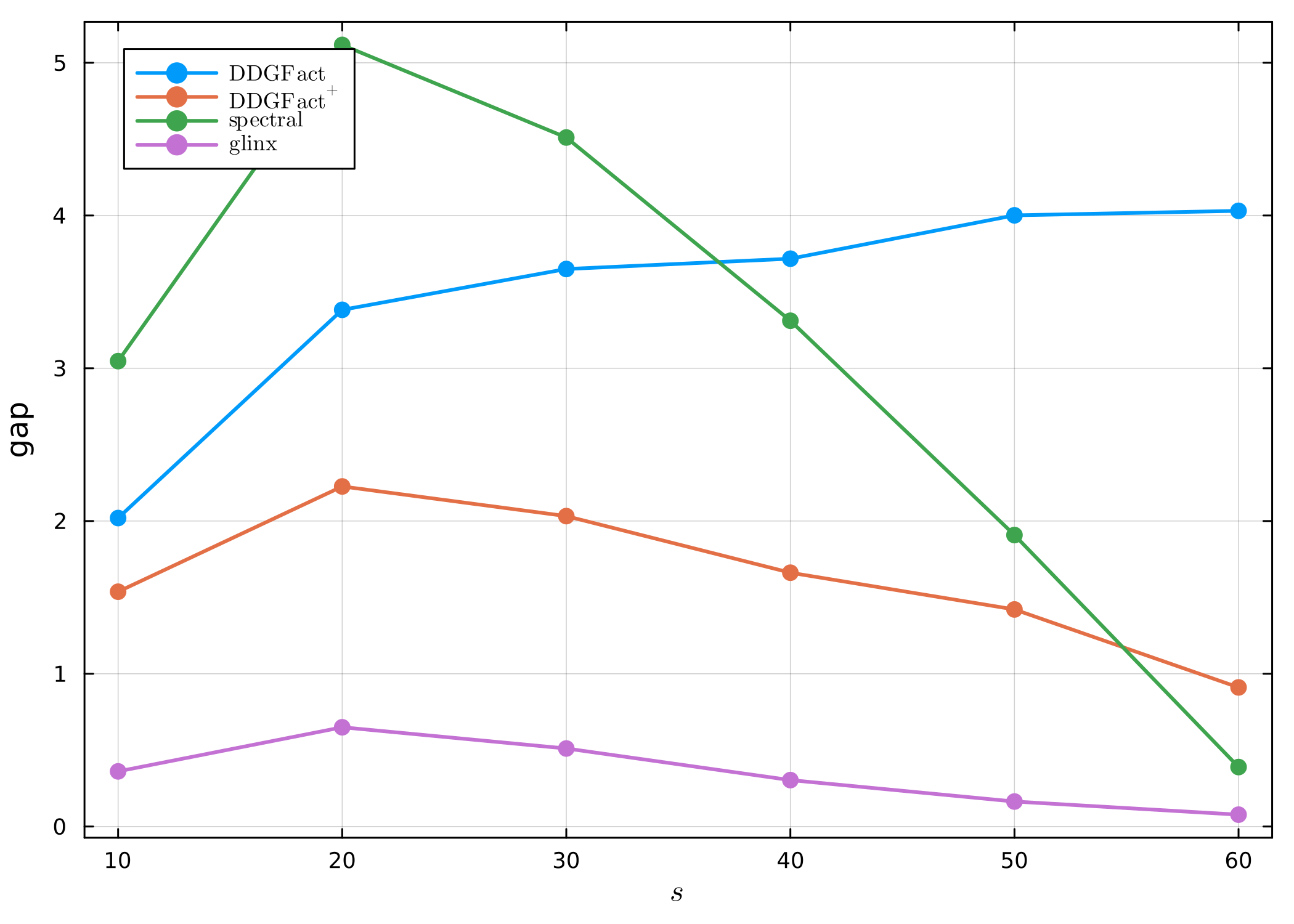}
        \caption{$\kappa = 6$.}
        \label{fig:plot_n63_kappa6}
    \end{subfigure}
    \hfill
    \begin{subfigure}[t]{0.49\textwidth}
        \centering
        \includegraphics[width=\textwidth]{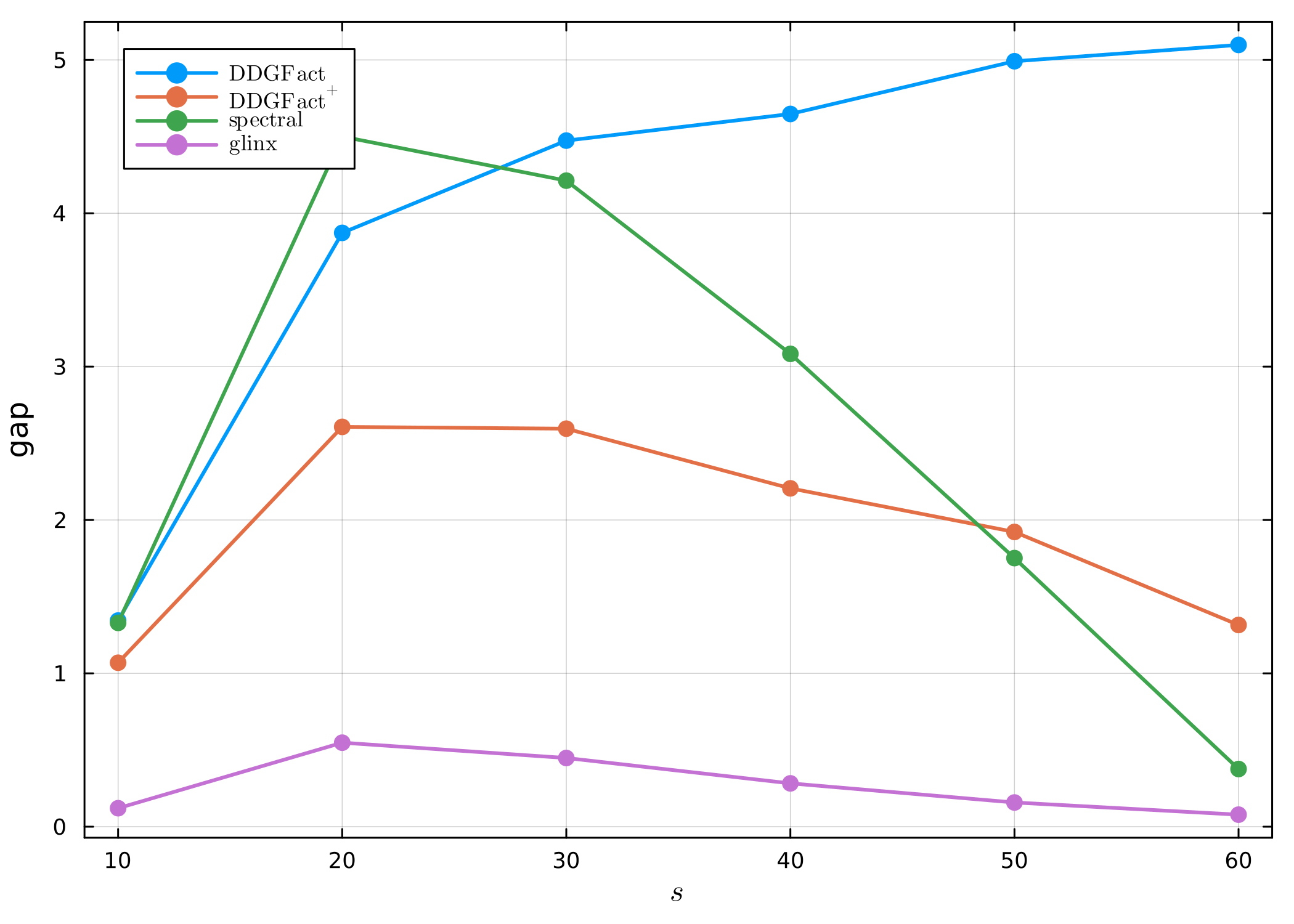}
        \caption{$\kappa = 8$.}
        \label{fig:plot_n63_kappa8}
    \end{subfigure}

    \caption{Root gaps for $n = 63$ and $\kappa \in \{2, 4, 6, 8\}$.}
    \label{fig:plots_n63_kappa2_4_6_8}
\end{figure}

\begin{figure}[htbp]
    \centering
    \begin{subfigure}[t]{0.49\textwidth}
        \centering
        \includegraphics[width=\textwidth]{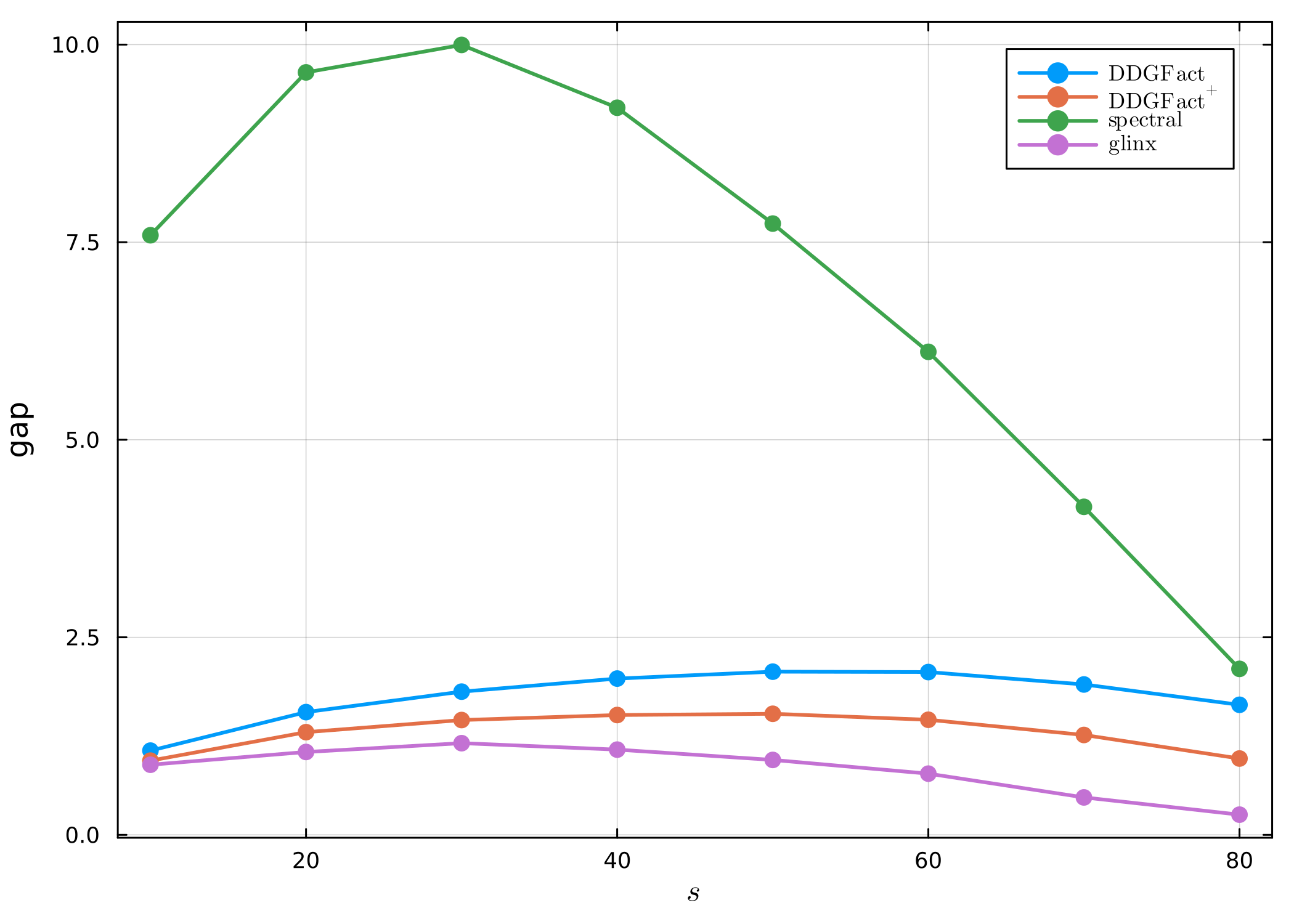}
        \caption{$\kappa = 2$.}
        \label{fig:plot_n90_kappa2}
    \end{subfigure}
    \hfill
    \begin{subfigure}[t]{0.49\textwidth}
        \centering
        \includegraphics[width=\textwidth]{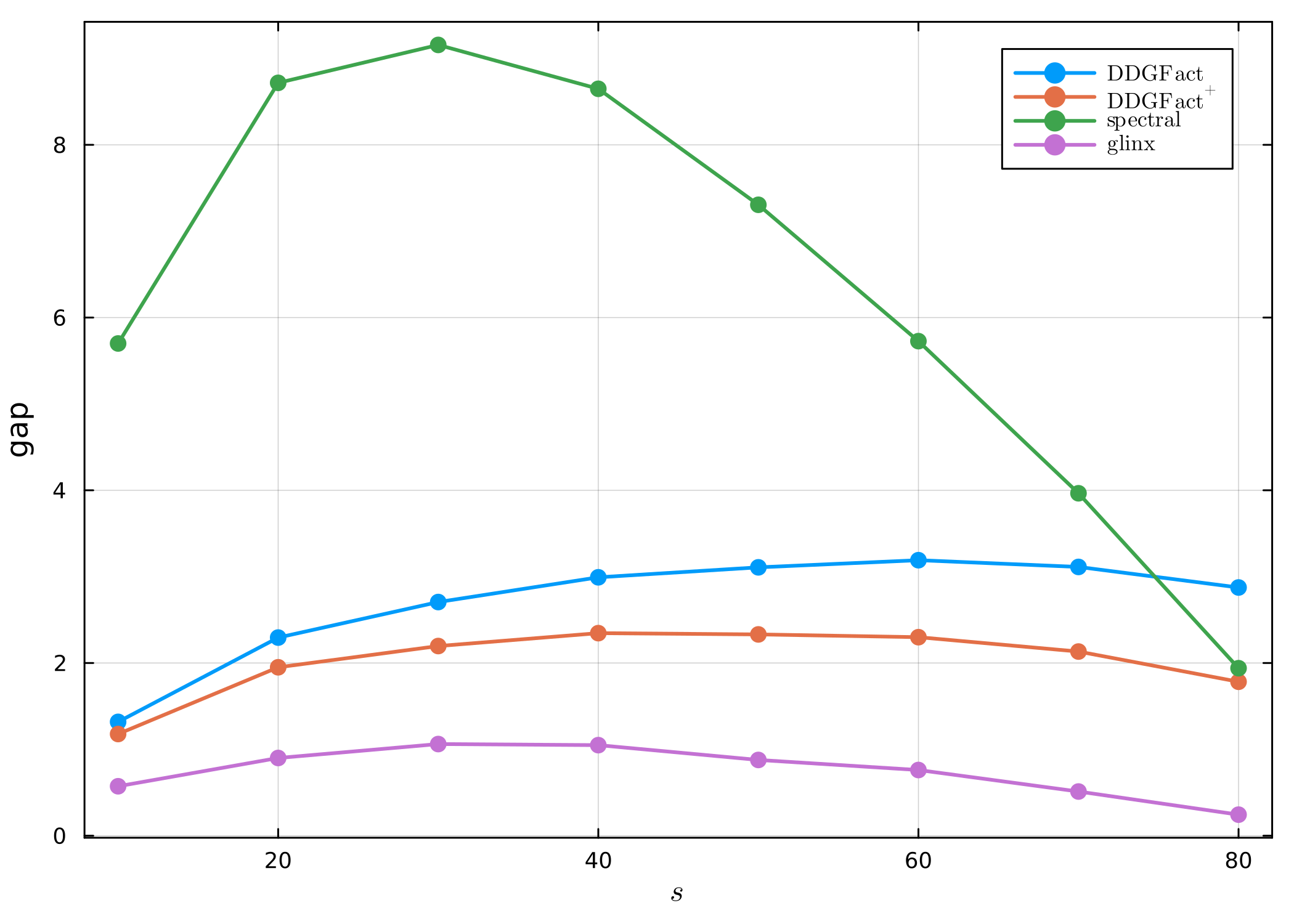}
        \caption{$\kappa = 4$.}
        \label{fig:plot_n90_kappa4}
    \end{subfigure}

    \medskip

    \begin{subfigure}[t]{0.49\textwidth}
        \centering
        \includegraphics[width=\textwidth]{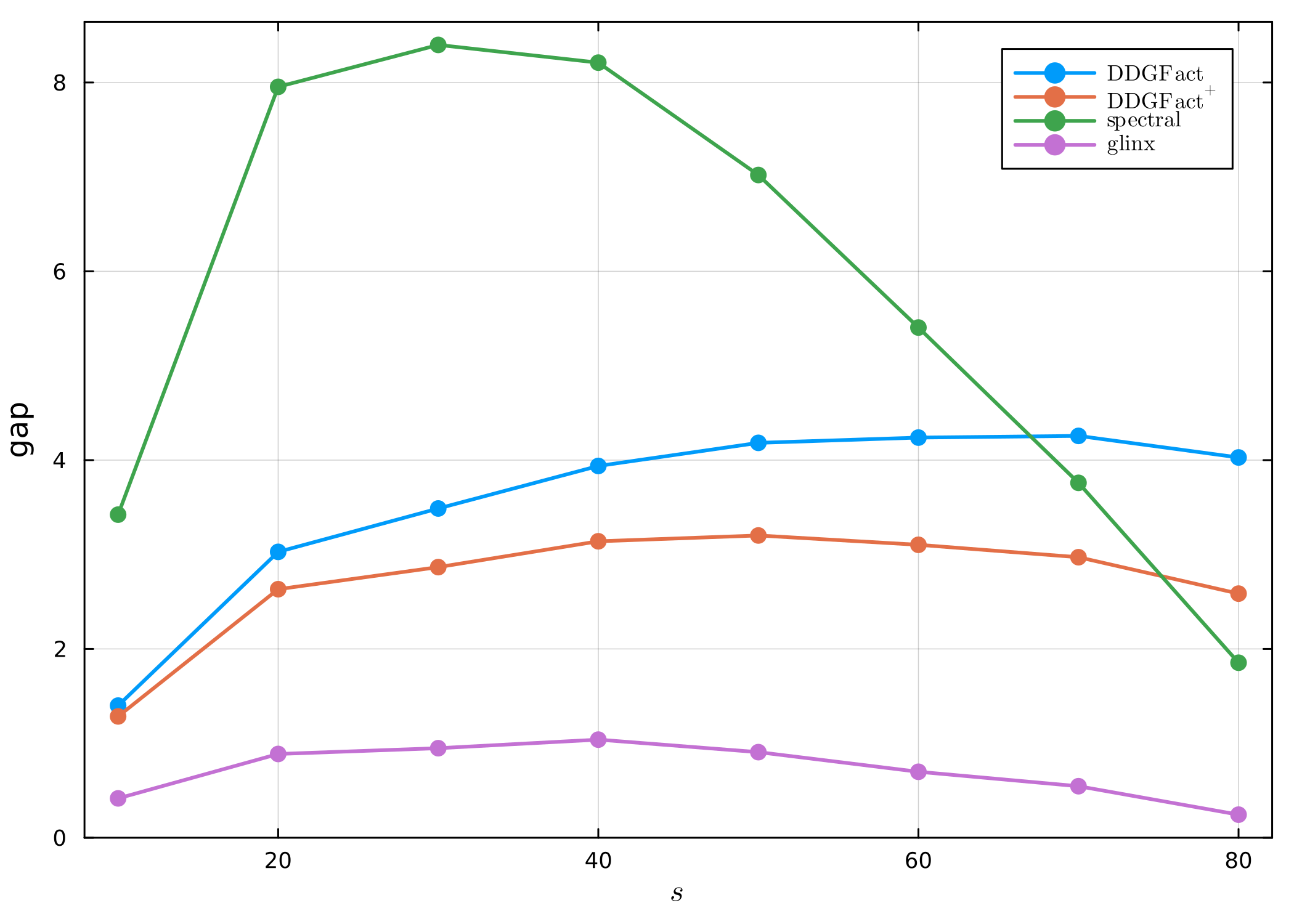}
        \caption{$\kappa = 6$.}
        \label{fig:plot_n90_kappa6}
    \end{subfigure}
    \hfill
    \begin{subfigure}[t]{0.49\textwidth}
        \centering
        \includegraphics[width=\textwidth]{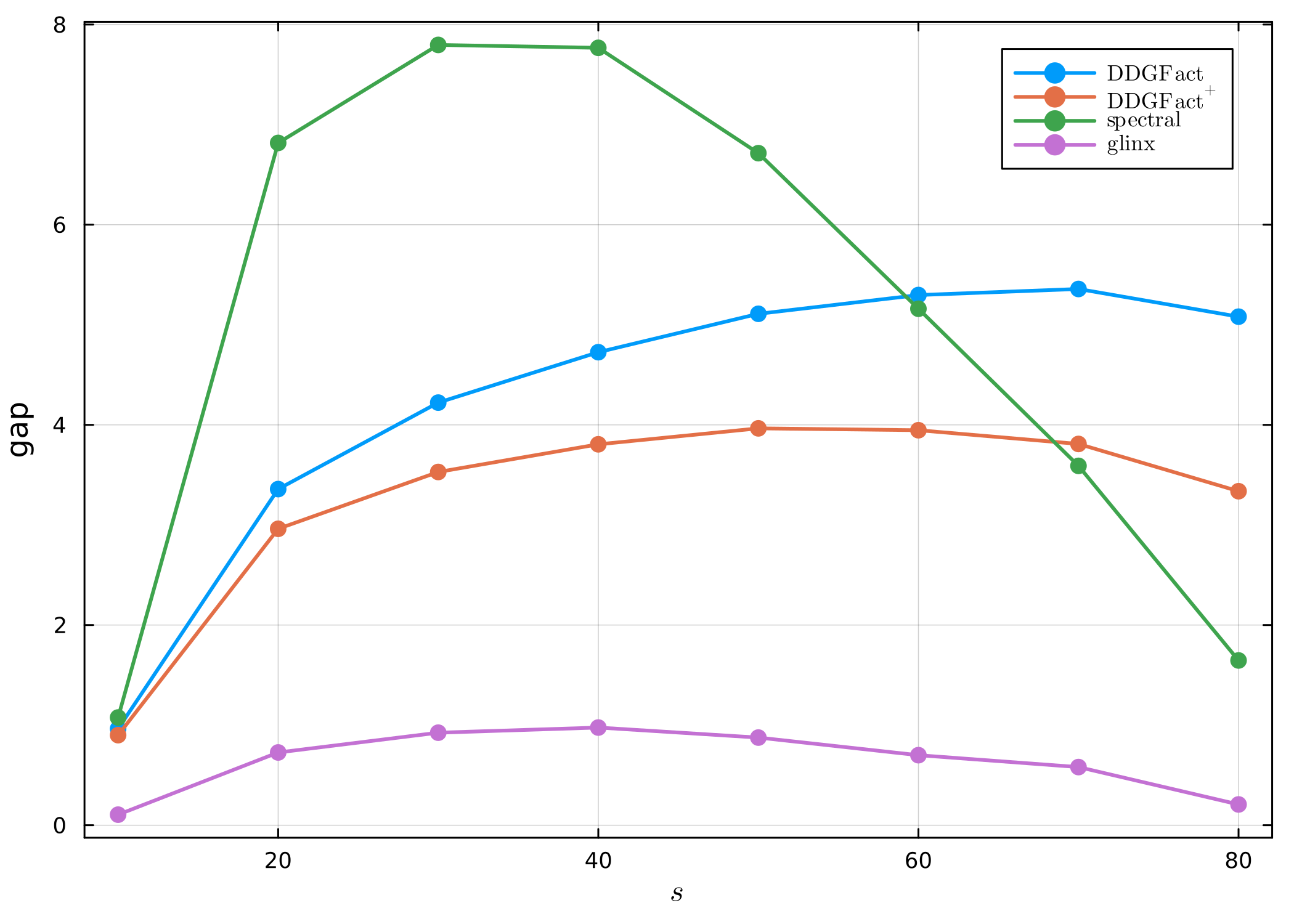}
        \caption{$\kappa = 8$.}
        \label{fig:plot_n90_kappa8}
    \end{subfigure}

    \caption{Root gaps for $n = 90$ and $\kappa \in \{2, 4, 6, 8\}$.}
    \label{fig:plots_n90_kappa2_4_6_8}
\end{figure}

In Figures~\ref{fig:plots_n63_kappa2_4_6_8} and~\ref{fig:plots_n90_kappa2_4_6_8}, the blue and orange curves correspond to the \ref{DDGFactplus} bound with $\psi=0$, which is the DDGFact bound of \citep*{SEA_proceedings,GMESP_Alg}, and with $\psi=\lambda_{\min}(C)$, respectively; the green curve corresponds to the spectral bound, and the violet curve to the glinx bound.

Figure~\ref{fig:plots_n63_kappa2_4_6_8} shows the root gaps for $n=63$. The glinx bound gives the smallest gaps on all instances. On most instances, the \ref{DDGFactplus} bound with $\psi=\lambda_{\min}(C)$ comes second; the exceptions are the instances with $s=60$ and $\kappa\in\{4,6,8\}$, and with $s=50$ and $\kappa=8$, for which the spectral bound gives a smaller gap. The augmentation substantially reduces the gaps of the DDGFact bound, especially for large $s$: for $s=60$ and $\kappa=8$, the gap decreases by about 4, which can considerably reduce the number of nodes explored in B\&B. The spectral bound is not competitive for small and moderate values of $s$, but its gap decreases sharply as $s$ grows. We also note that, for every $\kappa$, the gap of DDGFact tends to increase with $s$, whereas the gap of \ref{DDGFactplus} decreases for $s\ge 20$; hence, the benefit of the augmentation grows with $s$.

Figure~\ref{fig:plots_n90_kappa2_4_6_8} shows that the results for $n=90$ are similar: the glinx bound again gives the smallest gaps, and the spectral bound is stronger than \ref{DDGFactplus} only for large $s$, namely for $s=80$ with $\kappa\in\{6,8\}$ and for $s=70$ with $\kappa=8$. The gap reduction obtained by the augmentation is, however, smaller than for $n=63$: for $\kappa=8$, it is at most about $1.7$ (for $s=80$), against about $4$ for $n=63$. Although $\lambda_{\min}(C)$ is much larger for $n=90$ than for $n=63$, this does not contradict the smaller improvement, as the improvement depends on $\lambda_{\min}(C)$ relative to the eigenvalues of the matrices $M_0(x^*)$ and on the solution $x^*$ of the relaxation. For both values of $n$, the gap of \ref{DDGFactplus} grows considerably with $\kappa$, whereas the gap of the glinx bound remains below about $1.2$ on all instances.

Table~\ref{tab:runtime_stats_n63_n90} reports statistics of the runtimes, in seconds, of the DDGFact bound, the \ref{DDGFactplus} bound, and the glinx bound, computed over all the instances described above for each value of $n$; the spectral bound is omitted, as it only requires an eigenvalue computation. For both values of $n$, the mean runtimes of the DDGFact and \ref{DDGFactplus} bounds are below $0.1$ seconds, whereas those of the glinx bound are about $11.3$ seconds for $n=63$ and $65.9$ seconds for $n=90$; hence, both factorization bounds are two to three orders of magnitude faster on average, and \ref{DDGFactplus} is, on average, $1.5$ to $1.7$ times faster than DDGFact. The results illustrate the trade-off between the quality of bounds and their computational cost. Moreover, as the computation of the glinx bound requires solving a semidefinite program, its runtime grows quickly with $n$. This is relevant in B\&B, where a faster bound allows more nodes to be explored within a time limit (see Subsection~\ref{subsec:bnb_cf_trace}). The standard deviations are small relative to the means for the \ref{DDGFactplus} bound and the glinx bound. For the DDGFact bound, the standard deviation is also small for $n=63$, but for $n=90$ it exceeds the mean, indicating a larger variability of its runtimes across these instances.

\begin{table}[htbp]
\centering
\caption{Runtime statistics (in seconds) for $n = 63$ and $n = 90$.}
\label{tab:runtime_stats_n63_n90}
\begin{tabular}{lrrrrrr}
\toprule
 & \multicolumn{3}{c}{$n = 63$} & \multicolumn{3}{c}{$n = 90$} \\
\cmidrule(lr){2-4} \cmidrule(lr){5-7}
 & DDGFact & DDGFact$^+$ & glinx & DDGFact & DDGFact$^+$ & glinx \\
\midrule
Mean      & 0.0432 & 0.0262 & 11.3385 & 0.0821 & 0.0542 & 65.9029 \\
Std. Dev. & 0.0090 & 0.0083 & 0.7132  & 0.0948 & 0.0148 & 3.4855 \\
\bottomrule
\end{tabular}
\end{table}

\subsection{B\&B with the closed-form certificate and the trace upper bound}\label{subsec:bnb_cf_trace}

The computational cost of evaluating upper bounds is an important factor in the efficiency of B\&B algorithms for \ref{GMESP}. In particular, solving a continuous relaxation at every node can be substantially more expensive than evaluating an upper bound in closed form. Lemma \ref{lem:gamma_trace_upper_bound} provides a computationally inexpensive upper bound, which we refer to as the trace upper bound, given by
\begin{equation}
    \ztracething(C, s, t;\psi)
    := t\log\left(
        \tfrac{\tr((C-\psi I)[\bar{S},\bar{S}])}{t}+\psi
    \right),
    \label{trace_upper_bound}
    \tag{Trace$^+$}
\end{equation}
where $\bar{S}$ is a subset of cardinality $s$ corresponding to the $s$ greatest diagonal entries of $C$. When some variables are fixed at one, the trace upper bound is computed by including their corresponding indices and selecting the remaining indices associated with the greatest diagonal entries of $C$ among the free variables.

By Lemma~\ref{lem:gamma_trace_upper_bound}, the trace upper bound is never stronger than the \ref{DDGFactplus} bound, but Theorem~\ref{thm:ddgfact_plus_closed_form_iota_zero_certificate} and Corollary~\ref{cor:ddgfact_plus_closed_form_iota_zero_node} give a sufficient condition under which the two bounds coincide, at the root and at the B\&B subproblems, respectively; in this case, the relaxation does not need to be solved. Note that, as $(C-\psi I)[\bar S,\bar S]=C[\bar S,\bar S]-\psi I_s$\,, the condition of Theorem~\ref{thm:ddgfact_plus_closed_form_iota_zero_certificate} is equivalent to
\begin{equation}\label{eq:cf_condition_psi}
\tr(C[\bar S,\bar S])-t\,\lambda_{\max}(C[\bar S,\bar S])\ \ge\ (s-t)\,\psi,
\end{equation}
and an analogous equivalence holds for Corollary~\ref{cor:ddgfact_plus_closed_form_iota_zero_node}. Therefore, when $t<s$, the condition becomes more restrictive as $\psi$ increases; moreover, as its left-hand side decreases with $t$, it is more likely to hold for small values of $t$.

We compare the following strategies. For DDGFact and \ref{DDGFactplus}, we either solve the relaxation at every node, or use a hybrid strategy, denoted by $\mathrm{DDGFact}\text{(CF)}$ and $\mathrm{DDGFact}^{+}\text{(CF)}$, respectively, which first tries to prune the node with the trace upper bound, then checks the sufficient condition, and solves the relaxation only if the condition does not hold. As the trace upper bound is never smaller than the bound given by the relaxation, every node pruned by it would also be pruned by the relaxation; hence, the hybrid strategy prunes exactly the nodes that the relaxation would prune and, at the other nodes, computes the same upper bound (up to the solver tolerance), so that the two strategies differ essentially in the time spent per node. We do not use the glinx bound in these experiments: although it gives the smallest root gaps (see Subsection~\ref{subsec:root-gap_comparison}), it is two to three orders of magnitude more expensive to compute (see Table~\ref{tab:runtime_stats_n63_n90}), which makes its use at every node not competitive.

Regarding the branching rule, fractional branching was used in \citep*{GMESP_Alg}. However, at the nodes where the closed form applies, $\mathrm{DDGFact}\text{(CF)}$ and $\mathrm{DDGFact}^{+}\text{(CF)}$ yield a binary solution, namely the characteristic vector of $\bar{S}$ (of $J\cup\bar S$, at a subproblem). We therefore extend the fractional branching rule as follows: when the solution is binary, we branch on the free variable with the largest diagonal entry of $C$.

All B\&B experiments use the matrix with $n=63$ and a time limit of 2 hours. We consider two sets of instances. In the first one, $t = 2$ and $s \in \{5, 10, 15, 20, 25\}$; by \eqref{eq:cf_condition_psi}, small values of $t$ favor the closed-form certificate. In the second one, $\kappa = 5$, $t := s-\kappa$, and $s \in \{10, 20, 30, 40, 50\}$; as discussed in Subsection~\ref{subsec:root-gap_comparison}, instances with larger $\kappa$ are harder and are the most relevant for PCA. In the tables, ``Root Gap'' and ``Gap'' are the gaps at the root node and at the end of the run, ``Nodes'' is the number of nodes explored, ``Fixed'' is the number of variables fixed, ``Time'' is the runtime in seconds, ``CF'' is the number of nodes at which the relaxation was solved in closed form, and ``Prunes'' is the number of nodes pruned by the trace upper bound.

Table~\ref{tab:bnb_comparison_n63_t2} reports the results for $t=2$. Only the instance with $s=5$ was solved within the time limit, by every strategy except DDGFact. For this instance, the augmentation is decisive: DDGFact stops with a gap of $0.199$, whereas \ref{DDGFactplus} closes the gap in $2698.4$ seconds, exploring less than half of the nodes explored by DDGFact. The hybrid strategy accelerates both bounds: $\mathrm{DDGFact}\text{(CF)}$ solves the instance in $382.5$ seconds, and $\mathrm{DDGFact}^{+}\text{(CF)}$ is $4.3$ times faster than \ref{DDGFactplus}. Note that $\mathrm{DDGFact}\text{(CF)}$ is faster than $\mathrm{DDGFact}^{+}\text{(CF)}$, despite its larger root gap. This is explained by \eqref{eq:cf_condition_psi}: the sufficient condition holds more often with $\psi=0$, and indeed it was satisfied at $71.0\%$ of the nodes explored by $\mathrm{DDGFact}\text{(CF)}$, against $42.5\%$ for $\mathrm{DDGFact}^{+}\text{(CF)}$. The same holds for all values of $s$: the closed form was applied at $64.5\%$ to $81.4\%$ of the nodes of $\mathrm{DDGFact}\text{(CF)}$, and at $42.5\%$ to $56.9\%$ of the nodes of $\mathrm{DDGFact}^{+}\text{(CF)}$, and, for $s\ge10$, these percentages increase with $s$, so the closed form remains effective as the instances become harder. Pruning with the trace upper bound was effective only for $\mathrm{DDGFact}\text{(CF)}$ with $s=5$, where it pruned $800875$ nodes, $41.3\%$ of the nodes explored, contributing to the speed-up of this strategy. Variables were fixed only for $s=5$, the instance with the smallest gaps; for $s\ge10$, the gaps remained too large to fix any variable.

For the instances with $s\ge10$, all runs stopped at the time limit, and the number of nodes explored measures how fast each strategy processes the nodes. In this case, a large number of nodes is a desirable feature: it shows that the upper bound is computed quickly, and it increases the chance of reducing the gap within the time limit. Our results confirm that, among strategies that compute the same upper bounds, the reduction of the gap during the run, given by Root Gap $-$ Gap, increases with the number of nodes explored. As noted above, a hybrid strategy computes the same bound as the corresponding relaxation at every node. $\mathrm{DDGFact}\text{(CF)}$ explored $19$ to $31$ times more nodes than DDGFact, and $\mathrm{DDGFact}^{+}\text{(CF)}$ explored $9$ to $11$ times more nodes than \ref{DDGFactplus}; correspondingly, their gap reductions were larger by factors of $1.6$ to $1.7$ and $1.3$ to $1.5$, respectively. However, speed does not always compensate for a weaker bound: $\mathrm{DDGFact}\text{(CF)}$ explored $1.2$ to $1.6$ times more nodes than $\mathrm{DDGFact}^{+}\text{(CF)}$, but obtained smaller gap reductions and final gaps larger by $0.34$ to $0.41$. When we compare DDGFact with \ref{DDGFactplus}, the augmentation helps in both respects: \ref{DDGFactplus} explored $1.5$ to $1.9$ times more nodes, in agreement with its smaller runtimes in Table~\ref{tab:runtime_stats_n63_n90}, its gap reductions were $1.6$ to $1.8$ times larger, and its root and final gaps were smaller by $0.28$ to $0.38$ and by $0.35$ to $0.40$, respectively.

\begin{table}[htbp]
\centering
\resizebox{\columnwidth}{!}{%
\begin{tabular}{cc|ccccc|cc}
\hline
$s$ & Relaxation & Root Gap & Gap & Nodes & Fixed & Time & CF & Prunes \\
\hline
\multirow{4}{*}{5} & $\mathrm{DDGFact}$ & 0.988 & 0.199 & 2011329 & 412884 & 7200.2 & -- & -- \\
 & $\mathrm{DDGFact}\text{(CF)}$ & 0.988 & 0.000 & 1940040 & 85765 & 382.5 & 1377541 & 800875 \\
 & $\mathrm{DDGFact}^{+}$ & 0.786 & 0.000 & 824341 & 34566 & 2698.4 & -- & -- \\
 & $\mathrm{DDGFact}^{+}\text{(CF)}$ & 0.786 & 0.000 & 823977 & 32374 & 625.2 & 349908 & 0 \\
\hline
\multirow{4}{*}{10} & $\mathrm{DDGFact}$ & 1.344 & 1.217 & 896049 & 0 & 7200.3 & -- & -- \\
 & $\mathrm{DDGFact}\text{(CF)}$ & 1.344 & 1.125 & 27924407 & 0 & 7200.5 & 18009093 & 0 \\
 & $\mathrm{DDGFact}^{+}$ & 1.069 & 0.863 & 1613345 & 0 & 7200.2 & -- & -- \\
 & $\mathrm{DDGFact}^{+}\text{(CF)}$ & 1.069 & 0.781 & 17410155 & 0 & 7200.2 & 9255824 & 0 \\
\hline
\multirow{4}{*}{15} & $\mathrm{DDGFact}$ & 1.667 & 1.598 & 712673 & 0 & 7200.6 & -- & -- \\
 & $\mathrm{DDGFact}\text{(CF)}$ & 1.667 & 1.557 & 19820701 & 0 & 7201.0 & 14618777 & 0 \\
 & $\mathrm{DDGFact}^{+}$ & 1.348 & 1.231 & 1339689 & 0 & 7200.5 & -- & -- \\
 & $\mathrm{DDGFact}^{+}\text{(CF)}$ & 1.348 & 1.195 & 12655203 & 0 & 7200.5 & 6850056 & 0 \\
\hline
\multirow{4}{*}{20} & $\mathrm{DDGFact}$ & 1.848 & 1.804 & 654419 & 0 & 7202.0 & -- & -- \\
 & $\mathrm{DDGFact}\text{(CF)}$ & 1.848 & 1.779 & 14841519 & 0 & 7201.1 & 11716426 & 0 \\
 & $\mathrm{DDGFact}^{+}$ & 1.495 & 1.417 & 1206997 & 0 & 7200.9 & -- & -- \\
 & $\mathrm{DDGFact}^{+}\text{(CF)}$ & 1.495 & 1.393 & 10922463 & 0 & 7201.0 & 6026841 & 0 \\
\hline
\multirow{4}{*}{25} & $\mathrm{DDGFact}$ & 1.993 & 1.970 & 611659 & 0 & 7201.7 & -- & -- \\
 & $\mathrm{DDGFact}\text{(CF)}$ & 1.993 & 1.953 & 11665719 & 0 & 7201.7 & 9496805 & 0 \\
 & $\mathrm{DDGFact}^{+}$ & 1.609 & 1.567 & 894779 & 0 & 7201.4 & -- & -- \\
 & $\mathrm{DDGFact}^{+}\text{(CF)}$ & 1.609 & 1.547 & 9363831 & 0 & 7201.4 & 5326242 & 0 \\
\hline
\end{tabular}
}
\caption{Branch-and-bound experiments for $n = 63$ and $t = 2$.}
\label{tab:bnb_comparison_n63_t2}
\end{table}

Table~\ref{tab:bnb_comparison_n63_kappa5} reports the results for $\kappa=5$. No instance was solved within the time limit. The sufficient condition for the closed form did not hold, and no node was pruned by the trace upper bound, at any of the nodes explored, so the columns CF and Prunes are omitted. This is expected from \eqref{eq:cf_condition_psi}: with $t=s-5$, the condition requires the sum of the $s$ eigenvalues of $(C-\psi I)[\bar S,\bar S]$ to be at least $s-5$ times the greatest of them. Consequently, the hybrid strategies behave essentially as the corresponding relaxations, with node counts within $4\%$ and the same gaps up to $0.002$, which also shows that the cost of checking the condition is negligible.

The augmentation has a strong effect on these harder instances: \ref{DDGFactplus} reduces the root gap of DDGFact by $0.52$ to $2.30$, and the final gap by $36\%$ to $73\%$, with larger improvements for larger $s$; it also explores $9\%$ to $46\%$ more nodes. Only \ref{DDGFactplus} and $\mathrm{DDGFact}^{+}\text{(CF)}$ fixed variables, and only for $s\in\{40,50\}$, the two instances with the smallest root gaps in Table~\ref{tab:bnb_comparison_n63_kappa5}. A small gap is necessary, but not sufficient, for fixing, as the fixing rule compares the gap with the dual values $\hat\upsilon_i$ and $\hat\nu_i$: for $s=10$, \ref{DDGFactplus} attains a smaller final gap than for $s=40$, but it does not fix any variable.

In summary, the hybrid strategy based on the trace upper bound and the closed-form certificate is very effective for small $t$, where the sufficient condition holds at most nodes, whereas for larger $t$ the best results are obtained with \ref{DDGFactplus}, which clearly improves on DDGFact in all our B\&B experiments.

\begin{table}[htbp]
\centering
\resizebox{\columnwidth}{!}{%
\begin{tabular}{ccc|ccccc}
\hline
$s$ & $t$ & Relaxation & Root Gap & Gap & Nodes & Fixed & Time \\
\hline
\multirow{4}{*}{10} & \multirow{4}{*}{5} & $\mathrm{DDGFact}$ & 2.009 & 1.715 & 916041 & 0 & 7200.2 \\
 &  & $\mathrm{DDGFact}\text{(CF)}$ & 2.009 & 1.715 & 908485 & 0 & 7200.2 \\
 &  & $\mathrm{DDGFact}^{+}$ & 1.489 & 1.090 & 1340749 & 0 & 7200.2 \\
 &  & $\mathrm{DDGFact}^{+}\text{(CF)}$ & 1.489 & 1.091 & 1327565 & 0 & 7200.2 \\
\hline
\multirow{4}{*}{20} & \multirow{4}{*}{15} & $\mathrm{DDGFact}$ & 3.004 & 2.746 & 637049 & 0 & 7201.0 \\
 &  & $\mathrm{DDGFact}\text{(CF)}$ & 3.004 & 2.747 & 633203 & 0 & 7201.0 \\
 &  & $\mathrm{DDGFact}^{+}$ & 1.855 & 1.566 & 864657 & 0 & 7201.0 \\
 &  & $\mathrm{DDGFact}^{+}\text{(CF)}$ & 1.855 & 1.567 & 851127 & 0 & 7201.1 \\
\hline
\multirow{4}{*}{30} & \multirow{4}{*}{25} & $\mathrm{DDGFact}$ & 3.192 & 2.942 & 634287 & 0 & 7202.2 \\
 &  & $\mathrm{DDGFact}\text{(CF)}$ & 3.192 & 2.942 & 644729 & 0 & 7202.3 \\
 &  & $\mathrm{DDGFact}^{+}$ & 1.738 & 1.516 & 693287 & 0 & 7202.2 \\
 &  & $\mathrm{DDGFact}^{+}\text{(CF)}$ & 1.738 & 1.517 & 675987 & 0 & 7202.2 \\
\hline
\multirow{4}{*}{40} & \multirow{4}{*}{35} & $\mathrm{DDGFact}$ & 3.244 & 2.974 & 527521 & 0 & 7204.7 \\
 &  & $\mathrm{DDGFact}\text{(CF)}$ & 3.244 & 2.972 & 547145 & 0 & 7203.7 \\
 &  & $\mathrm{DDGFact}^{+}$ & 1.415 & 1.221 & 634479 & 14354 & 7203.5 \\
 &  & $\mathrm{DDGFact}^{+}\text{(CF)}$ & 1.415 & 1.222 & 623215 & 13966 & 7204.0 \\
\hline
\multirow{4}{*}{50} & \multirow{4}{*}{45} & $\mathrm{DDGFact}$ & 3.485 & 3.058 & 578805 & 0 & 7202.6 \\
 &  & $\mathrm{DDGFact}\text{(CF)}$ & 3.485 & 3.058 & 577199 & 0 & 7202.1 \\
 &  & $\mathrm{DDGFact}^{+}$ & 1.189 & 0.836 & 628361 & 932225 & 7202.2 \\
 &  & $\mathrm{DDGFact}^{+}\text{(CF)}$ & 1.189 & 0.837 & 620689 & 925397 & 7202.9 \\
\hline
\end{tabular}
}
\caption{Branch-and-bound experiments for $n = 63$ and $t = s - \kappa$, with $\kappa = 5$.}
\label{tab:bnb_comparison_n63_kappa5}
\end{table}

\section{Outlook}\label{sec:outlook}

We are currently working on a ``g-scaling'' technique (in the spirit of \citep*{gscale}; also see \citep*{ExtendedGMESP})
for \ref{DGFactplus},
which holds promise for bound improvement.

\backmatter

\bmhead{Acknowledgements}
A.~Machado was supported in part by CNPq grant 131219/2026-0.
M.~Fampa was supported in part by CNPq grant 304340/2026-0.
J.~Lee was supported in part by AFOSR grant FA9550-22-1-0172.

\newcounter{maingame}
\setcounter{maingame}{\value{equation}}

\begin{appendices}

\makeatletter
\renewcommand{\theequation}{\arabic{equation}}
\setcounter{equation}{\value{maingame}}
\makeatother

\section{Simple lemmas and postponed proofs}\label{sec:A1}

\begin{lemma}\label{lem:subproblem_W_univariable}
    Let $\psi \ge 0$. The optimization problem
    $
        \textstyle\sup_{y \geq 0} \{ \log(y + \psi) - wy \}
$
    is unbounded above if $w \leq 0$. If $w > 0$, it admits a unique optimal solution
$
        y^* = \max\left\{0, 1/w - \psi \right\},
$
    and $y^* > 0$ if and only if $0 < w < 1/\psi$.
\end{lemma}

\begin{proof}
    If $w \leq 0$, the objective clearly tends to $+\infty$ as $y \to +\infty$. For $w > 0$, the objective is strictly concave. Setting the derivative of the unconstrained objective to zero yields the unconstrained maximizer $y = 1/w - \psi$. Projecting this solution onto the feasible region $y \geq 0$ yields the unique optimal solution $y^* = \max\{0, 1/w - \psi\}$. Finally, $y^* > 0$ holds if and only if $1/w - \psi > 0$, which simplifies to $w < 1/\psi$ (a condition that holds for every $w>0$ when $\psi=0$).
\end{proof}

\begin{lemma}~\!\!\![KKT verification for Theorem~\ref{thm:subproblem_Theta}]
\label{lem:kkt_verification}
Let $\lambda$, $\psi$, $\hat y$, $\hat\iota$, $\bar y$, $\hat w$ and $\hat\eta$ be as in the statement and proof of Theorem~\ref{thm:subproblem_Theta}. Then $(\hat w,\hat\eta)$ satisfies the KKT conditions of step~2 of that proof, that is, $\hat w\in\mathcal{W}$, $\hat\eta_\ell\ge0$, for $\ell=1,\ldots,k-1$, stationarity and complementary slackness.
\end{lemma}

\begin{proof}
Observe that substituting the definition of $\bar y$ into the expression for  $\hat\eta_k$ yields  $\hat\eta_k=0$, as required by the  boundary condition.

\emph{Primal feasibility.}
As $\hat y_1\ge\cdots\ge\hat y_k$\,, we have
$
\hat w_1\le\cdots\le\hat w_{\hat\iota}\,.
$
Moreover, Lemma~\ref{lem:index_iota} gives
$
\hat y_{\hat\iota}>\bar{y}\ge \hat y_{\hat\iota+1}\,,
$
which implies
$
\hat w_{\hat\iota}
<
1/\bar{y}
=
\hat w_{\hat\iota+1}
=
\cdots
=
\hat w_k\,,
$ where the identities follow from the definition of $\hat w$.
Moreover, $\hat w_1>0$, as $\hat y_1>0$ and $\bar y>0$. Hence $\hat w\in\mathcal{W}$.

\emph{Dual feasibility.}
By definition, $\hat\eta_\ell=0$ for $\ell\le\hat\iota$.
Lemma~\ref{lem:index_iota} implies
$\hat y_j\le \bar y$ for all $j\ge \hat\iota+1$, therefore,
$\hat\eta_\ell\ge0$, for $\ell=\hat\iota+1,\dots,t$.
For $\ell\ge t$, using $\hat\eta_k=0$, we note that
$
\hat\eta_\ell
=
\textstyle\sum_{j=\ell+1}^{k}\hat y_j\,,
$
which is nonnegative because $\hat y\ge0$. Hence $\hat\eta$ is dual feasible.

\emph{Stationarity.}
For $\ell\le\hat\iota$, we have
$\hat\eta_\ell=\hat\eta_{\ell-1}=0$
and
$\hat w_\ell=1/\hat y_\ell$\,, so the stationarity equations hold.
For $\hat\iota<\ell\le t$, we have 
$
\hat\eta_\ell-\hat\eta_{\ell-1}
=
\bar{y} -\hat y_\ell\,,
$
while for $\ell>t$,
$
\hat\eta_\ell-\hat\eta_{\ell-1}
=
-\hat y_\ell\,.
$
Substituting these identities into the stationarity equations verifies them for every $\ell$.

\emph{Complementary slackness.}
If $\ell\le\hat\iota$, then $\hat\eta_\ell=0$. If $\hat\iota<\ell \le k-1$, then
$\hat w_\ell=\hat w_{\ell+1}=1/\bar y$. Hence,
$
\hat\eta_\ell(\hat w_\ell-\hat w_{\ell+1})=0$, for 
$\ell=1,\dots,k-1$.
\end{proof}

\begin{theorem}~\!\!\![Variable fixing based on dual information]
\label{thm:variable_fixing_dgfact_plus}
    Let {\rm LB} denote the objective value of any feasible solution for
    \ref{GMESP}. Let $(\hat\Theta,\hat\upsilon,\hat\nu,\hat\tau)$ be a feasible solution of \ref{DGFactplus}, with objective value $\hat\zeta$. Then every optimal solution $x^*$ of \ref{GMESP} satisfies
    \begin{align*}
        &x_i^*=0,
        \quad
        \forall i\in N
        \ \text{such that}\
        \hat\zeta-{\rm LB}<\hat\upsilon_i\,,\\
        &x_i^*=1,
        \quad
        \forall i\in N
        \ \text{such that}\
        \hat\zeta-{\rm LB}<\hat\nu_i\,.
    \end{align*}
\end{theorem}

\begin{proof}
Let $x^*$ be an optimal solution of \ref{GMESP}, and let $z^*$
denote its objective value. Because ${\rm LB}$ is the value of a feasible solution,
${\rm LB}\le z^*$.

Because $x^*$ is binary and satisfies
$\mathbf e^\top x^*=s$, it is feasible for \ref{GFactplus}. Setting
$W^*=M_\psi(x^*)$ and evaluating the Lagrangian at
$(x^*,W^*)$ gives
\[
\mathcal L(x^*,W^*,\hat\Theta,\hat\upsilon,\hat\nu,\hat\tau)
=
z^*
+
\hat\upsilon^\top x^*
+
\hat\nu^\top(\mathbf e-x^*),
\]
because the equality constraints are satisfied { and, by Theorem~\ref{Gmesp2GFact}, $\sum_{i=1}^t\log(\lambda_i(W^*)+\psi)=z^*$}.
By weak duality,
$
\mathcal L(x^*,W^*,\hat\Theta,\hat\upsilon,\hat\nu,\hat\tau)
\le
\hat\zeta$,
and therefore
\[
\hat\upsilon^\top x^*
+
\hat\nu^\top(\mathbf e-x^*)
\le
\hat\zeta-z^*
\le
\hat\zeta-{\rm LB}.
\]

Suppose that
$\hat\zeta-{\rm LB}<\hat\upsilon_i$\,.
If $x_i^*=1$, then, because
$\hat\upsilon,\hat\nu\ge0$,
\[
\hat\upsilon^\top x^*
+
\hat\nu^\top(\mathbf e-x^*)
\ge
\hat\upsilon_i\,,
\]
contradicting the previous inequality. Hence $x_i^*=0$.

The proof of the second statement is identical.
Indeed, if
$\hat\zeta-{\rm LB}<\hat\nu_i$
and
$x_i^*=0$,
then
$
\hat\upsilon^\top x^*
+
\hat\nu^\top(\mathbf e-x^*)
\ge
\hat\nu_i$\,,
again yielding a contradiction.
Therefore $x_i^*=1$.
\end{proof}

\begin{proof}~\!\!\![of Lemma~\ref{lem:von_neumann}]
    The right-hand inequality in \eqref{eq:von_neumann} is the (unnumbered) theorem of \citep*[p.~265]{Theobald1975}. Applying it to $A$ and $-B$, and using $\lambda_i(-B)=-\lambda_{n-i+1}(B)$, gives the left-hand inequality. For the attainment claim, as $Q$ is orthogonal,
    $\bar B=Q\Diag(\beta)Q^\top$ is symmetric with the same spectrum as
    $\Diag(\beta)$, hence $\lambda(\bar B)=\lambda(B)$. Moreover,
    \begin{equation*}
        A\bullet\bar B
        =
        \tr\!\left(Q\Diag(\lambda(A))Q^\top Q\Diag(\beta)Q^\top\right)
        =
        \tr\!\left(\Diag(\lambda(A))\Diag(\beta)\right)
        =
        \textstyle\sum_{i=1}^{n}\lambda_i(A)\,\beta_i\,. \qedhere
    \end{equation*}
\end{proof}

\begin{proof}~\!\!\![of Lemma~\ref{lem:subproblem_W_multivariable}]
    We consider the optimization problem \eqref{eq:solve_y_ordered_given_z_ordered}. If $w_1 \leq 0$, fix $y_2 = \dots = y_t = 1$ and let $y_1\ge 1$. These values satisfy the ordering constraints, and the objective of \eqref{eq:solve_y_ordered_given_z_ordered} becomes $\log(y_1+\psi)-w_1y_1+\sum_{i=2}^{t}\left(\log(1+\psi)-w_i\right)$, where the sum is a finite constant also when $\psi=0$. As $w_1\le 0$, the objective tends to $+\infty$ as $y_1\to+\infty$; thus, the problem is unbounded.

    Now assume $w_1 > 0$. Because the components of $w$ are sorted in nondecreasing order, $w_i > 0$ for all $i$. Dropping the ordering constraints $y_1 \geq \dots \geq y_t$ yields a separable relaxation. By Lemma~\ref{lem:subproblem_W_univariable}, its unique optimal solution is
    \begin{equation*}
        y_i^* = \max\left\{0, 1/w_i - \psi \right\}, \quad i=1,\dots,t.
    \end{equation*}
    Because $w_1 \leq \dots \leq w_t$, it follows that $y_1^* \geq \dots \geq y_t^*$. Because this relaxed solution satisfies the omitted ordering constraints, it is  optimal for the original problem.
\end{proof}

{
\begin{proof}~\!\!\![of Lemma~\ref{lem:phi_t_subgrad}]
Let $z\in\mathbb{R}^k_+$ with $z_1\ge\cdots\ge z_k$ and $z_t=0$, let $r<t$ be the number of positive components of $z$, and let $\iota$ be the index given by Lemma~\ref{lem:index_iota} for $z$. If $\iota<r$, then $z_{\iota+1}>0$ and $\sum_{j=\iota+1}^k z_j=\sum_{j=\iota+1}^r z_j\le(r-\iota)z_{\iota+1}<(t-\iota)z_{\iota+1}$\,, contradicting Lemma~\ref{lem:index_iota}. Hence $\iota\ge r$, so that $\sum_{j=\iota+1}^k z_j=0$, and, as $t-\iota\ge1$, $\phi_t(z)=-\infty$.

Now apply Theorem~\ref{thm:subproblem_Theta} with $\psi:=0$ and $\lambda:=y$, which is allowed because $y_t>0$. Then $\hat y=y$ and $\hat w=g$, so that $\bar y>0$, $g$ is well defined and positive, and, as $y$ is nonincreasing and $y_{\hat\iota}>\bar y$, $g$ is nondecreasing. Moreover, $\phi_t(y)=\sum_{\ell=1}^t h_0(g_\ell)+y^\top g$. Let $z\in\mathbb{R}^k_+$ with $z_1\ge\cdots\ge z_k$\,. If $z_t=0$, the claimed inequality holds trivially, as $\phi_t(z)=-\infty$. Otherwise, Theorem~\ref{thm:subproblem_Theta} with $\psi:=0$ and $\lambda:=z$ shows that $\phi_t(z)$ is the optimal value of \eqref{probinw} with $\lambda:=z$; as $g$ is feasible for that problem,
\[
\textstyle\phi_t(z)\le\sum_{\ell=1}^t h_0(g_\ell)+z^\top g=\phi_t(y)+g^\top(z-y). \qedhere
\]
\end{proof}
}

{
\begin{proof}~\!\!\![of Lemma~\ref{lem:gamma_t_directional}]
Let $\hat M:=M_\psi(\hat x)=Q\Diag(\lambda)Q^\top$ be as in Proposition~\ref{prop:dgfactplus_associated_solution_fixed}, and write $\hat\lambda:=\lambda=\lambda(\hat M)$ and $\hat y:=\hat\lambda+\psi\mathbb{I}_t$\,. Let $\hat\iota$, $\bar y$ and $\hat w$ be as in Theorem~\ref{thm:subproblem_Theta}, so that $\bar y>0$, $\hat w_\ell=1/\hat y_\ell$ for $\ell=1,\ldots,\hat\iota$, and $\hat w_\ell=1/\bar y$ for $\ell=\hat\iota+1,\ldots,k$. For $\iota\in\{0,\ldots,t-1\}$ and $y\in\mathbb{R}^k$, let
$m_\iota(y):=\frac{1}{t-\iota}\sum_{j=\iota+1}^{k}y_j$ and, wherever $y_1,\ldots,y_\iota>0$ and $m_\iota(y)>0$,
\[
\textstyle G_\iota(y):=\sum_{j=1}^{\iota}\log y_j+(t-\iota)\log m_\iota(y).
\]
By Definition~\ref{def:phi_t}, $\phi_t(y)=G_{\iota(y)}(y)$ for every $y\in\mathbb{R}^k_+$ with $y_1\ge\cdots\ge y_k$, where $\iota(y)$ denotes the index given by Lemma~\ref{lem:index_iota} (with the convention $y_0=+\infty$).

{1. The candidate indices.}
Let $\mathcal{I}$ be the set of $\iota\in\{0,\ldots,t-1\}$ such that $\hat y_\iota\ge m_\iota(\hat y)\ge\hat y_{\iota+1}$ (with $\hat y_0=+\infty$). We show that, for every $\iota\in\mathcal{I}$, $G_\iota$ is continuously differentiable on a neighborhood of $\hat y$, $G_\iota(\hat y)=\phi_t(\hat y)$ and $\nabla G_\iota(\hat y)=\hat w$. Recall that $\sum_{j=\hat\iota+1}^{k}\hat y_j=(t-\hat\iota)\bar y$ and, by Lemma~\ref{lem:index_iota}, $\hat y_j>\bar y$ for $j\le\hat\iota$ and $\hat y_j\le\bar y$ for $j>\hat\iota$. Let $\iota\in\mathcal{I}$.

If $\iota<\hat\iota$, then $m_\iota(\hat y)=\frac{1}{t-\iota}\big(\sum_{j=\iota+1}^{\hat\iota}\hat y_j+(t-\hat\iota)\bar y\big)$ is a weighted average of the numbers $\hat y_{\iota+1},\ldots,\hat y_{\hat\iota}$\,, all at most $\hat y_{\iota+1}$, and of $\bar y<\hat y_{\iota+1}$, the latter with positive weight, as $t-\hat\iota\ge1$. Hence $m_\iota(\hat y)<\hat y_{\iota+1}$, contradicting $\iota\in\mathcal{I}$. Thus $\iota\ge\hat\iota$.

If $\iota>\hat\iota$, then, as $\hat y_j\le\bar y$ for $j>\hat\iota$,
\[
\textstyle m_\iota(\hat y)=\frac{1}{t-\iota}\big((t-\hat\iota)\bar y-\sum_{j=\hat\iota+1}^{\iota}\hat y_j\big)\ge\bar y\ge\hat y_\iota\ge m_\iota(\hat y),
\]
where the last inequality holds because $\iota\in\mathcal{I}$. Hence all these inequalities hold with equality, so that $m_\iota(\hat y)=\bar y$ and $\hat y_j=\bar y$, for $j=\hat\iota+1,\ldots,\iota$. This also holds trivially when $\iota=\hat\iota$.

Therefore, for every $\iota\in\mathcal{I}$, we have $\hat y_j\ge\bar y>0$ for $j\le\iota$, and $m_\iota(\hat y)=\bar y>0$. Hence $G_\iota$ is continuously differentiable on a neighborhood of $\hat y$, with $\partial G_\iota(y)/\partial y_j=1/y_j$ for $j\le\iota$ and $\partial G_\iota(y)/\partial y_j=1/m_\iota(y)$ for $j>\iota$, so that $\nabla G_\iota(\hat y)=\hat w$. Moreover,
$G_\iota(\hat y)=\sum_{j=1}^{\hat\iota}\log\hat y_j+(\iota-\hat\iota)\log\bar y+(t-\iota)\log\bar y=\phi_t(\hat y)$.

{2. First-order expansion of $\phi_t$ at $\hat y$.}
We show that, for every $\varepsilon>0$, there exists a neighborhood $\mathcal{U}$ of $\hat y$ such that
\begin{equation}\label{eq:phi_t_expansion}
|\phi_t(y)-\phi_t(\hat y)-\hat w^\top(y-\hat y)|\le\varepsilon\|y-\hat y\|,
\end{equation}
for all $y\in\mathbb{R}^k_+$ with $y_1\ge\cdots\ge y_k$ and $y\in\mathcal{U}$. Otherwise, there exist $\varepsilon>0$ and a sequence $\{y^r\}$ of such vectors converging to $\hat y$ for which \eqref{eq:phi_t_expansion} fails. As $\iota(y^r)\in\{0,\ldots,t-1\}$, passing to a subsequence we may assume that $\iota(y^r)=\iota$ for all $r$. By Lemma~\ref{lem:index_iota}, $y^r_\iota>m_\iota(y^r)\ge y^r_{\iota+1}$, and letting $r\to\infty$ we obtain $\iota\in\mathcal{I}$. By step~1, $\phi_t(y^r)=G_\iota(y^r)$ for every $r$, $\phi_t(\hat y)=G_\iota(\hat y)$ and $\nabla G_\iota(\hat y)=\hat w$. As $G_\iota$ is differentiable at $\hat y$, \eqref{eq:phi_t_expansion} holds for $y=y^r$ when $r$ is sufficiently large, a contradiction.

{3. Differentiability of $\hat w^\top\lambda(\cdot)$ at $\hat M$.}
The vector $\hat\lambda$ block-refines $\hat w$ in the sense of \citep*{Lewis1996derivatives}, that is, $\hat w_i=\hat w_j$ whenever $\hat\lambda_i=\hat\lambda_j$\,. Indeed, if $i,j\le\hat\iota$, then $\hat w_i=1/(\hat\lambda_i+\psi)=1/(\hat\lambda_j+\psi)=\hat w_j$; if $i,j>\hat\iota$, then $\hat w_i=\hat w_j=1/\bar y$; and $i\le\hat\iota<j$ is impossible, because $\hat\lambda_i\ge\hat\lambda_{\hat\iota}>\hat\lambda_{\hat\iota+1}\ge\hat\lambda_j$\,, where the strict inequality follows from $\hat y_{\hat\iota}>\bar y\ge\hat y_{\hat\iota+1}$ and $\hat\iota+1\le t$. Hence, by \citep*[Lem.~2.2]{Lewis1996derivatives}, in its version for real symmetric matrices (see the end of \citep*[Sec.~1]{Lewis1996derivatives}), the function $X\mapsto\hat w^\top\lambda(X)$ is differentiable at $\Diag(\hat\lambda)$, with gradient $\Diag(\hat w)$. As $\hat w^\top\lambda(X)=\hat w^\top\lambda(Q^\top XQ)$ for all $X\in\mathbb{S}^k$, and $Q^\top\hat MQ=\Diag(\hat\lambda)$, the chain rule shows that $X\mapsto\hat w^\top\lambda(X)$ is differentiable at $\hat M$, with gradient $Q\Diag(\hat w)Q^\top=\hat\Theta$.

{4. Conclusion.}
Let $x\in\mathbb{R}^n_+$, $X:=M_\psi(x)\in\mathbb{S}^k_+$ and $y:=\lambda(X)+\psi\mathbb{I}_t$\,, so that $y\in\mathbb{R}^k_+$, $y_1\ge\cdots\ge y_k$\,, $f(x)=\phi_t(y)$, $f(\hat x)=\phi_t(\hat y)$ and $y-\hat y=\lambda(X)-\hat\lambda$. Moreover, as in the proof of Lemma~\ref{lem:gamma_t_supergradient}, $\hat d^\top(x-\hat x)=\hat\Theta\bullet(X-\hat M)$. Therefore,
\[
f(x)-f(\hat x)-\hat d^\top(x-\hat x)
=
\big[\phi_t(y)-\phi_t(\hat y)-\hat w^\top(y-\hat y)\big]
+
\big[\hat w^\top\lambda(X)-\hat w^\top\hat\lambda-\hat\Theta\bullet(X-\hat M)\big].
\]
Here, $\|\cdot\|$ denotes the Euclidean norm for vectors and the Frobenius norm for matrices. The function $\lambda$ is locally Lipschitz (see \citep*[p.~579]{Lewis1996derivatives}), and $x\mapsto M_\psi(x)$ is linear, so there exists $K>0$ such that $\|\lambda(X)-\hat\lambda\|\le K\|X-\hat M\|\le K^2\|x-\hat x\|$ for all $x$ sufficiently close to $\hat x$ (enlarging $K$ if necessary). Given $\varepsilon>0$, step~2 bounds the first bracket by $\varepsilon K^2\|x-\hat x\|$, and step~3 bounds the second bracket by $\varepsilon\|X-\hat M\|\le\varepsilon K\|x-\hat x\|$, for all $x\in\mathbb{R}^n_+$ sufficiently close to $\hat x$. This proves the limit in the statement. Finally, if $\hat x+\delta v\in\mathbb{R}^n_+$ for all sufficiently small $\delta>0$, taking $x:=\hat x+\delta v$ and letting $\delta\downarrow0$ gives $f'(\hat x;v)=\hat d^\top v$.
\end{proof}
}

\begin{proof}~\!\!\![of Theorem~\ref{thm:dgfactplus_associated_solution_optimal_fixed}]
By Proposition~\ref{prop:dgfactplus_associated_solution_fixed},
$(\hat\Theta,\hat\upsilon,\hat\nu,\hat\tau)$
is feasible for \ref{DGFactpluslc}. Thus, it suffices to show that
its objective value coincides with the optimal value of
\ref{DDGFactpluslc}.

Let $\hat d$, $\pi$, $q$, $P$ and $R$ be defined as in Proposition~\ref{prop:dgfactplus_associated_solution_fixed}.
We first establish that
\begin{equation}
\label{eq:linear_identity}
\hat x^\top\hat d
=
-\hat\upsilon^\top l
+\hat\nu^\top c
+\hat\tau s.
\end{equation}

For
$j=0,\ldots,n$, let
\[
\sigma_j:=
\textstyle\sum_{i=1}^{j}c_{\pi(i)}
+
\sum_{i=j+1}^{n}l_{\pi(i)}\,.
\]
If $q<n$, then
$\sigma_q\le s<\sigma_{q+1}$\,.
As
\[
\sigma_{q+1}-\sigma_q
=
c_{\pi(q+1)}
-
l_{\pi(q+1)}
\in\{0,1\},
\]
we necessarily have
$\sigma_q=s$ and
$l_{\pi(q+1)}=0$.
Consequently,
\[
\textstyle \sum_{i\in P}c_i+\sum_{i\in R}l_i=s.
\]
The same holds trivially when $q=n$.

Using the definitions of
$\hat\upsilon$
and
$\hat\nu$,
we obtain
\[
-\hat\upsilon^\top l
+\hat\nu^\top c
+\hat\tau s
=
\textstyle\sum_{i\in P}\hat d_i c_i
+
\sum_{i\in R}\hat d_i l_i\,.
\]
On the other hand,
let $f(x):=\Gamma_t^+(M_\psi(x);\psi)$, and let $x$ be any feasible solution of \ref{DDGFactpluslc}. As the feasible set of \ref{DDGFactpluslc} is convex and $\hat x$ is optimal, $f(\hat x+\delta(x-\hat x))\le f(\hat x)$ for all $\delta\in(0,1]$, so that, by Lemma~\ref{lem:gamma_t_directional}, $\hat d^\top(x-\hat x)=f'(\hat x;x-\hat x)\le0$. Hence, $\hat x$ is an optimal solution of
\[
\max\{
x^\top\hat d:
l\le x\le c,\;
\mathbf e^\top x=s
\},
\]
whose optimal value is
\[
\textstyle\sum_{i\in P}\hat d_i c_i
+
\sum_{i\in R}\hat d_i l_i\,.
\]
Therefore,
\eqref{eq:linear_identity}
holds. Finally,  Theorem~\ref{thm:subproblem_Theta} yields
\[
\Gamma_t^+(M_{\psi}(\hat x);\psi)
=
\textstyle\sum_{\ell=k-t+1}^{k}
h_\psi(\lambda_\ell(\hat\Theta))
+
M_{\psi}(\hat x)\bullet\hat\Theta.
\]
Moreover,
$
M_{\psi}(\hat x)\bullet\hat\Theta
=
\hat x^\top\hat d.
$
Combining these two  identities with
\eqref{eq:linear_identity}
gives
\[
\Gamma_t^+(M_{\psi}(\hat x);\psi)
=
\textstyle\sum_{\ell=k-t+1}^{k}
h_\psi(\lambda_\ell(\hat\Theta))
-
\hat\upsilon^\top l
+\hat\nu^\top c
+\hat\tau s.
\]
Therefore,  the objective values of
\ref{DDGFactpluslc}
and
\ref{DGFactpluslc}
coincide at
$\hat x$
and
$(\hat\Theta,\hat\upsilon,\hat\nu,\hat\tau)$,
respectively.
As
$(\hat\Theta,\hat\upsilon,\hat\nu,\hat\tau)$
is feasible for 
\ref{DGFactpluslc}, it is optimal.
\end{proof}

\begin{proof}~\!\!\![of Theorem~\ref{thm:dgfactplus_spectral_comparison}]
By Theorem~\ref{thm:dgfactplus_factorization_invariance}, we may assume
$
F=(C-\psi I)^{1/2}$.
Let
\[
\alpha:=\textstyle\frac ts,
\qquad
\hat\Theta:=\alpha C^{-1},
\qquad
\hat\tau:=\alpha,
\qquad
\hat\nu:=0,
\qquad
\hat\upsilon
:=
\alpha\mathbf e-
\diag(F\hat\Theta F^\top).
\]

We first verify that
$(\hat\Theta,\hat\upsilon,\hat\nu,\hat\tau)$
is feasible for \ref{DGFactplus}.
As
\[
F\hat\Theta F^\top
=
\alpha(C-\psi I)^{1/2}C^{-1}(C-\psi I)^{1/2},
\]
its eigenvalues are
\[
\textstyle\alpha\!\left(1-\frac{\psi}{\lambda_\ell(C)}\right),
\qquad
\ell=1,\ldots,n.
\]
Because
$0<\psi\le\lambda_{\min}(C)$,
these eigenvalues lie in $[0,\alpha]$, implying
$
F\hat\Theta F^\top
\in \mathbb{S}^n_+$ 
and 
$\alpha I_n - F\hat\Theta F^\top \in \mathbb{S}^n_+$\,.
Therefore
$
\diag(F\hat\Theta F^\top)\le\alpha\mathbf e$,
so that
$\hat\upsilon\ge0$ and
$
\diag(F\hat\Theta F^\top)
+\hat\upsilon
-\hat\nu
-\hat\tau\mathbf e
=0.$
Therefore
$(\hat\Theta,\hat\upsilon,\hat\nu,\hat\tau)$
is feasible.

It remains to evaluate its objective value.
The eigenvalues of
$\hat\Theta=\alpha C^{-1}$ are
\[
\textstyle\frac{\alpha}{\lambda_n(C)}
\ge\cdots\ge
\frac{\alpha}{\lambda_1(C)}{.}
\]
As
$\alpha\le1$ and
$\psi\le\lambda_\ell(C)$,
\[
\textstyle\frac{\alpha}{\lambda_\ell(C)}
\le
\frac1\psi,
\qquad
\ell=1,\ldots,t{,}
\]
and therefore
\[
\textstyle h_\psi\!\left(
\frac{\alpha}{\lambda_\ell(C)}
\right)
=
-\log\!\left(
\frac{\alpha}{\lambda_\ell(C)}
\right)
+
\frac{\alpha\psi}{\lambda_\ell(C)}
-
1.
\]
Thus,
\[
\begin{aligned}
&
\textstyle\sum_{\ell=k-t+1}^{k}
h_\psi(\lambda_\ell(\hat\Theta))
+\hat\nu^\top\mathbf e
+\hat\tau s
=
\textstyle\sum_{\ell=1}^{t}
\left(
-\log\!\left(
\frac{\alpha}{\lambda_\ell(C)}
\right)
+
\frac{\alpha\psi}{\lambda_\ell(C)}
-
1
\right)
+t
\\
&\qquad =
\textstyle\sum_{\ell=1}^{t}
\log\lambda_\ell(C)
-
t\log\alpha
+
\alpha\psi
\textstyle \sum_{\ell=1}^{t}
\frac1{\lambda_\ell(C)}.
\end{aligned}
\]
As
$\alpha=t/s$,
this equals
\[
\zspectralthing(C,t)
+
t\log(s/t)
+
\textstyle\frac{t\psi}{s}
\sum_{\ell=1}^{t}
\frac1{\lambda_\ell(C)}.
\]
Because \ref{DGFactplus} is a minimization problem, its optimal value
is no greater than the objective value of any feasible solution,
yielding the desired inequality.
\end{proof}

\end{appendices}

\bibliography{Bib} 

\end{document}